\documentclass[11pt]{article}
\usepackage{graphicx}
\usepackage{subcaption}
\usepackage{amsmath}
\usepackage{amssymb}
\usepackage{amsthm}
\usepackage{bbm}
\usepackage{comment}
\newtheorem{theorem}{Theorem}[section]
\newtheorem{remark}{Remark}[section]
\newtheorem{lemma}[theorem]{Lemma}

\usepackage{scalerel,stackengine}
\usepackage[ruled,vlined,noresetcount]{algorithm2e}
\usepackage{algorithmicx}
\usepackage{algpseudocode}
\usepackage{float}
\usepackage{ulem}
\usepackage{cancel}
\usepackage[title]{appendix}
\stackMath

\newcommand\reallywidehat[1]{%
\savestack{\tmpbox}{\stretchto{%
  \scaleto{%
    \scalerel*[\widthof{\ensuremath{#1}}]{\kern.1pt\mathchar"0362\kern.1pt}%
    {\rule{0ex}{\textheight}}%WIDTH-LIMITED CIRCUMFLEX
  }{\textheight}% 
}{2.4ex}}%
\stackon[-6.9pt]{#1}{\tmpbox}%
}
\newcommand{\angflux}{f}

\newcommand{\dt}{\Delta t}

\newcommand{\dx}{h}

\newcommand{\velangle}{v}

\newcommand{\vweight}{\omega}

\newcommand{\spaceset}{\Omega_x}

\newcommand{\velangleset}{\Omega_v}

\newcommand{\vareps}{\varepsilon}

\newcommand{\scat}{\sigma_s}

\newcommand{\absorp}{\sigma_a}

\newcommand{\scatlower}{\sigma_{s,m}}

\newcommand{\absorpupper}{\sigma_{a,M}}

\newcommand{\source}{G}

\newcommand{\lvg}{\langle}
\newcommand{\rvg}{\rangle}

\newcommand{\I}{\mathbf{I}}

\newcommand{\V}{\Pi}

\newcommand{\macro}{\rho}

\newcommand{\micro}{g}

\newcommand{\testmacro}{\phi}

\newcommand{\testmicro}{\psi}

\newcommand{\imag}{\mathcal{I}}

\newcommand{\iden}{J}

\newcommand{\IMEXBDForder}{s}

\newcommand{\pd}{r}

\newcommand{\Dupwindop}{\mathcal{D}_h^{up}}

\newcommand{\Dplusop}{\mathcal{D}_h^+}

\newcommand{\Dminusop}{\mathcal{D}_h^-}

\newcommand{\pdspacebasis}{e}

\newcommand{\M}{M}

\newcommand{\Dplus}{D^+}

\newcommand{\Dminus}{D^-}

\newcommand{\scatmat}{\Sigma_s}

\newcommand{\absorpmat}{\Sigma_a}

\newcommand{\IMEXBDFexpop}{\mathcal{F}}

\newcommand{\IMEXBDFimpop}{\mathcal{G}}

\newcommand{\IMEXBDFprevcoef}{a}

\newcommand{\IMEXBDFexpcoef}{b}

\newcommand{\IMEXBDFimpcoef}{c}

\newcommand{\pdspace}{U_{\dx}}

\newcommand{\Legendre}{\phi}

\newcommand{\Mhat}{\widehat{M}}

\newcommand{\Mhatblock}{\mathbf{M}}

\newcommand{\Dplushat}{\widehat{D}^+}

\newcommand{\Dplushatblock}{\mathbf{{D}^+}}

\newcommand{\Dminushat}{\widehat{D}^-}

\newcommand{\Dminushatblock}{\mathbf{{D}^-}}

\newcommand{\Uhat}{\widehat{U}}

\newcommand{\Phat}{\widehat{P}}

\newcommand{\indicator}{\mathbbm{1}}

\newcommand{\mS}{\mathcal{S}}
\newcommand{\mH}{\mathcal{H}}
\newcommand{\mO}{\mathcal{O}}
\newcommand{\mP}{\mathcal{P}}
\newcommand{\mI}{\mathcal{I}}
\newcommand{\by}{\bf{y}}
\newcommand{\bz}{\bf{z}}
\newcommand{\eps}{\varepsilon}

\newcommand{\UPhatblock}{\mathbf{D}^{up}}

\usepackage{bbding}
\usepackage{pifont}
\usepackage{xcolor}

\newcommand\revc[1]{\textcolor{black}{#1}}

\numberwithin{figure}{section}
\numberwithin{equation}{section}
\numberwithin{table}{section}
\usepackage{multirow}
\usepackage{pdflscape}
\usepackage{diagbox}

\usepackage{color}
\usepackage{pdfpages}
\usepackage{url}
\usepackage{chngcntr}
\counterwithout{algocf}{section}
\counterwithin*{algocf}{subsection}

\usepackage{listings}
\graphicspath{{figures/fourier_analysis/strategy_1}{figures/fourier_analysis/strategy_2}{figures/fourier_analysis/strategy_3}}

\begin{document}
\baselineskip=1.2pc
\title{Multi-scale kinetic simulation: asymptotic preserving IMEX-BDF-DG schemes with three implicit-explicit partitionings}
\author{
Kimberly Matsuda\thanks{Department of Mathematical Sciences, Rensselaer Polytechnic Institute, Troy, NY 12180, U.S.A. Email: 
{\tt  kimberlyamatsuda@gmail.com}}
\and 
Fengyan Li\thanks{Department of Mathematical Sciences, Rensselaer Polytechnic Institute, Troy, NY 12180, U.S.A. Email: 
{\tt lif@rpi.edu}.}
}

\maketitle

\abstract{Kinetic transport models are mesoscopic mathematical descriptions of the transport of particles as well as their interactions with the background media or among themselves, and they have wide applications in many areas of mathematical physics such as nuclear and biomedical engineering, rarefied gas dynamics, {and} plasma physics. They are often multi-scale, with different characteristics (e.g. hyperbolic, diffusive) depending on the  material properties. As our continuing effort in a sequence of works  to design and analyze numerical methods for  accurate and robust simulation of the multi-scale kinetic transport models, in this work,  we consider a linear kinetic transport model, a simplified radiative transfer equation, in a diffusive scaling, and propose and analyze three families of asymptotic preserving (AP) methods. Numerical methods with the AP property, that is to preserve the asymptotic behavior of the models at the discrete level on under-resolved meshes, can work uniformly well to simulate multi-scale models across a wide range of scales.  The proposed methods start from the micro-macro decomposition of the model, and  involve discontinuous Galerkin methods in space, the discrete ordinates method (i.e. $S_N$ method) in velocity, and implicit-explicit (IMEX) BDF methods in time, with three different implicit-explicit partitionings. A systematic study, both analytically and computationally, is presented regarding their difference  in stability, accuracy, computational complexity and AP property. These methods, with multi-step time integrators, are also compared  in terms of their accuracy and efficiency with the ones that only differ in using certain IMEX Runge-Kutta (RK) methods in time. Together with our previous developments in \cite{jang2015high,jang2014analysis,xiong2015high,peng2020stability,peng2021stability,peng2021asymptotic} using IMEX-RK methods in time, the present work further contributes to  high order discontinuous Galerkin AP methods for multi-scale kinetic simulation, especially by utilizing the structure of the micro-macro decomposition of the models.  
}

%------------------------------------
% Section: Introduction
%--------------------------------------

\section{Introduction} 
\label{sec:introduction}

Kinetic transport models (e.g. radiative transfer equation, Boltzmann equation, Vlasov-type models)
are mesoscopic mathematical descriptions of the transport of particles such as photons, molecules, electrons as well as their
interactions with the background media or among themselves, and they have wide applications in many areas of
mathematical physics such as nuclear and biomedical engineering, rarefied gas dynamics, {and} plasma physics. Like many differential equation based mathematical models, numerical simulations are  still one primary approach to understand the solutions hence the  underlying physics.   Though recent years have seen quite  a lot of progress in developing computational algorithms for kinetic transport models, it continues to be an active subject to further enrich and advance the numerical methods, both algorithmically and analytically, for solving these models.   Numerical challenges in the kinetic simulation can arise due to the intrinsic high dimensionality (with the unknown probability distribution functions defined on the phase space), nonlinear or nonlocal nature, intrinsic  structures (e.g. conservation, positivity), and complex multi-scale property.

This paper presents our recent development in a sequence of works to design and analyze high order methods that perform uniformly well in the {\it multi-scale} kinetic simulation. %
We consider a linear kinetic transport equation under a diffusive scaling, 
\begin{equation} \label{eq:kinetic transport equation}
\vareps \partial_t \angflux + \velangle \partial_x \angflux = \frac{\scat}{\vareps} (\lvg \angflux \rvg - \angflux) - \vareps \absorp \angflux +  \vareps \source,
\end{equation}
with  initial and suitable boundary conditions (e.g. periodic, Dirichlet inflow boundary conditions). This is a simplified model of the linear radiative transfer equation with one energy group and  isotropic scattering, and it has the key characteristics when we come to the multi-scale aspect of the full physical model.
In this integro-differential equation, $\angflux = \angflux(x,\velangle,t)$ is the probability distribution function of particles (or the angular flux in the setting of radiative transfer), and it depends on the spatial variable $x \in \spaceset$, the velocity variable  $\velangle \in \velangleset$, and  time $t$. The operator on the left models the free streaming of particles, and the ones on the right model the interaction of particles with the background medium, through the scattering and absorption processes, along with the {external source $\source=\source(x)$.} Related, $\scat = \scat(x) \geq 0$, $\absorp = \absorp(x) \geq 0$, and they are the scattering and absorption cross sections, respectively. In addition, the total cross section  $\scat(x)+\absorp(x)$ is positive on $\spaceset$, \revc{and we only consider the non-degenerate case with  $\scat(x)$ being a nonzero function.}   In the scattering term, $\lvg \angflux \rvg = \int_{\velangleset} \angflux d \nu$  and it gives the macroscopic density (or the scalar flux in the setting of radiative transfer), where $\nu$ is some measure of the velocity space, satisfying {$\int_{\velangleset} 1 d \nu = 1$.}
With the diffusive scaling, we consider a system over a long time period under the assumption  that both the absorption and the source are relatively weak  and comparable in size in the case of strong scattering \cite{larsen1987asymptotic}. The dimensionless parameter $\vareps > 0$ is the Knudsen number, defined as the ratio of the mean free path of particles and the characteristic length of the system.  Though the methodologies in this paper are applicable to more general cases,  two specific models are considered here in our numerical  studies: (i) {the} one-group transport equation in slab geometry, with   $\velangleset = [-1,1]$ and $\lvg \angflux \rvg = \frac{1}{2} \int_{-1}^{1} \angflux(x,\velangle,t) d\velangle$ (associated with the standard Lebesgue measure); (ii) the telegraph equation (also called Goldstein-Taylor equation), with   $\velangleset = \{-1, 1\}$ and $\lvg \angflux \rvg = \frac{1}{2} \big( \angflux(x,-1,t)+\angflux(x,1,t)\big)$.

The model \eqref{eq:kinetic transport equation} is multi-scale in nature. When the scattering is relatively weak (e.g. $\vareps=O(1)$),  the model is in the transport regime and it is more hyperbolic, while with  relatively strong scattering (e.g. $\vareps\ll 1$), the model is more parabolic/diffusive. More specifically, when $\vareps\ll 1$ and $\scat>0$, the density $\macro=\lvg \angflux \rvg$ satisfies 
\begin{equation} \label{eq:diff-limit:0}
\partial_t \macro = \lvg \velangle^2 \rvg \partial_x \big(\scat^{-1}
\partial_x \macro\big) - \absorp \macro +\source + O(\vareps),
\end{equation}
also see Section \ref{sec:micro-macro decomposition}. In practice, the model can be hyperbolic or diffusive in different subregions of the physical domain, due to the spatially dependent material properties  hence the different magnitudes of the effective  $\vareps$.  Standard explicit schemes  applied to \eqref{eq:kinetic transport equation} will require the time step condition  $\dt = O(\vareps \dx)$, where $\dx$ is the spatial mesh size, and this is a stringent 
condition in the diffusive regime with $\vareps \ll 1$ and when the model is stiff.
On the other hand, while fully implicit schemes can allow larger time step sizes and can be efficient,  they may or may not  capture faithfully the correct diffusive limit as $\vareps \rightarrow 0$. {One example to illustrate this is given in Figure \ref{fig:non-AP}, with the scheme defined and analyzed in Appendix \ref{sec:appendA}. }

\begin{figure}[ht]
\centering
{
\begin{subfigure}{0.42\textwidth}
\includegraphics[width=\linewidth]{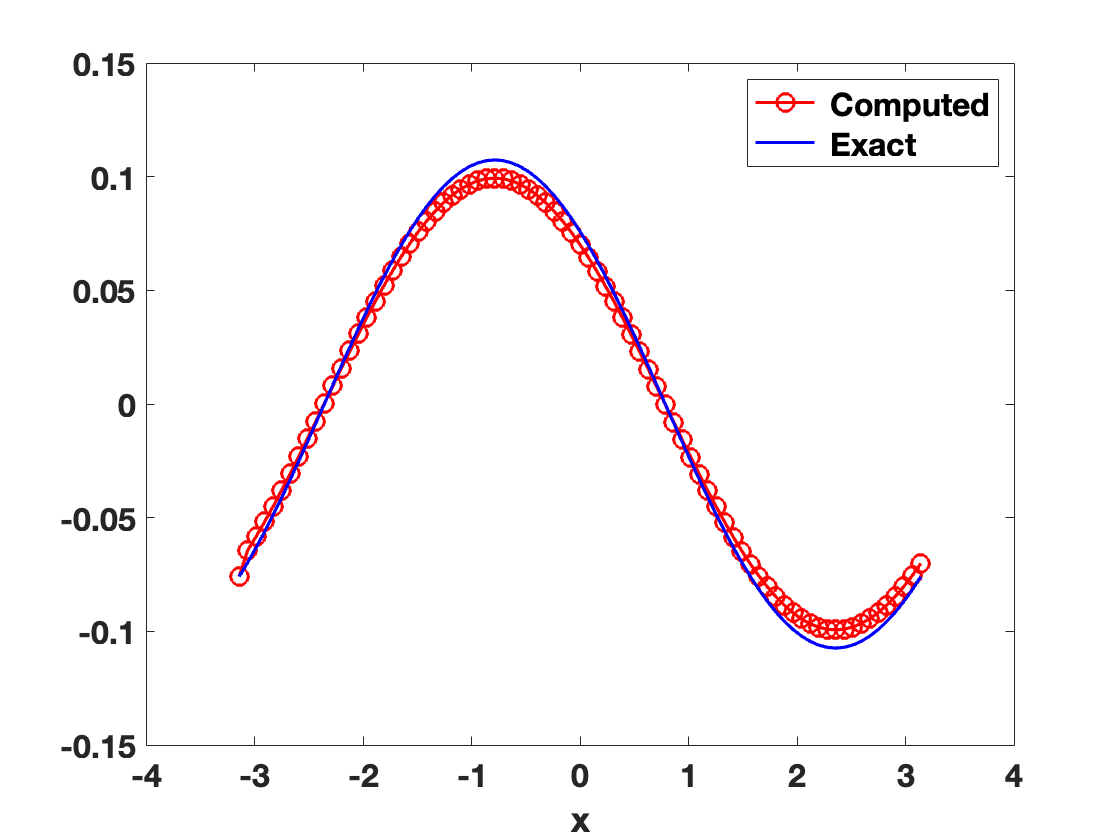}
\caption{$\vareps = 0.5$}
\end{subfigure}
\begin{subfigure}{0.42\textwidth}
\includegraphics[width=\linewidth]{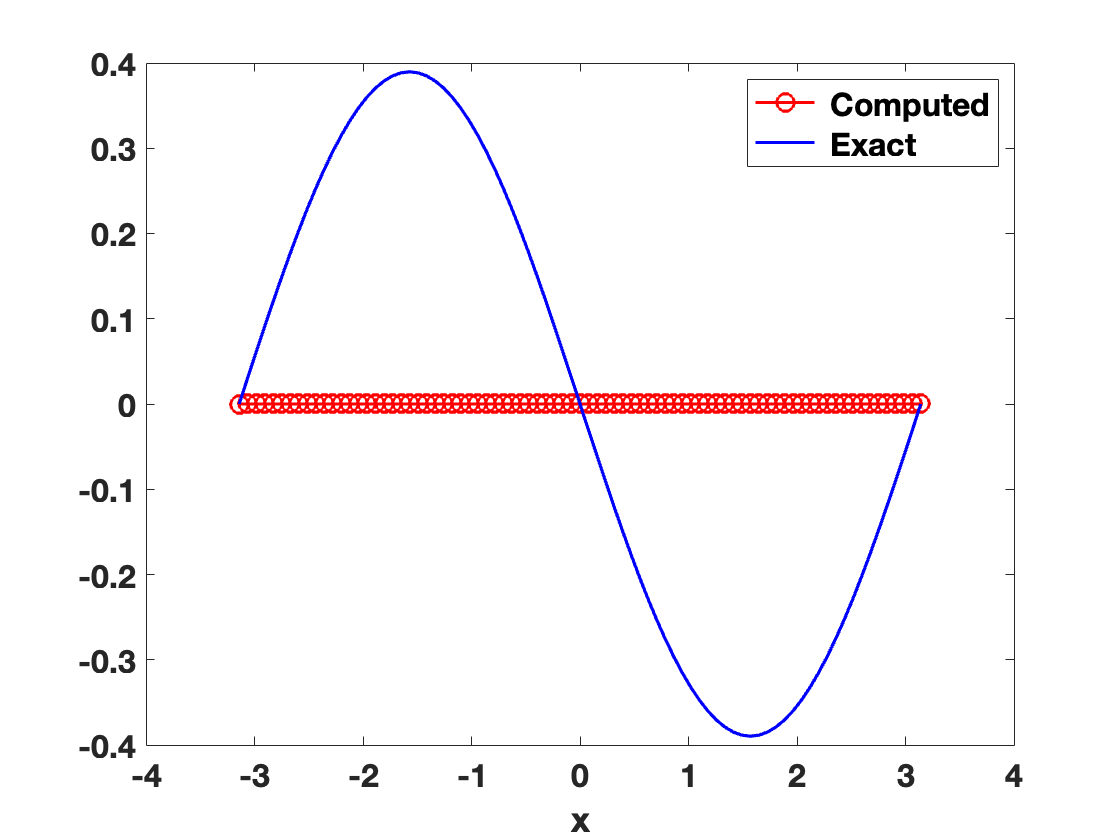}
\caption{$\vareps = 10^{-6}$}
\end{subfigure}
}
\caption{Telegraph equation: $f(x, \velangle=1, T=1)$ on $[-\pi, \pi]$ computed by the first order upwind finite difference (also $P^0$ upwind discontinuous Galerkin) method in space  with the backward Euler method in time for \eqref{eq:kinetic transport equation}. Here $\scat=1, \absorp=0, \source=0$ and the exact solution $f(x, \velangle, t)=\frac{1}{r}e^{rt}\sin(x) + \eps \velangle e^{rt}\cos(x)$, $r= \frac{-2}{1+\sqrt{1-4\epsilon^2}}$,   $\dx=\dt=\pi/40$. The discrepancy in the computed and exact solutions for $\vareps=10^{-6}$ evidences that this fully implicit method is not asymptotic preserving (AP).  }
\label{fig:non-AP}
\end{figure}

To efficiently and reliably simulate the multi-scale kinetic transport equations, it is important that the numerical methods work uniformly well for different regimes. Interested readers can refer to the introduction of \cite{adams2001discontinuous} for a well-argued justification.   One established framework to achieve this is to design numerical methods that are asymptotic preserving (AP).  A numerical scheme is AP for \eqref{eq:kinetic transport equation} if, as $\vareps\rightarrow 0$, it becomes a consistent and stable numerical discretization of the limiting equation \eqref{eq:diff-limit:0}. Though it is often involved to  show rigorously,  such methods provide a pathway to achieve uniformly good performance when the model is multi-scale. 

With a series of works in  \cite{jang2015high,jang2014analysis, peng2020stability, peng2021stability,peng2021asymptotic}, we design and analyze several families of high-order AP methods for the kinetic transport equation in \eqref{eq:kinetic transport equation}. What is common among these methods is that they all start with a reformulation of the model, namely, the micro-macro decomposition in \eqref{eq:micro-macro}, apply (local) discontinuous Galerkin (DG) methods in the physical space with suitably chosen numerical fluxes, the discrete ordinates method (i.e. $S_N$ method \cite{lewis1984computational}) in velocity, and implicit-explicit Runge-Kutta (IMEX-RK) methods  (of certain types \cite{boscarino2013implicit}) in time. These methods differ in how they achieve the AP property, their computational complexity, and the type of stability attained in different regimes. For example, with the limiting schemes as $\vareps\rightarrow 0$ in \cite{jang2015high,jang2014analysis} being explicit, the methods there require a parabolic time step condition $\dt=O(\dx^2)$ in the diffusive regime for stability, and the computational complexity is comparable to some fully explicit methods for \eqref{eq:diff-limit:0}.
The methods in \cite{peng2020stability,peng2021stability,peng2021asymptotic} on the other hand were developed to be unconditionally stable in the diffusive regime, and this is achieved by following two  different ideas, with an additional yet reasonable cost of  inverting a discrete diffusive operator per stage of each RK step: the methods in \cite{peng2020stability,peng2021stability} rely on a second reformulation by adding and subtracting to \eqref{eq:micro-macro macro} a diffusive operator $\lvg \velangle^2 \rvg \partial_x \big(\scat^{-1}(x)\partial_x \macro\big)$ (motivated by the limiting equation  \eqref{eq:diff-limit:0}), with one diffusive operator treated explicitly and one treated implicitly, while \cite{peng2021asymptotic} utilizes {an implicit-explicit (IMEX) partitioning (namely, a strategy to tell which term is treated implicitly or explicitly)} different from that in \cite{jang2015high, jang2014analysis}.  Energy{-}based stability is established when the methods are of first order accuracy in time, and Fourier-based stability analysis is carried out for the methods up to  third order accuracy. DG methods are an ideal candidate for spatial discretizations, with their many attractive properties, especially since they suit well to discretize various differential operators,  including the convective  and diffusive operators in our model.

Note that RK methods are multi-stage time integrators, involving multiple evaluations or inversions of the discrete spatial operator per time step. {In addition, RK methods may suffer from order reduction for  multi-scale problems.}  Motivated by seeking methods with better cost efficiency {and more uniform accuracy} for multi-scale kinetic simulations,  in this work, we further our endeavor in designing high{-}order AP methods for \eqref{eq:kinetic transport equation} by applying {\bf IMEX linear multi-step methods} in time, specifically IMEX-BDF methods in time \cite{ascher1995implicit,wang2008variable}, that involve only one evaluation or inversion of the discrete spatial operator per time step. While {the DG discretization similar as in \cite{jang2015high, jang2014analysis,peng2021asymptotic} is used in space here}  along with the discrete ordinates method \cite{lewis1984computational}  in velocity, we will consider {\bf three} families of methods that differ {in their adopted IMEX partitionings (\cite{lemou2008new, jang2015high, peng2021asymptotic}), also referred to as (IMEX) Strategy $k$ ($k=1, 2, 3$) in this paper.} We will systematically study how the methods differ in their stability, accuracy, computational complexity, and the AP property. In addition, we compare the proposed methods with the ones that only differ in using IMEX-RK schemes in time,  in terms of their accuracy and cost efficiency.
With the multi-step methods,  numerical initialization approaches are also examined. For the proposed methods that require an inversion of the discrete diffusive operator (i.e. {with IMEX Strategy 3, defined in Sections \ref{sec:temporal discretization} and \ref{sec:fully discrete scheme}}),  the condition number is analyzed for the associated  linear system.  Such analysis will inform us about the computational complexity of the methods applied to different regimes, and it is also relevant for understanding the methods with IMEX-RK methods in time from  \cite{peng2020stability, peng2021asymptotic}.  

With the present work together with \cite{jang2015high,jang2014analysis,peng2020stability,peng2021stability,peng2021asymptotic}, our contribution can be summarized as a systematic design and study of   high{-}order AP methods based on DG spatial discretizations for the kinetic transport model \eqref{eq:kinetic transport equation} based on its micro-macro reformulation and IMEX temporal discretizations. We choose to work with the micro-macro reformulation  to address the multi-scale aspect of the model as it directly reveals the role and contribution of each term to the diffusive or transport effects in different regimes and their different levels of stiffness. Alternatively, one can  work with the model in its original form \eqref{eq:kinetic transport equation}, as in the pioneering work of upwind DG methods for stationary radiative transport equations \cite{larsen1983numerical,adams2001discontinuous}. Readers can also refer to the numerical analysis in \cite{guermond2010asymptotic} especially to understand the AP property of the methods and how it is related to the choices of discrete spaces and numerical fluxes. From the implementation point of view, one main difference between the methods based on the micro-macro decomposition and those directly based on \eqref{eq:kinetic transport equation}  comes from the linear solvers. Using our methods as an example, they are either essentially explicit (e.g. those with IMEX Strategy 1 and 2 as defined in  Section  \ref{sec:fully discrete scheme}), or need to invert a discrete diffusive operator (e.g. methods here with   Strategy 3, similar to solving \eqref{eq:diff-limit:0} using the backward Euler method), hence one can utilize classical linear solvers (e.g. Krylov-subspace based iterative solvers, multi-grid solvers or preconditioning techniques) that are available for elliptic/parabolic operators. The methods in \cite{adams2001discontinuous,guermond2010asymptotic} (or the associated fully implicit methods for time{-}dependent models) are mainly solved by source iterations with synthetic accelerations (i.e. pre-conditioning), coupled with transport sweeps along characteristics of the model \cite{adams2002fast}. It will be a meaningful task to perform  a systematic numerical comparison between the two types of AP methods above.

There has been a long history to design and analyze deterministic methods for the robust simulation of multi-scale transport equations. 
These can be methods that are asymptotic preserving \cite{golse1999convergence,jin2010asymptotic,pareschi2011efficient},  or  based on domain decomposition methodologies \cite{bourgat1994coupling,degond2005smooth,filbet2015hierarchy}. 
The methods can involve various spatial discretizations such as finite difference or finite volume methods, and apply in time  splitting methods \cite{jin2010asymptotic}, implicit-explicit Runge-Kutta methods \cite{boscarino2013implicit, jang2015high}, or linear multi-step  methods \cite{ascher1995implicit,dimarco2017implicit}.  They are often based on reformulations of the models, such as   micro-macro decomposition \cite{lemou2008new}, even-odd parity \cite{lewis1984computational,jin2000uniformly}, or adding/subtracting a diffusive operator 
\cite{boscarino2013implicit, peng2020stability}.  {Our proposed methods are stable as $\vareps$ varies from $0$ to $O(1)$, with the type of time step conditions one would expect in both the transport and the diffusive regimes (see e.g., Remark \ref{rem:stab:energy:1}). In  the intermediate regime, our methods require $\dt=O(\vareps\dx)$,} and such time step condition was improved  in  \cite{zhang2023asymptotic} by using a characteristic tracking technique. It is worth noting that there are quite some activities in recent years to develop AP methods that also address the high dimensionality of various kinetic transport models,  e.g., tensor-based low-rank methods  
\cite{einkemmer2021asymptotic, guo2022low}, reduced order models \cite{tencer2016reduced, peng2024micro} and (probabilistic) AP Monte Carlo methods \cite{dimarco2018asymptotic, crestetto2019asymptotically}.

The rest of the paper is organized as follows. In Section \ref{sec:numerical scheme}, the proposed three families of numerical methods are formulated. 
In Section \ref{sec:stability}, stability analysis is established, including  energy-based analysis in Section \ref{sec:stab:energy} for the first order methods and Fourier-based analysis in Section \ref{sec:fourier} for methods up to  third order accuracy, with the difference among the methods highlighted. {Numerical validation and comparison are further performed in Section \ref{sec:stab:numV}.}
Section \ref{sec:alg} presents the algebraic form of the methods. {The condition number of the matrix, that needs to be inverted in the methods with IMEX Strategy 3,  will be estimated regarding its dependence on the model parameters $\vareps, \scat, \absorp$ and discretization parameters $h, \dt$.} Two initialization approaches are discussed in Section \ref{sec:init}, and a formal analysis  of the AP property is presented in Section \ref{sec:AP}, once again 
highlighting the differences among the proposed methods. In Section \ref{sec:num}, a collection of numerical experiments is reported, to showcase the accuracy, AP property, and initialization approaches of the proposed methods, 
along with their cost efficiency in comparison with methods using RK time integrators, through {smooth examples with constant or spatially varying material properties,} a two-material problem, 
a Riemann problem, and an example with {a non-well-prepared initial condition.} We conclude in Section \ref{sec:conclusions}. {For better readability, proofs of all theorems and one lemma are provided in Appendix.}

%----------------------------------------------------------------------------------------------------------------
% Section: Numerical scheme
%----------------------------------------------------------------------------------------------------------------

\section{Proposed numerical schemes}
\label{sec:numerical scheme}

In this section, we will formulate three families of numerical schemes, all based on the micro-macro decomposition of the model \eqref{eq:kinetic transport equation}. The methods mainly differ in their temporal discretizations, especially regarding their adopted IMEX {partitionings}.

%----------------------------------------------------------------------------------------------------------------
% Section: Micro-macro decomposition
%----------------------------------------------------------------------------------------------------------------

\subsection{Micro-macro decomposition}
\label{sec:micro-macro decomposition}

Following \cite{liu2004boltzmann, lemou2008new}, we start with reformulating \eqref{eq:kinetic transport equation} into its micro-macro decomposition.
To this end, we define an orthogonal\footnote{This is with respect to the inner product $\lvg f_1f_2\rvg$  of $f_1, f_2$. } projection operator $\V: \angflux \mapsto \V \angflux  = \lvg \angflux \rvg$, and decompose $\angflux$ into $\angflux = \macro + \vareps \micro$, with the macroscopic part $\macro=\macro(x,t)=\V \angflux$ and the microscopic part $\micro = \micro(x,\velangle,t) 
= \frac{1}{\vareps} (\I - \V) \angflux$. Here $\I$ is the identity operator, {and it is easy to check $\langle \micro\rangle=0$.} By applying $\V$ and $\I - \V$ to \eqref{eq:kinetic transport equation}, our model turns to 
\begin{subequations} \label{eq:micro-macro}
\begin{align}
\partial_t \macro + \partial_x \lvg \velangle \micro \rvg & = -\absorp \macro + G, \label{eq:micro-macro macro} \\
 \partial_t \micro +\frac{1}{\vareps} (\I - \V) (\velangle \partial_x \micro) + \frac{1}{\vareps^2}\velangle \partial_x \macro &= -\frac{1}{\vareps^2}\scat \micro - \absorp \micro. 
\label{eq:micro-macro micro}
\end{align}
\end{subequations}
Formally as $\vareps \rightarrow 0$ and under the assumption $\scat > 0$,  \eqref{eq:micro-macro micro} gives the {\it local equilibrium}  
\begin{equation}
    \scat \micro  = -\velangle \partial_x \macro,
    \label{eq:local-equi}
\end{equation} and this leads to the limiting diffusion equation,
\begin{equation} \label{eq:diff-limit}
\partial_t \macro = -\partial_x \lvg \velangle \micro \rvg - \absorp \macro+ G = \lvg \velangle^2 \rvg \partial_x (\scat^{-1}{\partial_x \macro})
- \absorp \macro +G.
\end{equation}
The multi-scale nature of the model can be seen here. When the scattering is relatively weak (e.g. $\vareps=O(1)$),  the model is more hyperbolic, while with relatively strong scattering (e.g. $\vareps\ll 1$), the model is more parabolic/diffusive. In practice, the model can be hyperbolic or diffusive in different subregions of the physical domain, due to the spatially dependent material properties  hence the different magnitudes of the effective  $\vareps$.

%----------------------------------------------------------------------------------------------------------------
% Subsection: Temporal discretization
%----------------------------------------------------------------------------------------------------------------

\subsection{Temporal discretizations} 
\label{sec:temporal discretization}

To address the stiffness of the model when $\vareps\ll 1$ and to capture the correct asymptotic diffusion limit as $\vareps\rightarrow 0$, while achieving reasonable computational efficiency, in time we  apply implicit-explicit (IMEX) BDF methods \cite{ascher1995implicit}. These are linear multi-step methods, designed to numerically solve an initial value problem 
of the {additively partitioned form}, namely
\begin{equation} \label{eq:IMEX-BDF ode}
\partial_t u = \IMEXBDFexpop(u) + \IMEXBDFimpop(u),
\end{equation}
where $\IMEXBDFexpop(u)$ is a non-stiff operator and treated explicitly, and 
$\IMEXBDFimpop(u)$ is  a stiff operator and treated implicitly.
This can avoid the stringent time step conditions required for the stability of fully explicit time discretizations in the presence of stiff terms while achieving better computational efficiency in comparison with  fully implicit  time discretizations.

Let $0=t^0<t^1<\dots<t^{N_t}=T$ be a uniform mesh in time, with the time step size $\Delta t$ and $t^n=n\Delta t$. An $\IMEXBDForder$th-order IMEX-BDF method of $s$ steps seeks an approximate solution $u^{n+\IMEXBDForder}\approx u(t^{n+\IMEXBDForder})$ to \eqref{eq:IMEX-BDF ode}, by using the numerical solution $u^{n+k}$ at $t^{n+k}$, $k=0,\dots, s-1$, as follows,
\begin{equation}
\label{eq:imex:BDF}
u^{n+\IMEXBDForder} = \sum_{k=0}^{\IMEXBDForder-1} {\IMEXBDFprevcoef}_k u^{n+k} + \Delta t \sum_{k=0}^{\IMEXBDForder-1} {\IMEXBDFexpcoef}_k \IMEXBDFexpop(u^{n+k}) + \Delta t {\IMEXBDFimpcoef}_{\IMEXBDForder} \IMEXBDFimpop(u^{n+\IMEXBDForder}).
\end{equation}
In this work, we will apply the first, second, and third order IMEX-BDF time integrators \cite{ascher1995implicit}  with the parameters $\mathbf{a}=[a_0, a_1, \dots, a_{\IMEXBDForder-1}]^T$, $\mathbf{b}=[b_0, b_1, \dots, b_{\IMEXBDForder-1}]^T$ and $c_s$ given in Table \ref{tab:IMEX-BDF coefficients}. The order condition requires $c_s=\sum_{k=0}^{s-1}b_k$.

% IMEX-BDF: Table of Coefficients
%----------------------------------------------------------------------------------------------------------------
\begin{table}[ht]
\caption{Coefficients of the $\IMEXBDForder$th-order IMEX-BDF method in \eqref{eq:imex:BDF}, $\IMEXBDForder = 1, 2, 3$.}
\centering 
\begin{tabular}{|c|c|c|c|}\hline 
$\IMEXBDForder$ & $\mathbf{a}$ & $\mathbf{b}$ & $c_{\IMEXBDForder}$ \\ \hline 
1 & 1 & 1 & 1 \\ 
\hline 
2 & $[-\frac{1}{3},\frac{4}{3}]^T$ & $[-\frac{2}{3},\frac{4}{3}]^T$ & $\frac{2}{3}$ \\ 
\hline 
3 & $[\frac{2}{11},-\frac{9}{11},\frac{18}{11}]^T$ & $[\frac{6}{11},-\frac{18}{11},\frac{18}{11}]^T$ & $\frac{6}{11}$ \\ 
\hline 
\end{tabular} 
\label{tab:IMEX-BDF coefficients}
\end{table}

Now we return to \eqref{eq:micro-macro}. To apply an IMEX-BDF time integrator to this system, one needs to  specify an {\it IMEX partitioning}, i.e. a strategy that delineates which terms are deemed stiff and 
treated implicitly and which terms are 
deemed non-stiff and 
treated explicitly. With the multi-scale nature of our model and the stiffness being relative, 
{different strategies of IMEX partitionings  can be adopted.  In this work,  we consider three partitioned forms of our model}  (\cite{lemou2008new, jang2015high, peng2021asymptotic}),  specified in Table \ref{imex-stragegy-FG} and referred to as  (IMEX) \textbf{Strategy $k$}, $k=1, 2, 3$, 
and the respective first order in time IMEX-BDF methods of our model are also given here:    
given $\macro^n\approx  \macro(\cdot,t^n)$, $\micro^n\approx \micro(\cdot,\cdot,t^n)$, the numerical solutions $\macro^{n+1}$, $\micro^{n+1}$ at $t^{n+1}$ are sought using one of the following semi-discrete in time methods,

%------------------------------------------
% First Order in Time: 
%----------------------------------------------------------------------------------------------------------------
\medskip
\noindent \textbf{Strategy 1:}
\begin{subequations} \label{eq:strategy 1 temporal discretization}
\begin{align}
\frac{\macro^{n+1} - \macro^n}{\dt} + \partial_x \lvg \velangle \micro^{n+1} \rvg & = -\absorp \macro^n +G, \\
\frac{\micro^{n+1} - \micro^n}{\dt} + \frac{1}{\vareps} (\I - \V) (\velangle \partial_x \micro^n) + \frac{1}{\vareps^2} \velangle \partial_x \macro^{n} & = -\frac{1}{\vareps^2}\scat \micro^{n+1} - \absorp \micro^n,
\end{align}
\end{subequations}
%----------

%----------------------------------------------------------------------------------------------------------------
% First Order in Time: Strategy 1
%----------------------------------------------------------------------------------------------------------------
\noindent\textbf{Strategy 2:}
\begin{subequations} \label{eq:strategy 2 temporal discretization}
\begin{align}
\frac{\macro^{n+1} - \macro^n}{\dt} + \partial_x \lvg \velangle \micro^n \rvg & = -\absorp \macro^n+G, \\
\frac{\micro^{n+1} - \micro^n}{\dt} + \frac{1}{\vareps} (\I - \V) (\velangle \partial_x \micro^n) + \frac{1}{\vareps^2} \velangle \partial_x \macro^{n+1} & = -\frac{1}{\vareps^2}\scat \micro^{n+1} - \absorp \micro^n,
\end{align}
\end{subequations}
%----------------------------------------------------------------------------------------------------------------------------------------------------------------------------
% First Order in Time: Strategy 3
%----------------------------------------------------------------------------------------------------------------
%
\textbf{Strategy 3:}
\begin{subequations} \label{eq:strategy 3 temporal discretization}
\begin{align}
\frac{\macro^{n+1} - \macro^n}{\dt} + \partial_x \lvg \velangle \micro^{n+1} \rvg & = -\absorp \macro^{n+1} +G, \\
\frac{\micro^{n+1} - \micro^n}{\dt} + \frac{1}{\vareps} (\I - \V) (\velangle \partial_x \micro^n) + \frac{1}{\vareps^2} \velangle \partial_x \macro^{n+1} & = -\frac{1}{\vareps^2}\scat \micro^{n+1} - \absorp \micro^{n+1}.
\end{align}
\end{subequations}
The second and third order in time IMEX-BDF methods for our model can similarly follow from \eqref{eq:imex:BDF} and
Tables \ref{tab:IMEX-BDF coefficients}-\ref{imex-stragegy-FG}. 

% IMEX Strategies
%-----------------------
\begin{table}[ht]
\caption{{Three partitioned forms of our model  $\partial_t u = \IMEXBDFexpop(u) + \IMEXBDFimpop(u)$, with $u=[\rho, g]^T$, $\IMEXBDFexpop=\IMEXBDFexpop_1+\IMEXBDFexpop_2$, $\IMEXBDFimpop=\IMEXBDFimpop_1+\IMEXBDFimpop_2$, and $\IMEXBDFexpop_2(u)=[0,- \vareps^{-1} (\I - \V) (\velangle \partial_x \micro) ]^T$, $\IMEXBDFimpop_2(u)=[0,    - 
   \vareps^{-2} \scat \micro]^T$.}}
\centering 
\begin{tabular}{c|c|c}\hline 
& $\IMEXBDFexpop_1(u)$ & {$\IMEXBDFimpop_1(u)$} \\ \hline 
Strategy 1 \cite{lemou2008new}&
$\begin{bmatrix}
    -\absorp \macro\\
 \vareps^{-2}\velangle \partial_x \macro- \absorp \micro
\end{bmatrix}$&
$\begin{bmatrix}
    -\partial_x \lvg \velangle \micro \rvg\\
   0
\end{bmatrix}$\\ 
\hline 
Strategy 2 \cite{jang2015high}&
    $\begin{bmatrix}
    -\partial_x \lvg \velangle \micro \rvg -\absorp \macro\\
    - \absorp \micro
\end{bmatrix}$
&$\begin{bmatrix}
    0\\
    -\vareps^{-2}\velangle \partial_x \macro
\end{bmatrix}$ \\ 
\hline 
Strategy 3 \cite{peng2021asymptotic}&
$\begin{bmatrix}
     0\\
     0
\end{bmatrix}$&
$\begin{bmatrix}
    -\partial_x \lvg \velangle \micro \rvg -\absorp \macro\\
   -
   \vareps^{-2}\velangle \partial_x \macro - \absorp \micro
\end{bmatrix}$ \\ 
\hline 
\end{tabular} 
\label{imex-stragegy-FG}
\end{table}

We would like to make a few remarks about the adopted IMEX partitionings. First of all, as the most stiff term in the diffusive regime when $\vareps \ll 1$, the scattering term $\vareps^{-2}\scat \micro$ is always treated implicitly, while the less stiff and non-symmetric free-streaming operator  $\vareps^{-1}(\I - \V) (\velangle \partial_x \micro)$ is  always treated explicitly.  For the absorption terms $\absorp \macro$ and $\absorp \micro$, it is generally sufficient to treat them explicitly, unless $\absorp$ is very large (see Section \ref{sec:fourier}  for its effect on the stability).
In fact, even when these terms are treated implicitly, the impact on the computational efficiency is minimal given the spatial discretizations to be discussed in {the} next section.
The major differences among the three IMEX  partitionings come from the treatments to $\partial_x \lvg \velangle \micro \rvg$ and $\vareps^{-2}\velangle \partial_x \macro$, as well as the resulting properties of the methods. The resulting differences and similarities will be discussed in the following  sections.

%----------------------------------------------------------------------------------------------------------------
% Subsection: Spatial discretization
%----------------------------------------------------------------------------------------------------------------

\subsection{Spatial discretization} 
\label{sec:spatial discretization}

In space, we apply discontinuous Galerkin (DG) discretizations. Periodic boundary conditions are assumed for now, with general boundary conditions discussed later in {Section \ref{sec:num}.}
We start with a spatial mesh 
 $\Omega_h = \{ I_j=[x_{j-\frac{1}{2}},x_{j+\frac{1}{2}}],j=1,\dots, N_x\}$ for $\spaceset = [x_L,x_R]$, where $x_{\frac{1}{2}} = x_L$ and $x_{N_x+\frac{1}{2}} = x_R$. Let $x_j$ be the midpoint of $I_j$, $h_j$ be its length, and  $h=\max_j{h_j}$. Given a nonnegative integer $r$, we define a finite-dimensional discrete space
$U_h^r=\{u\in L^2(\Omega_x): u|_{I_j}\in P^r(I_j),\forall j=1,\dots, N_x\}$, where $P^r(I_j)$ is the set of  polynomials with degree at most $r$ on $I_j$. Note that functions in $U_h^r$ are double-valued at mesh nodes, and their one-sided traces at $x_{j+\frac{1}{2}}$ are denoted as  $u^{\pm}_{j+\frac{1}{2}}=\lim_{\Delta x\rightarrow 0^{\pm}} u(x_{j+\frac{1}{2} }+\Delta x)$, with the jump {denoted} as $[u]_{j+\frac{1}{2}}=u^+_{j+\frac{1}{2} }-u^-_{j+\frac{1}{2} }$, $\forall j$. 

The DG discretizations as in \cite{jang2015high,peng2021asymptotic} will be adopted. They are based on the discrete derivative operators, $\Dplusop, \Dminusop: \pdspace^{\pd} \rightarrow \pdspace^{\pd}$, with each approximating the spatial differentiation  $\partial_x$, as well as the discrete derivative operator, $\Dupwindop(\cdot;v): \pdspace^{\pd} \rightarrow \pdspace^{\pd}$, approximating  $v\partial_x$ in an upwind fashion. More specifically, these  discrete operators are defined as follows, $\forall \testmacro, \testmicro \in \pdspace^{\pd}$,
\begin{subequations}
\begin{align}
(\Dminusop \testmacro, \testmicro) &:= -\int_{\spaceset} \testmacro \partial_x \testmicro dx -\sum_{j} \breve{\testmacro}_{j-\frac{1}{2}} [\testmicro]_{j-\frac{1}{2}},\\
(\Dplusop \testmacro, \testmicro) &:= -\int_{\spaceset} \testmacro \partial_x \testmicro dx -\sum_{j} \widehat{\testmacro}_{j-\frac{1}{2}} [\testmicro]_{j-\frac{1}{2}},\\
(\Dupwindop(\testmacro;\velangle),\testmicro)&:= -\int_{\spaceset} \velangle \testmacro \partial_x \testmicro dx - \sum_{j} \widetilde{(\velangle \testmacro)}_{j-\frac{1}{2}} [\testmicro]_{j-\frac{1}{2}}.
\end{align}
\end{subequations}
Here $(\cdot,\cdot)$ is the $L^2$ inner product for $L^2(\spaceset)$. The numerical fluxes at the element interfaces are taken as   $\breve{\testmacro}_{j-\frac{1}{2}}=\testmacro^-_{j-\frac{1}{2}}$, $\widehat{\testmacro}_{j-\frac{1}{2}}=\testmacro^+_{j-\frac{1}{2}}$, {and in addition, $\widetilde{(\velangle \testmacro)}$ is the upwind flux, defined as $\widetilde{(\velangle \testmacro)}_{j-\frac{1}{2}}=(\velangle \testmacro)_{j-\frac{1}{2}}^{-}$ if $v\geq 0$ and  $\widetilde{(\velangle \testmacro)}_{j-\frac{1}{2}}=(\velangle \testmacro)_{j-\frac{1}{2}}^{+}$ if $v<0$.}

With periodic boundary conditions, the following property can be verified directly:  {$(\Dplusop \testmacro, \testmicro)=-(\testmacro,\Dminusop \testmicro)$, $\forall \testmacro, \testmicro\in \pdspace^{\pd}$, or equivalently, in the operator notation}
\begin{equation}
   \Dplusop  = -(\Dminusop)^T.
   \label{eq:adjoint}
\end{equation}

In this work, we adopt the following discretizations in space, 
\begin{equation}
    \Dplusop \lvg \velangle \micro \rvg \approx \partial_x \lvg \velangle \micro \rvg,  \qquad \Dminusop \macro \approx \partial_x \macro, \qquad
    \Dupwindop(\micro;\velangle) \approx\velangle \partial_x \micro. 
\end{equation}
Alternatively, one can use  $\Dminusop \lvg \velangle \micro \rvg \approx \partial_x \lvg \velangle \micro \rvg$ and $\Dplusop \macro \approx \partial_x \macro$, along with $
    \Dupwindop(\micro;\velangle) \approx\velangle \partial_x \micro. $

%----------------------------------------------------------------------------------------------------------------
% Subsection: Velocity discretization
%----------------------------------------------------------------------------------------------------------------

\subsection{Velocity discretization} 
\label{sec:velocity discretization}

In velocity, we use the discrete ordinates method \cite{lewis1984computational}.
Let $\{\velangle_q\}_{q=1}^{N_{\velangle}}$, $\{\vweight_q\}_{q=1}^{N_{\velangle}}$ denote the sets of quadrature points and weights on $\velangleset = [-1,1]$, respectively. The numerical solution $\micro_h$ in the $\velangle$ variable will be sought at $\velangle_q$,  denoted as $\micro_{h,q}$, $ q=1, \dots, N_{\velangle}$,  while the integral operator 
$\lvg \cdot \rvg$ will be approximated by 
$\lvg \cdot \rvg_h$  using numerical integration, namely
\begin{equation} \label{eq:velocity integral approx}
\lvg \eta(\velangle) \rvg_h: = \sum_{q=1}^{N_{\velangle}} \vweight_q \eta(\velangle_q).
\end{equation}
We require
\begin{equation} \label{eq:v square integral exact}
\lvg \velangle^2 \rvg = \lvg \velangle^2 \rvg_h,    
\end{equation}
a property that is easy to satisfy,
and this is to ensure the correct asymptotic limit of the proposed methods as $\vareps \rightarrow 0$  (see Section \ref{sec:AP}).
%

%----------------------------------------------------------------------------------------------------------------
% Subsection: Fully discrete scheme
%----------------------------------------------------------------------------------------------------------------

\subsection{Fully discrete schemes}
\label{sec:fully discrete scheme}

By combining the temporal, spatial, and velocity discretizations, we now present the three families of fully discrete schemes proposed in this work with respect to three IMEX  partitionings in Section \ref{sec:temporal discretization}.  

We first consider the case with  constant scattering and absorption cross sections, $\scat$ and $\absorp$. Given $\macro_h^{n+k}$ and $\{\micro_{h,q}^{n+k}\}_{q=1}^{N_{\velangle}}$ in $\pdspace^{\pd}$, $k = 0, 1, \dots, \IMEXBDForder-1$, we seek $\macro_h^{n+s}$ and $\{\micro_{h,q}^{n+s}\}_{q=1}^{N_{\velangle}}$ in  $\pdspace^{\pd}$, such that

%----------------------------------------------------------------------------------------------------------------
% Strategy 1: Fully discrete scheme
%----------------------------------------------------------------------------------------------------------------
\medskip
\textbf{Strategy 1:}
\begin{subequations} \label{eq:fully:st1}
\begin{align}
\macro_h^{n+\IMEXBDForder} =& \sum_{k=0}^{\IMEXBDForder-1} \IMEXBDFprevcoef_k \macro_h^{n+k} - \dt \sum_{k=0}^{\IMEXBDForder-1} \IMEXBDFexpcoef_k (\absorp \macro_h^{n+k}-G_h) - \dt \IMEXBDFimpcoef_{\IMEXBDForder} \Dplusop \lvg \velangle \micro_h^{n+\IMEXBDForder} \rvg_h, \label{eq:fully:st1.a}\\
\micro_{h,q}^{n+\IMEXBDForder} = &\sum_{k=0}^{\IMEXBDForder-1} \IMEXBDFprevcoef_k \micro_{h,q}^{n+k} - \dt \sum_{k=0}^{\IMEXBDForder-1} \IMEXBDFexpcoef_k \Big[ \frac{1}{\vareps} \big(\Dupwindop(\micro_{h}^{n+k};\velangle_q)-\lvg 
\Dupwindop(\micro_{h}^{n+k};\velangle)\rvg_h \big) \notag \\
&+  \frac{1}{\vareps^2}\velangle_q \Dminusop \macro_h^{n+k} + \absorp \micro_{h,q}^{n+k} \Big] - \dt \IMEXBDFimpcoef_{\IMEXBDForder} \frac{1}{\vareps^2}\big(\scat \micro_{h,q}^{n+\IMEXBDForder}\big), \;\;q = 1, \dots, N_{\velangle},\label{eq:fully:st1.b}
\end{align}
\end{subequations}

%----------------------------------------------------------------------------------------------------------------
% Strategy 2: Fully discrete scheme
%----------------------------------------------------------------------------------------------------------------
\textbf{Strategy 2:}
\begin{subequations} \label{eq:fully:st2}
\begin{align}
\macro_h^{n+\IMEXBDForder} =& \sum_{k=0}^{\IMEXBDForder-1} \IMEXBDFprevcoef_k \macro_h^{n+k} - \dt \sum_{k=0}^{\IMEXBDForder-1} \IMEXBDFexpcoef_k (\Dplusop \lvg \velangle \micro_h^{n+k} \rvg_h+\absorp \macro_h^{n+k}-G_h),\label{eq:fully:st2.a} \\
\micro_{h,q}^{n+\IMEXBDForder} = &\sum_{k=0}^{\IMEXBDForder-1} \IMEXBDFprevcoef_k \micro_{h,q}^{n+k} - \dt \sum_{k=0}^{\IMEXBDForder-1} \IMEXBDFexpcoef_k \Big[ \frac{1}{\vareps} \big(\Dupwindop(\micro_{h}^{n+k};\velangle_q)-\lvg \Dupwindop(\micro_{h}^{n+k};\velangle)\rvg_h \big) \notag \\
&+ \absorp \micro_{h,q}^{n+k} \Big] - \dt \IMEXBDFimpcoef_{\IMEXBDForder} \frac{1}{\vareps^2}\big(\velangle_q \Dminusop \macro_h^{n+s}+\scat \micro_{h,q}^{n+\IMEXBDForder}\big), \;\;q = 1, \dots, N_{\velangle},\label{eq:fully:st2.b}
\end{align}
\end{subequations}

%----------------------------------------------------------------------------------------------------------------
% Strategy 3: Fully discrete scheme
%----------------------------------------------------------------------------------------------------------------
\textbf{Strategy 3:}
\begin{subequations} \label{eq:fully:st3}
\begin{align}
\macro_h^{n+\IMEXBDForder} =& \sum_{k=0}^{\IMEXBDForder-1} \IMEXBDFprevcoef_k \macro_h^{n+k} - \dt \IMEXBDFimpcoef_{\IMEXBDForder} \big(\Dplusop \lvg \velangle \micro_h^{n+\IMEXBDForder}\rvg_h+\absorp \macro_h^{n+\IMEXBDForder}-G_h \big), \label{eq:fully:st3.a}\\
\micro_{h,q}^{n+\IMEXBDForder} = &\sum_{k=0}^{\IMEXBDForder-1} \IMEXBDFprevcoef_k \micro_{h,q}^{n+k} - \dt \sum_{k=0}^{\IMEXBDForder-1} \IMEXBDFexpcoef_k \frac{1}{\vareps} \Big[  \Dupwindop(\micro_{h}^{n+k};\velangle_q)-\lvg \Dupwindop(\micro_{h}^{n+k};\velangle)\rvg_h \Big] \notag \\
& - \dt \IMEXBDFimpcoef_{\IMEXBDForder} \frac{1}{\vareps^2}\big(\velangle_q \Dminusop \macro_h^{n+s}+\scat \micro_{h,q}^{n+\IMEXBDForder}+\vareps^2\absorp \micro_{h,q}^{n+\IMEXBDForder} \big), \;\;q = 1, \dots, N_{\velangle}.\label{eq:fully:st3.b}
\end{align}
\end{subequations}
Again with the order condition $\IMEXBDFimpcoef_{\IMEXBDForder}=\sum_{k=0}^{\IMEXBDForder-1} \IMEXBDFexpcoef_k$, the {time-independent} source term $\source_h:=\mathcal{P}_h(\source)\in \pdspace^{\pd}$ can appear either in the explicit terms or the implicit terms in any of the schemes above, with the schemes unchanged.\footnote{{Our schemes can be easily adapted to the case with the time-dependent source term.}} Here $\mathcal{P}_h$ is the $L^2$ projection operator onto $\pdspace^{\pd}$. The initialization for $\macro_h^0(\cdot)$ and  $\micro_h^0(\cdot,v)$ can also be done via the $L^2$ projection  $\mathcal{P}_h$. For multi-step methods with $\IMEXBDForder>1$, numerical solutions at $t^j$ ($j=1, \dots, \IMEXBDForder-1$)
are also needed, and this will be examined in Section \ref{sec:init}.

We now consider the more general case with the scattering cross section $\scat=\scat(x)$ and the absorption cross section $\absorp=\absorp(x)$,  when the material properties are spatially dependent. Following the modal DG and nodal DG methodologies to treat the variable coefficient or nonlinear models \cite{hesthaven2007nodal}, the schemes can be formulated as above by  replacing the terms $\absorp \macro_h^{n+k}$, $\absorp \micro_h^{n+k}$, $\scat \micro_h^{n+k}$ in \eqref{eq:fully:st1}-\eqref{eq:fully:st3} either by 
\begin{equation}
    \mP_h(\absorp \macro_h^{n+k}),\quad \mP_h(\absorp \micro_h^{n+k}),
    \quad \mP_h(\scat \micro_h^{n+k})
    \label{eq:varying-crossS1}
\end{equation}
as in the modal DG setting, or by 
\begin{equation}
\mI_h(\absorp \macro_h^{n+k}),\quad \mI_h(\absorp \micro_h^{n+k}),
    \quad \mI_h(\scat \micro_h^{n+k})
    \label{eq:varying-crossS2}
\end{equation}
as in the nodal DG setting. While 
$\mP_h$ is the $L^2$ projection operator onto $\pdspace^{\pd}$, $\mI_h$ is the interpolation operator onto $\pdspace^{\pd}$ with respect to a set of interpolation points, e.g. the set of scaled $(r+1)$  Gauss-Legendre quadrature points (see Section 3 in \cite{peng2020asymptotic}).

From now on, IMEX-BDF$\IMEXBDForder$-DG$\pd$ will denote the proposed scheme using the $\IMEXBDForder$th order IMEX-BDF method ($\IMEXBDForder=1,2,3$) with the DG method involving the discrete space $\pdspace^{\pd-1}$ (and hence with the formal  $\pd$th order accuracy in space). We also use  IMEX-BDF$\IMEXBDForder$-DG$\pd$-$\mS$$k$ to indicate that the IMEX   Strategy $k$ is adopted. For later reference, we here explicitly write down the IMEX-BDF1-DG$r$-$\mS$1 scheme with the modal treatment \eqref{eq:varying-crossS1} for general cross sections $\scat(x)$ and $\absorp(x)$: given $\macro_h^n$ and $\{\micro_{h,q}^n\}_{q=1}^{N_{\velangle}}$ in $\pdspace^{\pd}$,  we seek $\macro_h^{n+1}$ and $\{\micro_{h,q}^{n+1}\}_{q=1}^{N_{\velangle}}$ in  $\pdspace^{\pd}$, such that
\begin{align*}
\frac{\macro^{n+1}_h - \macro^n_h}{\dt} + \Dplusop \lvg \velangle \micro^{n+1}_h \rvg_h & = -\mP_h(\absorp \macro^n_h) +G_h, \\
\frac{\micro^{n+1}_{h,q} - \micro^n_{h,q}}{\dt} + \frac{1}{\vareps} \Big[\Dupwindop(\micro_{h}^n; v_q)-\lvg\Dupwindop(\micro_h^n; v)\rvg_h\Big] + \frac{1}{\vareps^2} \velangle_q \Dminusop \macro^{n}_h & = -\frac{1}{\vareps^2}\mP_h(\scat \micro^{n+1}_{h,q}) - \mP_h(\absorp \micro^n_{h,q}).
\end{align*}
The scheme above is in the {\it strong} form, which has an equivalent {\it weak} form: 
$\forall \testmacro, \testmicro\in \pdspace^{\pd}$,
\begin{subequations} \label{eq:BDF1:DGr:st1:w}
\begin{align}
\Big(\frac{\macro_h^{n+1} - \macro_h^n}{\dt},\testmacro\Big) + l_h(\lvg \velangle \micro_h^{n+1} \rvg_h,\testmacro) & = -(\absorp \macro_h^n,\testmacro)+(G, \testmacro), \label{eq:BDF1:DGr:st1:w.a}\\
\vareps^2 \Big(\frac{\micro_{h,q}^{n+1} - \micro_{h,q}^n}{\dt},\testmicro\Big) + \vareps b_{h,\velangle_q} (\micro_{h}^n,\testmicro) - \velangle_q d_h(\macro_h^n,\testmicro) & = -(\scat \micro_{h,q}^{n+1},\testmicro) - \vareps^2 (\absorp \micro_{h,q}^n,\testmicro). \label{eq:BDF1:DGr:st1:w.b}
\end{align}
\end{subequations}
Here $l_h(\cdot,\cdot)$, $d_h(\cdot,\cdot)$ and $b_{h,\velangle}(\cdot,\cdot)$ are defined as follows:
%------------------------------------------
\begin{equation}
\label{eq:bilinear:a}
l_h(\testmacro, \testmicro)=(\Dplusop \testmacro, \testmicro),  \quad
d_h(\testmacro, \testmicro)=-(\Dminusop \testmacro, \testmicro),
\end{equation}
and
\begin{equation}\label{eq:bilinear:b}
b_{h,\velangle}(\micro_{h},\testmicro) = (\Dupwindop(\micro_{h};\velangle) - \lvg \Dupwindop(\micro_h;\velangle) \rvg_h,\testmicro).
\end{equation}

For each proposed scheme, the following property holds through a similar proof of Lemma 3.1 in \cite{jang2014analysis},
\begin{equation} \label{eq:micro v integral zero property}
\lvg \micro_h^n \rvg_h = \lvg \micro_h^n \rvg = 0, \quad \forall n \geq s,
\end{equation}
provided that the initialization procedure ensures $\lvg \micro_h^k \rvg_h=0$, $k=0, 1, \dots, s-1.$

%----------------------------------------------------------------------------------------------------------------
% Subsection: AP property
%----------------------------------------------------------------------------------------------------------------

%----------------------------------------------------------------------------------------------------------------
% Section: Stability
%----------------------------------------------------------------------------------------------------------------

\section{Stability}
\label{sec:stability}

{
In this section, we assume the boundary conditions in space are periodic and the source $G$ is zero, and   the stability of our schemes will be investigated, both analytically and numerically. This is the discrete analogue of an $L^2$ energy relation of the exact solution to \eqref{eq:kinetic transport equation} or to \eqref{eq:micro-macro}:
\begin{subequations}
    \label{eq:stab:model}
\begin{align}
&\frac{1}{2} 
    \frac{d}{dt} \int_{\spaceset\times \velangleset}\angflux^2dvdx=-\frac{1}{\vareps^2} \int_{\spaceset\times \velangleset}\scat\;\big(\angflux-\langle\angflux\rangle\big)^2 dvdx-\int_{\spaceset\times \velangleset}\absorp\;\angflux^2dvdx\leq 0,\\
\Leftrightarrow&\frac{1}{2} 
    \frac{d}{dt} \left(\int_{\spaceset}\macro^2dx+ \vareps^2\int_{\spaceset\times \velangleset}\micro^2dvdx\right)\notag\\
    &\hspace{0.5in}=- \int_{\spaceset\times \velangleset}\scat\;\micro^2 dvdx-\left(\int_{\spaceset}\absorp\; \macro^2dx+\vareps^2\int_{\spaceset\times \velangleset}\absorp\; \micro^2 dvdx\right)\leq 0.
\end{align}
\end{subequations}
Here $\langle \micro\rangle=0$ is used.} 
Particularly, energy-based stability analysis
is established for the first order IMEX-BDF1-DG1 schemes  in Section \ref{sec:stab:energy}, and Fourier-based stability analysis is carried out for the IMEX-BDF$\IMEXBDForder$-DG$\IMEXBDForder$ schemes, $s=1, 2, 3$ 
(up to third order accuracy) 
in Section \ref{sec:fourier}. 
 In  addition,  it is assumed that $ \scat, \absorp\in L^\infty(\spaceset)$, and they satisfy
\begin{equation}
    0 < \scatlower \leq \scat(x),\quad 0 \leq 
\absorp(x) \leq \absorpupper<\infty, \quad \text{a.e.  on}\;  \spaceset.
\label{eq:sig:bound}
\end{equation}

%%%%%%%%%%%%%%%%%
%\input{Stab1.2}
%%%%%%%%%%%%%%%%%%

%----------------------------------------------------------------------------------------------------------------
% Subsection: Energy analysis
%----------------------------------------------------------------------------------------------------------------

\subsection{Energy-based stability analysis}
\label{sec:stab:energy}

In this section, we perform energy-based stability analysis for the first order in space and time schemes. Similar analysis can be carried out for the schemes with  first order accuracy in time and higher order accuracy in space  as in \cite{jang2014analysis}, yet it is nontrivial to extend such analysis to higher order in time schemes  due to the multi-scale nature of the problem. 
Notation wise,  
we use 
$\|\phi\|=\sqrt{(\phi,\phi)}$, $\|\phi\|_a=\sqrt{(\absorp\phi,\phi)}$, $|||\psi||| = \sqrt{\lvg (\psi, \psi) \rvg_h}$, $|||\psi|||_s = \sqrt{\lvg (\scat\psi, \psi) \rvg_h}$, $|||\psi|||_a = \sqrt{\lvg (\absorp\psi, \psi) \rvg_h}$, whenever they are well-defined for $\phi(x)$ and $\psi(x,v)$.
Without loss of generality, we also assume the mesh to be uniform with $\dx_j = \dx$, $\forall j$.  The main results are given in  Theorem \ref{thm:stab:energy}, with the proof detailed in {Appendix \ref{sec:stab:energy:proof}.} With minor modification, our results can be extended
to general meshes when $\max_j{\dx_j}/\min_j{\dx_j}$ is uniformly bounded during the mesh refinement. In our presentation, we choose to explicitly state the dependence on $\velangleset$ such as through $||\velangle||_{h,\infty}:=\max_{1\leq q\leq N_v}|v_q|$. 

\begin{theorem} For the proposed first order in space and time schemes, the following hold regarding their stability. 
\begin{itemize}
\item [1.)] With Strategy 1, the IMEX-BDF1-DG1-$\mS$1 scheme is stable in the sense that the discrete energy $E_{h,\mS 1}^n$ does not grow in time, namely 
\begin{equation}\label{eq:st1:stab-ene}
    E_{h,\mS1}^{n+1}\leq E_{h,\mS1}^n,\quad \textrm{with}\;\; 
    {E_{h,\mS1}^n=|| \macro_h^n ||^2 + \vareps^2 ||| \micro_h^n |||^2,}
\end{equation}
under the time step condition
\begin{equation} \label{eq:st1:dt}
 \dt \leq \frac{\dt_s }{1+0.5\absorpupper \dt_s}.
\end{equation}
Here 
\begin{equation}
\label{eq:st1:dt.a}
\dt_s=\frac{2 || \velangle ||_{h,\infty} \vareps\dx + \scatlower \dx^2}{2||\velangle ||_{h,\infty} (|| \velangle ||_{h,\infty} + \lvg \left| \velangle \right| \rvg_h)}. 
\end{equation}
\item [2.)] With Strategy 2, the IMEX-BDF1-DG1-$\mS$2 scheme is stable in the sense that the discrete energy $E_{h,\mS 2}^n$ does not grow in time, namely 
\begin{equation}\label{eq:st2:stab-ene}
    E_{h,\mS2}^{n+1}\leq E_{h,\mS2}^n,\quad \textrm{with}\;\; 
    {E_{h,\mS2}^n=|| \macro_h^n ||^2 + \vareps^2 ||| \micro_h^{n-1} |||^2,}
\end{equation}
under the same time step condition as in \eqref{eq:st1:dt}-\eqref{eq:st1:dt.a}.
\item [3.)] With Strategy 3, and
\begin{equation}
    \label{eq:st3:ene}
E_{h,\mS3}^n=|| \macro_h^n ||^2 + \vareps^2 ||| \micro_h^n |||^2+\dt ||| \micro_h^n |||_s^2,\end{equation}
the IMEX-BDF1-DG1-$\mS$3 scheme is stable in the sense that 
\begin{equation}\label{eq:st3:stab}
E_{h,\mS3}^{n+1}+2\dt(|| \macro_h^{n+1} ||_a^2 + \vareps^2 ||| \micro_h^{n+1} |||_a^2)\leq E_{h,\mS3}^n, \quad \textrm{hence}\quad  E_{h,\mS3}^{n+1}\leq E_{h,\mS3}^n
\end{equation}
holds
\begin{itemize}
    \item[3.i)] unconditionally when $
\frac{\vareps}{\scatlower \dx} \leq \frac{1}{2 || \velangle ||_{h,\infty}}$, namely, with any time step size; otherwise, 

    \item[3.ii)] it  holds conditionally when $\frac{\vareps}{\scatlower \dx} >\frac{1}{2||\velangle ||_{h,\infty}}$, under the time step condition
\begin{equation} \label{eq:st3:dt}
\dt \leq \frac{2 \vareps^2 \dx}{2 || \velangle ||_{h,\infty}\vareps  - \scatlower \dx}.
\end{equation}

\end{itemize}
\end{itemize}
\label{thm:stab:energy}
\end{theorem}

\begin{remark} {There are two  factors that contribute to the stability of our methods. One is the dissipation  inherited from the scattering operator, as in  \eqref{eq:stab:model} and \eqref{eq:stab:1.b}. The other is the numerical dissipation due to the upwind treatment of the transport term $\velangle \partial_x \micro$, also see  \eqref{eq:stab:4.a}.}
\label{rem:stab:energy:1:rev}
\end{remark}
\begin{remark}
{Note that the kind of stability in Theorem \ref{thm:stab:energy}, namely some discrete energy does not grow in time, is referred to as ``monotonicity stability'' in literature \cite{xu20192}.  With any strategy of the three IMEX  partitionings, a hyperbolic time step condition $\dt=O(\vareps h)$ is required for stability in the transport regime with $\vareps=O(1)$. In the diffusive regime with $\vareps\ll 1$, IMEX Strategy 3 leads to unconditional stability, while a diffusion type time step condition $\dt=O(h^2)$ is required with Strategy 1 and 2.} {As to be elaborated in Section \ref{sec:alg}, the improvement in the stability of Strategy 3 comes at an expense of solving some linear system, while methods with Strategy 1 and 2 are essentially explicit.}
    \label{rem:stab:energy:1}
\end{remark}
%%%%%%%%%%%%%%%
\begin{remark} {With the model \eqref{eq:kinetic transport equation} under a diffusive scaling, $\absorp=O(1)$ is assumed, and treating $\absorp$-terms explicitly and implicitly in each proposed method involves comparable computational costs over one time step. If $\absorp$-terms are  treated} implicitly when either Strategy 1 or 2 is adopted, the stability result for the  IMEX-BDF1-DG1-$\mS$$k$ scheme ($k=1,2$) 
will be modified to 
\begin{equation*}
    E_{h,\mS k}^{n+1}+2\dt(|| \macro_h^{n+1} ||_a^2 + \vareps^2 ||| \micro_h^{n+1} |||_a^2) \leq E_{h,\mS k}^n,\;\;\textrm{hence}\; 
      E_{h,\mS k}^{n+1}\le   E_{h,\mS k}^{n}, 
\end{equation*}
under an improved time step condition  $\dt \leq \dt_s$, with $E_{h,\mS1}^n$, $E_{h,\mS2}^n$, $\dt_s$ defined in \eqref{eq:st1:stab-ene}, \eqref{eq:st2:stab-ene}, and \eqref{eq:st1:dt.a}, respectively. 
More insight will be gained about the  contribution  of $\absorp$ along with its numerical treatments to the stability  in the next section through Fourier-based stability analysis.
\end{remark}
%%%%%%%%%%%%%%
\begin{remark}
 Once the energy-based stability is available, one can further carry out error analysis similarly as in \cite{jang2014analysis} by combining the consistency of the methods and the approximation properties of the discrete space.  
\end{remark}

%----------------------------------------------------------------------------------------------------------------
% Subsection: Fourier analysis
%----------------------------------------------------------------------------------------------------------------

\subsection{Stability  by Fourier Analysis}
\label{sec:fourier}

\begin{figure}[ht]
\centering
\begin{subfigure}{0.32\textwidth}
\includegraphics[width=\linewidth,height=4cm]
{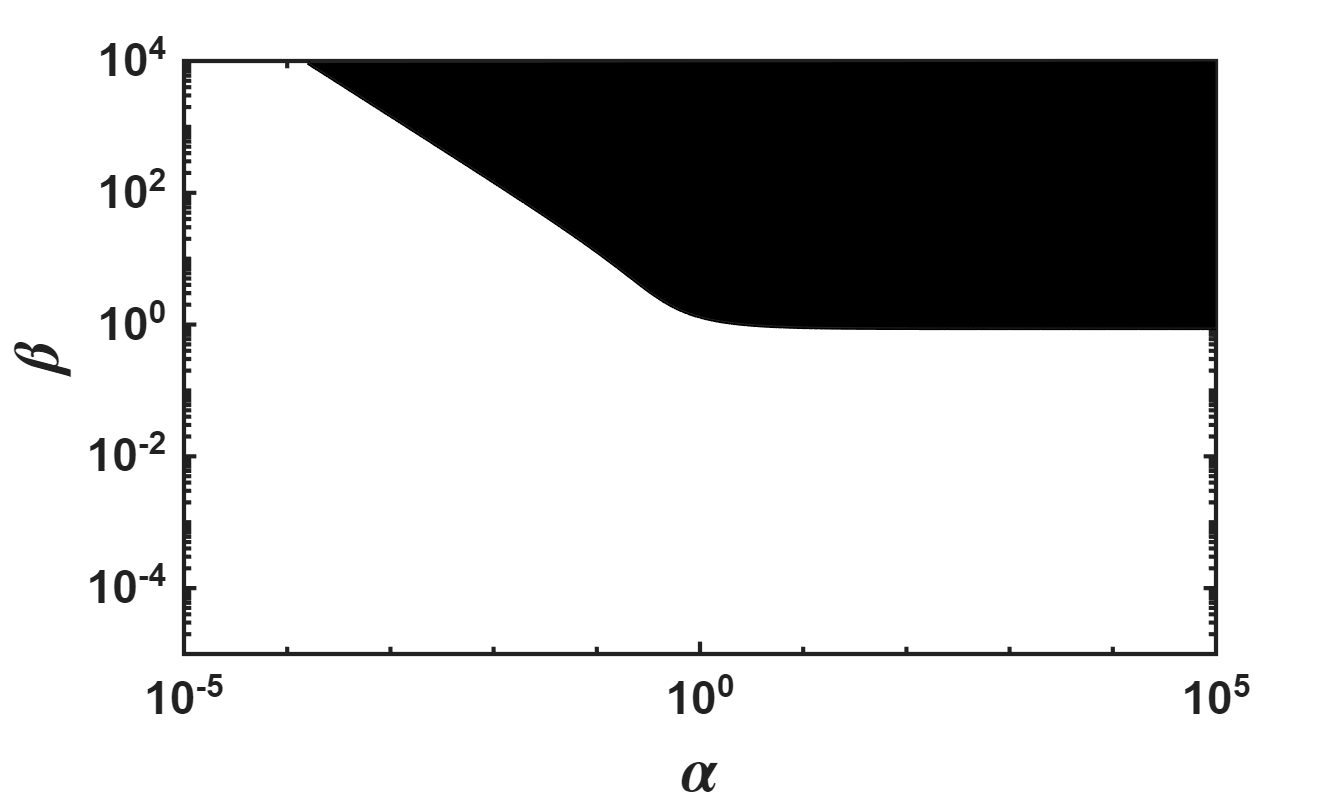}
\caption{IMEX-BDF1-DG1}
\end{subfigure}
\begin{subfigure}{0.32\textwidth}
\includegraphics[width=\linewidth,height=4cm]{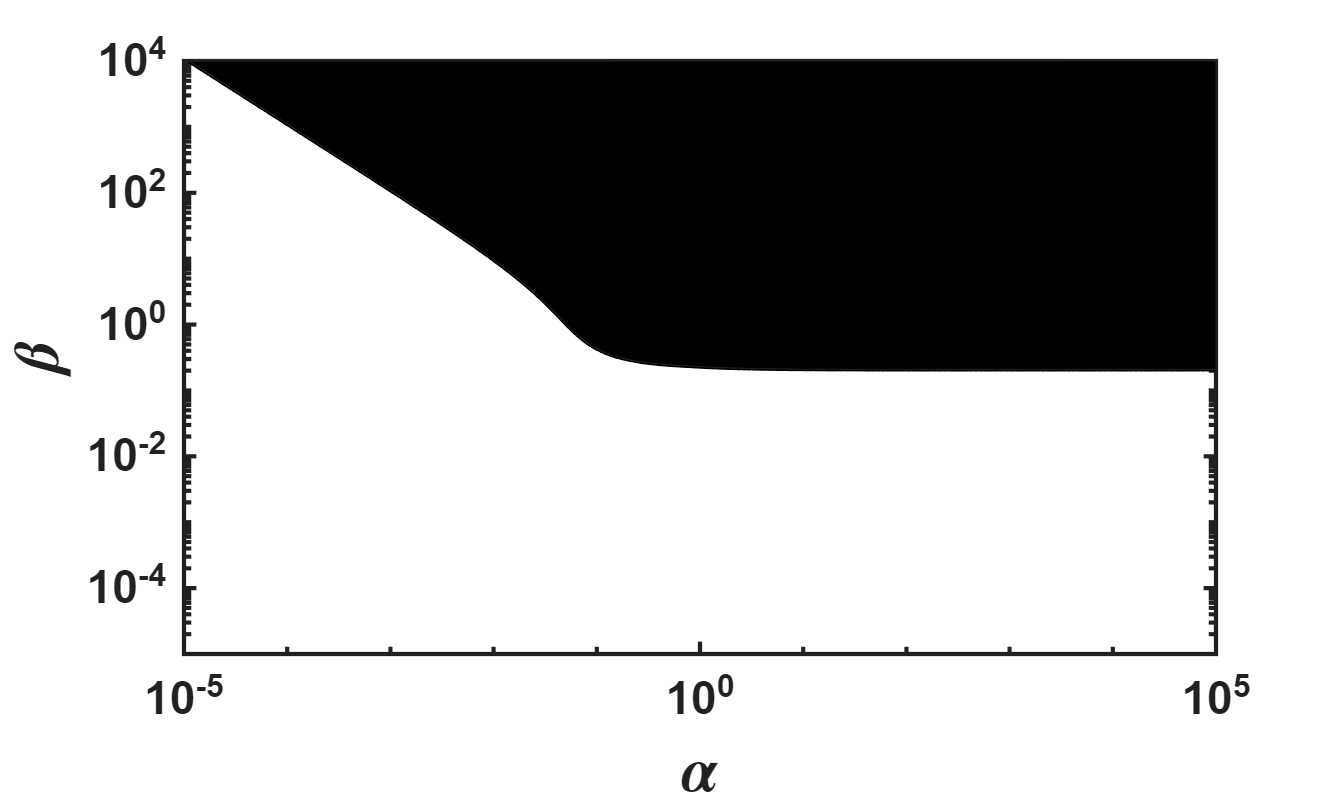}
\caption{IMEX-BDF2-DG2}
\end{subfigure}
\begin{subfigure}{0.32\textwidth}
\includegraphics[width=\linewidth,height=4cm]{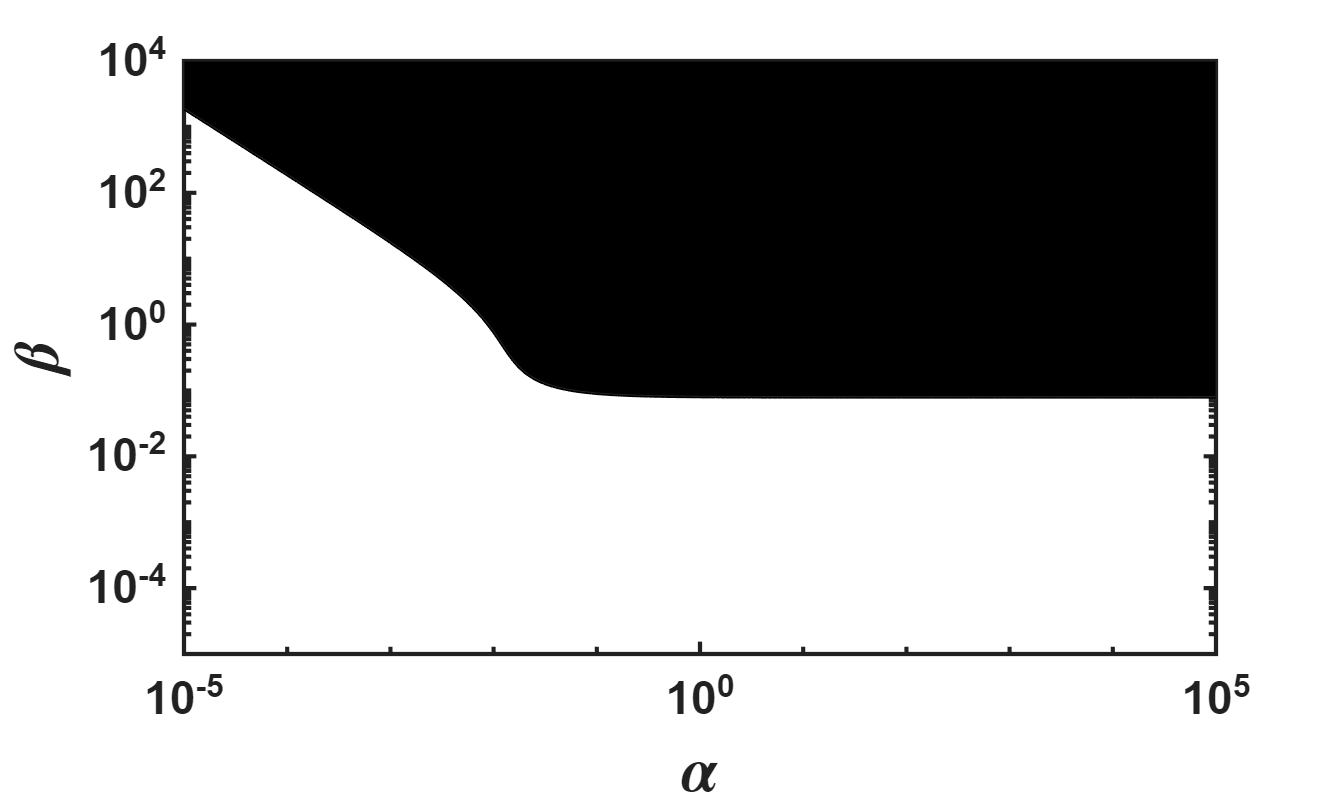}
\caption{IMEX-BDF3-DG3}
\end{subfigure}

\begin{subfigure}{0.32\textwidth}
\includegraphics[width=\linewidth,height=4cm]{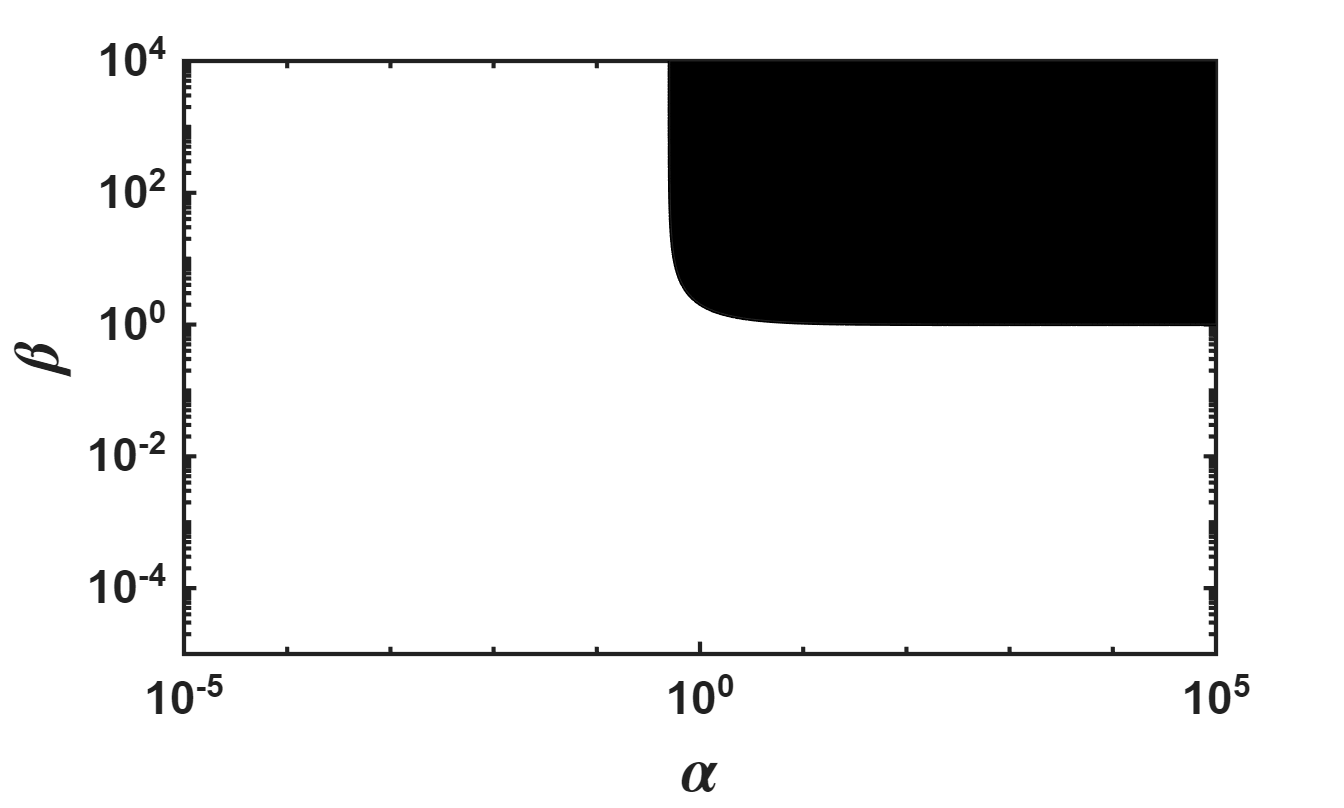}
\caption{IMEX-BDF1-DG1-$\mS$3}
\end{subfigure}
\begin{subfigure}{0.32\textwidth}
\includegraphics[width=\linewidth,height=4cm]{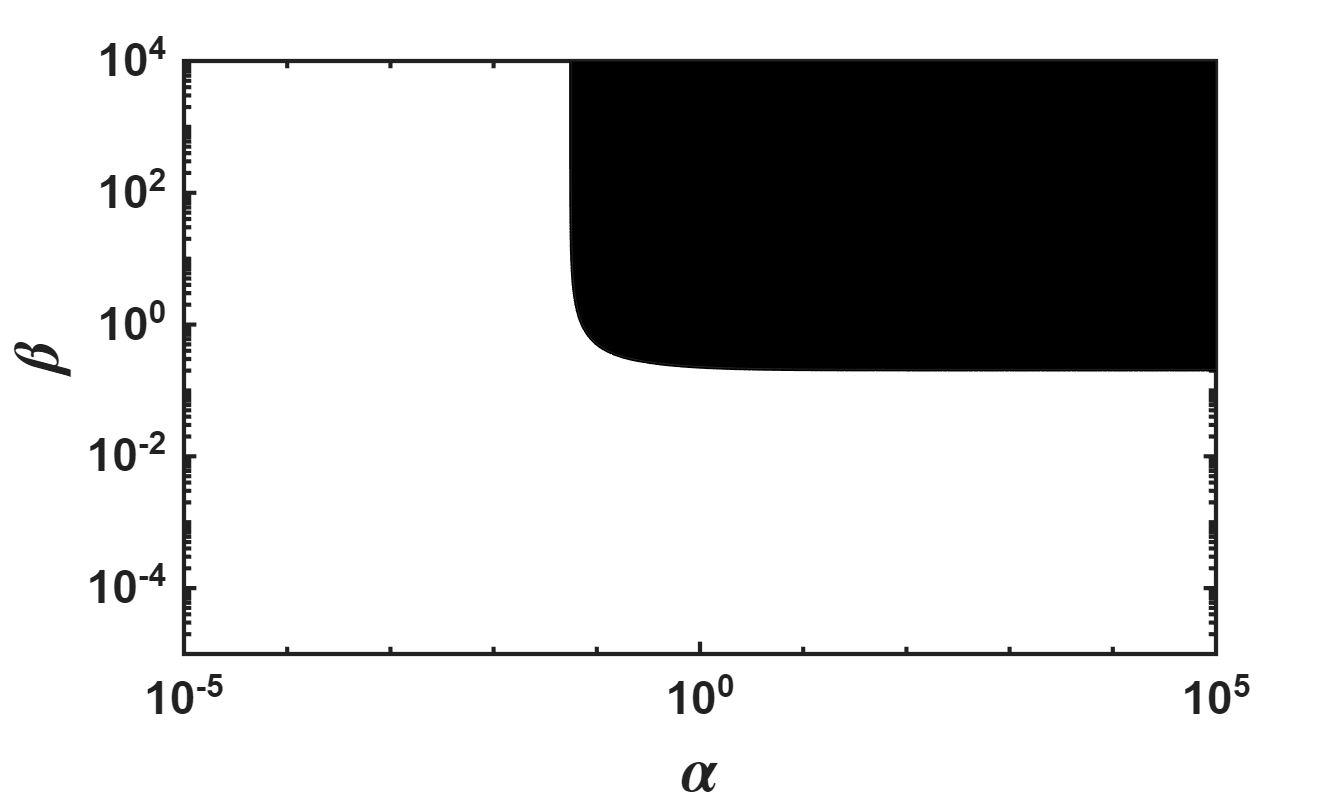}
\caption{IMEX-BDF2-DG2-$\mS$3}
\end{subfigure}
\begin{subfigure}{0.32\textwidth}
\includegraphics[width=\linewidth,height=4cm]{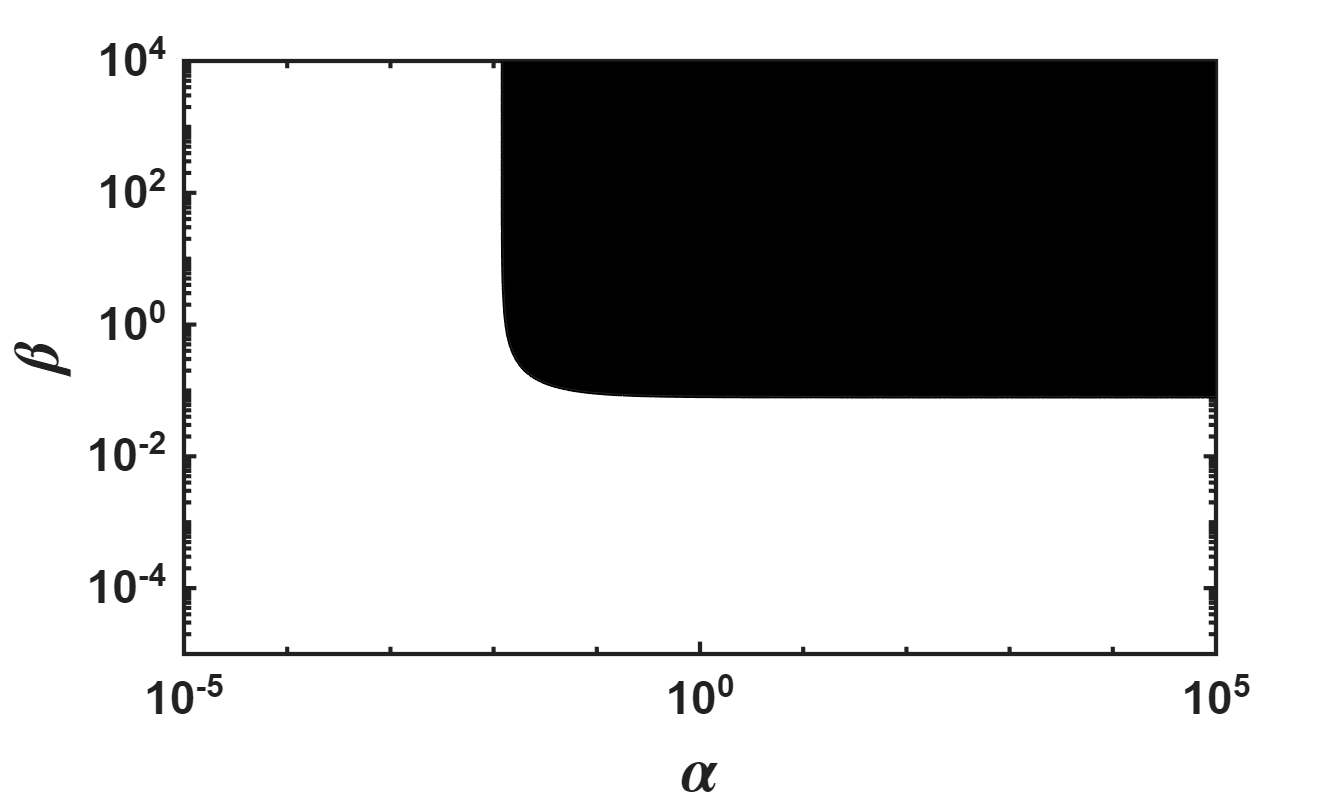}
\caption{IMEX-BDF3-DG3-$\mS$3}
\end{subfigure}

\caption{Stability regions of IMEX-BDF$\pd$-DG$\pd$ with Strategy 1 (and 2) (top row) and Strategy 3 (bottom row) for $\pd=1,2,3$; White: Stable, Black: Unstable; $\alpha = \vareps/(\scatlower  h), \beta = \dt/(\vareps h)$.}
\label{fig:fourier:st123}
\end{figure}

The energy analysis in the previous section provides  sufficient time step conditions for stability when the schemes are applied to the model in all regimes (e.g. scattering dominant, transport dominant) with constant or variable coefficients  $\scat(x), \absorp(x)$ and when the meshes are uniform. The analysis can  be easily extended to non-uniform meshes.  The energy-based stability analysis, however, is only available to the schemes with  first order accuracy in time. To get some insights for the stability of the schemes of higher order accuracy in time, in this section, we carry out Fourier-based stability analysis.
As is standard for such analysis, we assume the spatial mesh is uniform, the boundary conditions are periodic, and $\scat(x) = \scatlower>0$ {and $\absorp(x) = \absorpupper\geq 0$.}  We consider the  one-group
transport equation in slab geometry, and   $16$-point {normalized} Gauss-Legendre quadrature (with $N_{\velangle}=16$) is applied  to approximate $\lvg \angflux \rvg$.  {Note that the  quantitative results in this section,  e.g., Figure \ref{fig:fourier:st123}, \eqref{eq:fourier:dt:st12}-\eqref{eq:fourier:dt:st3}, depend on $N_v$, as previously revealed  by the energy-based stability results in Theorem  \ref{thm:stab:energy}.}

We will start with the case of $\absorp(x) = 0$. 
To set up the stage, let the numerical solution of the IMEX-BDF$\IMEXBDForder$-DG$\pd$ scheme on {the mesh element} $I_m$ be 
\begin{equation}
\macro_h^n(x) \Big|_{I_m} = \sum_{l=0}^{\pd-1} \macro_{m,l}^n \Legendre_l^m(x),\quad \micro_{h,q}^n(x) \Big|_{I_m}  = \sum_{l=0}^{\pd-1} \micro_{q,m,l}^n \Legendre_l^m(x),  \;\;q  = 1, 2, \dots, N_{\velangle},
\label{eq:fourier:0}
\end{equation}
where $\{\Legendre_l^m(x)=\Legendre_l(\frac{x-x_m}{\dx/2})\}_{l=0}^{\pd-1}$ is the basis of $P^\pd(I_m)$, with  $\Legendre_l(y)$ being the $l$-th order Legendre polynomial on $[-1,1]$, normalized via $\Legendre_l(1) = 1$. The unknown coefficients are collected into vectors,
\begin{equation}
\label{eq:fourier:1}
\boldsymbol{\macro}_m^n  = [\macro_{m,0}^n, \macro_{m,1}^n, \dots, \macro_{m,\pd-1}^n],\quad
% \\
\mathbf{\micro}_{q,m}^n  = [\micro_{q,m,0}^n, \micro_{q,m,1}^n, \dots, \micro_{q,m,\pd-1}^n].
\end{equation}
By taking the Fourier ansatz 
\begin{equation}
\label{eq:fourier:ansatz}
\boldsymbol{\macro}_m^n = \widehat{\boldsymbol{\macro}}^n \exp({\imag \kappa x_m}), \quad \mathbf{\micro}_{q,m}^n = \widehat{\mathbf{\micro}}_q^n \exp({\imag \kappa x_m}) 
\end{equation}
with the wave number $\kappa$ and {the imaginary unit $\imag$ (namely, $\imag^2 = -1$)},  
the proposed IMEX-BDF$\IMEXBDForder$-DG$\pd$ scheme will lead to 
\begin{equation}
\label{eq:fourier:2}
\mathbf{W}^{n+\IMEXBDForder} = \mathbf{G}^{(\IMEXBDForder,\pd)}(\vareps,\scatlower,\dx,\dt;\xi)\mathbf{W}^{n+\IMEXBDForder-1}.   
\end{equation}
Here $\mathbf{G}^{(\IMEXBDForder,\pd)}\in {\mathbb R}^{N_f\times N_f}$ is the amplification matrix, also depending on the IMEX  partitioning, with $N_f=\IMEXBDForder\cdot r\cdot(N_v+1)$;  
$\mathbf{W}^{n+\IMEXBDForder}\in {\mathbb R}^{N_f}$, and it consists of $\widehat{\boldsymbol{\macro}}^m,  \widehat{\mathbf{\micro}}_q^m, q=1,\dots, N_v$,  with $m=n+s, n+s-1,\dots,n+1$, {to be defined in Appendix \ref{sec:stab:fourier:proof}.} The following principle will be used to study the stability:

%----------------------------------------------------------------------------------------------------------------
% Fourier analysis: Principle of numerical stability
%----------------------------------------------------------------------------------------------------------------
\begin{itemize}
\item[]{\bf Principle for Numerical Stability.} 
For any $\vareps$, $\scatlower$, $\dx$, $\dt$, let the eigenvalues of the amplification matrix $\mathbf{G}^{(\IMEXBDForder,\pd)}(\vareps,\scatlower,\dx,\dt;\xi)$ for the IMEX-BDF$\IMEXBDForder$-DG$\pd$ scheme be $\{\lambda_i(\xi)\}_{i=1}^{N_f}$. The scheme is stable if for all $\xi \in [-\pi,\pi]$, there holds
$$\max_{1\leq i\leq N_f} \left| \lambda_i(\xi) \right| \leq 1.$$  
\end{itemize}

The principle is related to that in \cite{peng2020stability,peng2021asymptotic}.    The associated analysis provides some  mathematical insights regarding the stability  of the schemes. Particularly, given that the analysis only examines  the eigenvalues of  the amplification matrix 
 $\mathbf{G}^{(\IMEXBDForder,\pd)}$, {it provides the {\it necessary} time step conditions  for  $\mathbf{G}^{(\IMEXBDForder,\pd)}$ to be power bounded and thus for the monotonicity stability (in the discrete $L^2$ energy) of the solution.}
  In actual implementations, these time step conditions are observed to be  good choices for our numerical experiments.\footnote{Readers should be cautioned that the principle for numerical stability considered here in general  allows $\mathbf{G}^{(\IMEXBDForder,\pd)}$ to have
defective eigenvalues of magnitude 1, though such pathological cases do not appear to arise for our proposed methods.}

Before reporting the stability results, we will state a theorem that reveals an important structure of the amplification matrix $\mathbf{G}^{(\IMEXBDForder,\pd)}$ in terms of its dependence on the parameters $\vareps,\scatlower,\dx,\dt$, a finding similar to that in  \cite{peng2020stability, peng2021asymptotic}. 
{The existence of such structure is valuable to mitigate the complication in  quantifying the depedence of the stability on these parameters.  The proof is in Appendix \ref{sec:stab:fourier:proof}.}

%----------------------------------------------------------------------------------------------------------------
% Fourier analysis: GL, GR similarity matrix theorem
%----------------------------------------------------------------------------------------------------------------
\begin{theorem}[{\bf Structure of the amplification matrix}] 
{Assume $\absorp=0$.} For any $\IMEXBDForder$, $\pd \geq 1$, the amplification matrix $\mathbf{G}^{(\IMEXBDForder,\pd)}(\vareps,\scatlower,\dx,\dt;\xi)$ of the IMEX-BDF$\IMEXBDForder$-DG$\pd$ scheme is similar to another matrix $\mathbf{\hat{G}}^{(\IMEXBDForder,\pd)}(\frac{\vareps}{\scatlower \dx},\frac{\dt}{\vareps \dx};\xi)$.
This indicates that the eigenvalues of $\mathbf{G}^{(\IMEXBDForder,\pd)}$ depend on $\vareps$, $\scatlower$, $\dx$, $\dt$ only through $\frac{\vareps}{\scatlower \dx}$ and $\frac{\dt}{\vareps \dx}$.
\label{theorem:similar matrix}
\end{theorem}

\medskip
\noindent
\textbf{Fourier analysis results and discussions:}
From Theorem \ref{theorem:similar matrix} and the adopted numerical principle of stability, we know that the stability of the proposed schemes will be determined by the physical and discretization parameters  $\vareps$, $\scatlower$, $\dx$, $\dt$ only through two quantities $\alpha=\vareps/(\scatlower \dx)$ and $\beta=\dt/(\vareps \dx)$. Inspired by  the energy-based stability analysis for  the IMEX-BDF1-DG1 schemes in Section \ref{sec:stab:energy}
 and the findings in \cite{peng2020stability,peng2021asymptotic}, 
we proceed to identify the stability regions in the  $\alpha-\beta$  plane, associated with the  adopted stability principle, for each of the IMEX-BDF$\pd$-DG$\pd$  schemes with $\pd = 1, 2, 3$ and all  three IMEX  partitionings. It turns out the stability regions are separated from the unstable regions by some $\alpha-\beta$ curve, which will be found numerically as follows:  we sample $\alpha$ uniformly in a logarithmic scale, particularly, by taking $\log_{10}\alpha=-5+j/20\in[-5,5], j=0,\cdots,200$, and find   the respective $\beta$  on the interface curve between the stable and unstable regions via the standard bisection process. 
The variable $\xi$ is discretized uniformly over $[-\pi,\pi]$ with a spacing of $\pi/50$.
Using the logarithmic scale in both $\alpha$ and $\beta$, the stability results are plotted  for the IMEX-BDF$\pd$-DG$\pd$-$\mS$$k$ schemes with $\pd, k = 1, 2, 3$ in Figure \ref{fig:fourier:st123}. The main observations are summarized below.

\begin{enumerate}
\item[1.)] Similar to the findings in Remark \ref{rem:stab:energy:1} for the first order methods based on energy-based stability analysis,   the proposed schemes with  Strategy 1 or 2 are conditionally stable in all regimes, while with Strategy 3, the proposed schemes are unconditionally stable in the scattering dominant regime.
\item[2.)] More specifically, for each IMEX-BDF$\pd$-DG$\pd$-$\mS$$k$ scheme ($r,k=1,2,3$), there exist some positive constants $\alpha_{r,k}, \beta_{r,k}, \gamma_{r,k}$ and a constant $\eta_{r,k}$, such that 
\begin{itemize}
\item[2.a)] when $\alpha=\vareps/(\scatlower \dx) > \alpha_{r,k}$, the boundary between the stable and unstable regions is a straight line $\beta = \beta_{r,k}$. This implies when $\vareps/(\scatlower \dx)=O(1)$ and the numerical model (i.e. equation + scheme) is in the transport dominant regime, the scheme is (conditionally) stable under the condition, $\beta=\dt/(\vareps \dx)\leq \beta_{r,k}$, that is, under a hyperbolic time step condition $\dt \leq \beta_{r,k}\vareps \dx$.

\item [2.b)] When $\alpha=\vareps/(\scatlower \dx) \leq \alpha_{r,k}$, then
\begin{itemize}
\item [-] for the schemes with  Strategy 1 or 2,
the boundary between the stable and unstable regions (in a logarithmic scale) can be  bounded  below by a  line with a slope $-1$, namely $\log_{10}\beta = -\log_{10}\alpha + \eta_{r,k}$, or equivalently, $\alpha\beta=\dt/(\scatlower h^2)=\gamma_{r,k}$. In other words, when $\vareps/(\scatlower \dx)\ll 1$ and the numerical model  is in the scattering dominant regime, the scheme is (conditionally) stable under a parabolic time step condition $\dt \leq \gamma_{r,k}\scatlower \dx^2$.
    \item[-] for the scheme with  Strategy 3, 
    it is stable for any $\beta$. This implies that when $\vareps/(\scatlower \dx)\ll 1$ and the numerical model  is in the scattering dominant regime,  the scheme is unconditionally stable. 
    \end{itemize}
        \end{itemize}
    \end{enumerate}

{From the stability plots, we further extract numerically the explicit formulas   for the time step conditions,  $\dt\leq\dt_{CFLr}^{(k)}$.  These formulas are in similar ansatzes as in Theorem \ref{thm:stab:energy}, except that with Strategy 3, we set the time step condition to be $\Delta t = 0.5h$ when the methods are unconditionally stable.  These formulas  will be used in  Section \ref{sec:num}, 
unless otherwise specified. }

%----------------------------------------------------------------------------------------------------------------
% Strategy 1: Time step condition
%----------------------------------------------------------------------------------------------------------------
\begin{center}
\textbf{Strategy 1, 2} ($k=1,2$) 
\end{center}
\begin{subequations}
\label{eq:fourier:dt:st12}
\begin{align}
\text{IMEX-BDF1-DG1-}\mS k: \ \dt_{CFL1}^{(k)} & = 0.7 \vareps \dx + 0.62 \scatlower\dx^2, \\
\text{IMEX-BDF2-DG2-}\mS k : \ \dt_{CFL2}^{(k)} & = 0.151 \vareps \dx + 0.027 \scatlower\dx^2, \\
\text{IMEX-BDF3-DG3-}\mS k: \ \dt_{CFL3}^{(k)} & = 0.062 \vareps \dx + 0.0021 \scatlower\dx^2.
\end{align}
\end{subequations}
%---------------------------------------------------------------------------------------------------------------------------------
% Strategy 3: Time step condition
%----------------------------------------------------------------------------------------------------------------
\begin{center}
\textbf{Strategy 3}
\end{center}
\begin{subequations}
\label{eq:fourier:dt:st3}
\begin{align}
\text{IMEX-BDF1-DG1-}\mS3:
\ \dt_{CFL1}^{(3)} 
& =
\begin{cases}
0.5 \dx, & \vareps \leq 0.5 \scatlower \dx \\
\min\left(0.5 \dx, \frac{\vareps^2 \dx}{\vareps - 0.5 \scatlower \dx}\right), & \vareps > 0.5 \scatlower \dx
\end{cases},
\\
\text{IMEX-BDF2-DG2-}\mS3: 
\ \dt_{CFL2}^{(3)} 
& = 
\begin{cases}
0.5 \dx, & \vareps \leq 0.056 \scatlower \dx \\
\min\left(0.5 \dx, \frac{0.2 \vareps^2 \dx}{\vareps - 0.056 \scatlower \dx}\right), & \vareps >0.056 \scatlower \dx    
\end{cases},
\\
\text{IMEX-BDF3-DG3-}\mS3:
\ \dt_{CFL3}^{(3)} 
& = 
\begin{cases}
0.5 \dx, & \vareps \leq 0.01 \scatlower \dx \\
\min\left(0.5 \dx, \frac{0.07 \vareps^2 \dx}{\vareps - 0.01 \scatlower \dx}\right), & \vareps >0.01 \scatlower \dx    
\end{cases}.
\end{align}
\end{subequations}

\begin{figure}[h!]
\centering
\begin{subfigure}{0.32\textwidth}
\includegraphics[width=\linewidth,height=1.8in]
{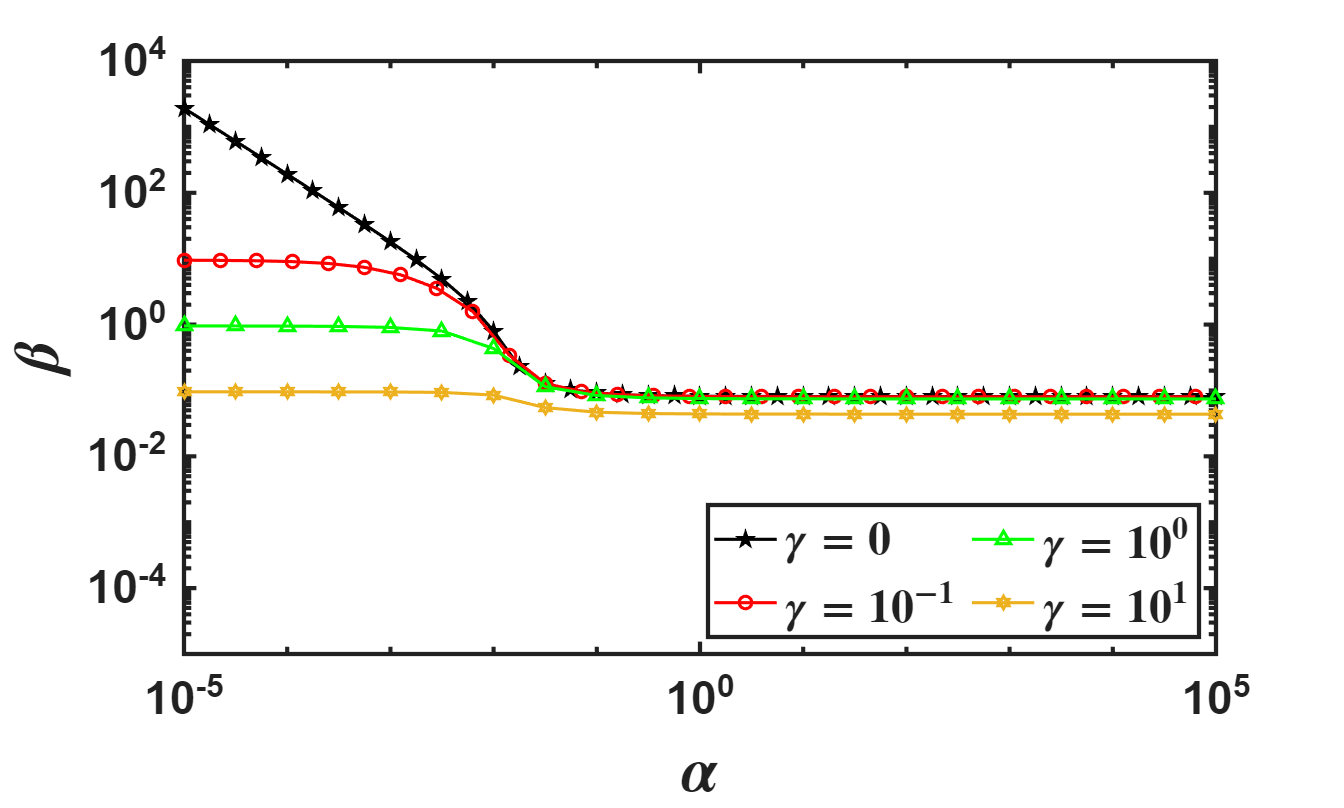}
\caption{Strategy 1, explicit $\sigma_a$}
\end{subfigure}
\begin{subfigure}{0.32\textwidth}
\includegraphics[width=\linewidth,height=1.8in]
{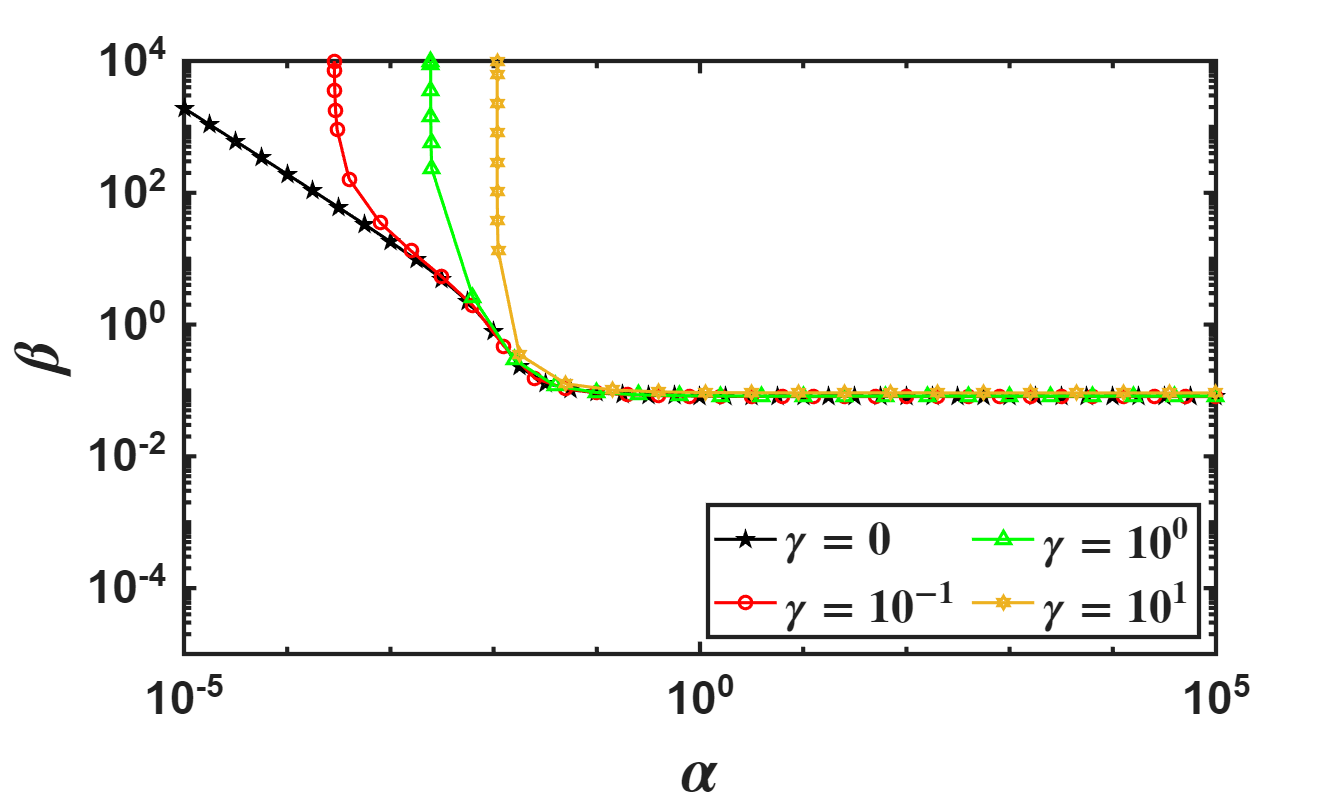}
\caption{Strategy 1, implicit $\sigma_a$}
\end{subfigure}
\begin{subfigure}{0.32\textwidth}
\includegraphics[width=\linewidth,height=1.8in]
{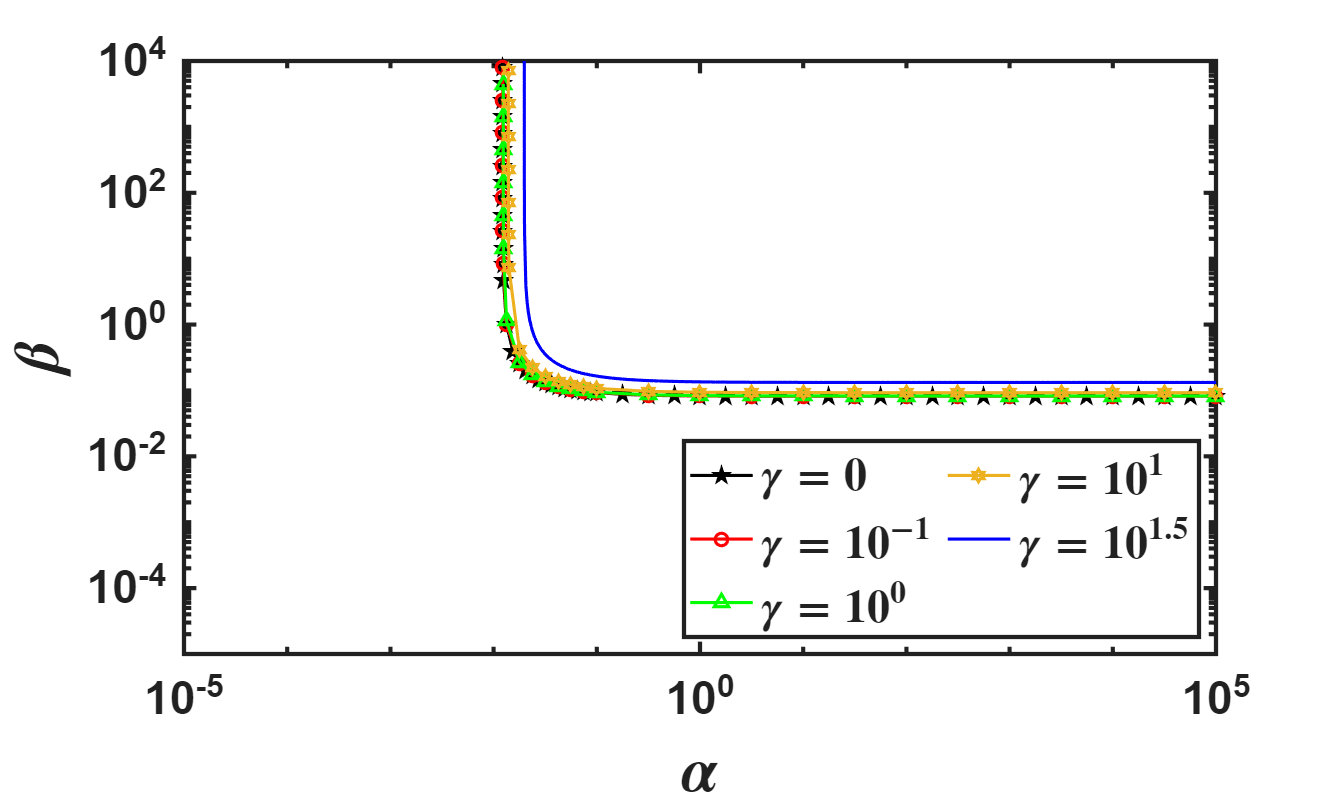}
\caption{Strategy 3}
\end{subfigure}
\caption{Stability curves of IMEX-BDF3-DG3, with the stable region below each curve. (a): with Strategy 1 and explicit treatment of the absorption terms; (b): with Strategy 1 and implicit treatment of the absorption terms; (c): with Strategy 3. $\alpha=\vareps/(\scatlower \dx)$, $\beta=\dt/(\vareps \dx)$, and $\gamma=\vareps\absorpupper h$. }
\label{fig:Fourier:absorp:3rdOrder}
\end{figure}

\medskip
\noindent
\textbf{Extension to the analysis with nonzero  $\absorp$:} Fourier analysis can be further extended to the case in the presence of the absorption terms when $\absorp(x)=\absorpupper>0$. {Similar to Theorem \ref{theorem:similar matrix}, one can show  that the stability, in the sense of the Principle for Numerical Stability stated earlier in this section,  will depend on  $\vareps$, $\scatlower$, $\absorpupper$, $\dx$, $\dt$ only through $\alpha=\vareps/(\scatlower \dx)$, $\beta=\dt/(\vareps \dx)$, and $\gamma=\vareps\absorpupper h$. Readers are referred to Remark 2.3.5 in \cite{Matsuda2025RPI} for the proof.} Given that the observations are 
qualitatively similar (not  quantitatively though), we only  present selectively three stability plots in Figure \ref{fig:Fourier:absorp:3rdOrder}, all for the third order IMEX-BDF3-DG3 schemes. There is no difference in the plots for Strategy 1 and Strategy 2. 
Again, what is being plotted are the $\alpha-\beta$ stability curves between the  stable and unstable regions,  with the stable  regions located below the curves,  {for representative values of $\gamma$ (recall $\sigma_a=O(1)$).}

{Note that with the model (1.1) in the dimensionless form under a diffusive scaling, one assumes $\sigma_s=O(1)$ and $\sigma_a=O(1)$, and  the strength of the scattering  is about $\sigma_s/\varepsilon$ and the strength of the absorption is about $\varepsilon\sigma_a$. To facilitate the discussion on stability based on Figure \ref{fig:Fourier:absorp:3rdOrder} in the presence of many parameters, $\sigma_{s,m}$ and $\alpha=\vareps/(\scatlower \dx)$ are assumed to be fixed. It is observed, as expected,  that with the explicit (resp. implicit) treatment of the absorption terms, the {region of stability} decreases (resp. increases) as $\gamma$ grows, e.g., when $\absorpupper$ increases, or when $\vareps$ and $h$ increase at the same rates. It is somewhat unexpected to observe that  when  the absorption terms are treated implicitly with  Strategy 1 (and 2), and  when  $\alpha\ll 1$ (i.e., the spatial meshes are under-resolved), the stability gets closer to being unconditional with relatively larger $\gamma$. }

\subsection{Further numerical validation and comparison}
\label{sec:stab:numV}

{In this subsection, we will perform some numerical tests to further assess the time step conditions predicted by the stability analysis. With qualitative similarity, we only present the results for first order methods. (Note that time step conditions by energy-based stability are only available for the first order methods in this work.) To this end, we consider the one-group transport equation in slab geometry on $\spaceset=[0,2\pi]$ with the material properties $\scat(x)=1$ and $\absorp(x)=0$ and zero source.
 The solution is smooth, with the initial data $\macro(x,0) = \sin(x)$ and $\micro(x,\velangle,0) = -\velangle \cos(x)$, and periodic boundary conditions.   Using this example, we implement each scheme on a uniform mesh of $N_x=N$ elements, with $N=20, 40, \cdots, 1280$, up to the final time $T=100$, when $\vareps=0.5, 10^{-2}, 10^{-6}$. We  numerically compute the maximal time step $\Delta t_\star$ that ensures the monotonicity $L^2$ stability  
for each scheme. This will be done via  
 a bisection search: we start with an interval $(\tau_1, \tau_2)=(0,1)$. In the first iteration,  the time step is taken as $\Delta t=(\tau_1+\tau_2)/2$ and we 
implement the method.  If the squared $L^2$ energy $E_h^n=|| \macro_h^n ||^2 + \vareps^2 ||| \micro_h^n |||^2$  of the numerical solution decreases in time\footnote{The monotonicity $L^2$ stability is examined when $E_h^n=E_{h,\mS1}^n$ is taken as the discrete energy for each scheme. 
There is no essential difference observed and concluded for Strategy 2 if $E_h^n=E_{h,\mS2}^n$ is used, and for Strategy 3 if $E_h^n=E_{h,\mS3}^n$ is used.},
in the sense that  
    $E_n^{n+1}-E_n^{n} \le  10^{-24}$
    for all $t^{n+1}\leq T$, we set $\tau_1=\Delta t$; otherwise, $\tau_2=\Delta t$. The search ends when  $|\tau_1-\tau_2|\le tol=10^{-6}$, and we set $\Delta t_\star=\Delta t$. }
    
\begin{figure}[ht]
\centering
\begin{subfigure}{0.32\textwidth}
\includegraphics[width=\linewidth,height=1.8in]
{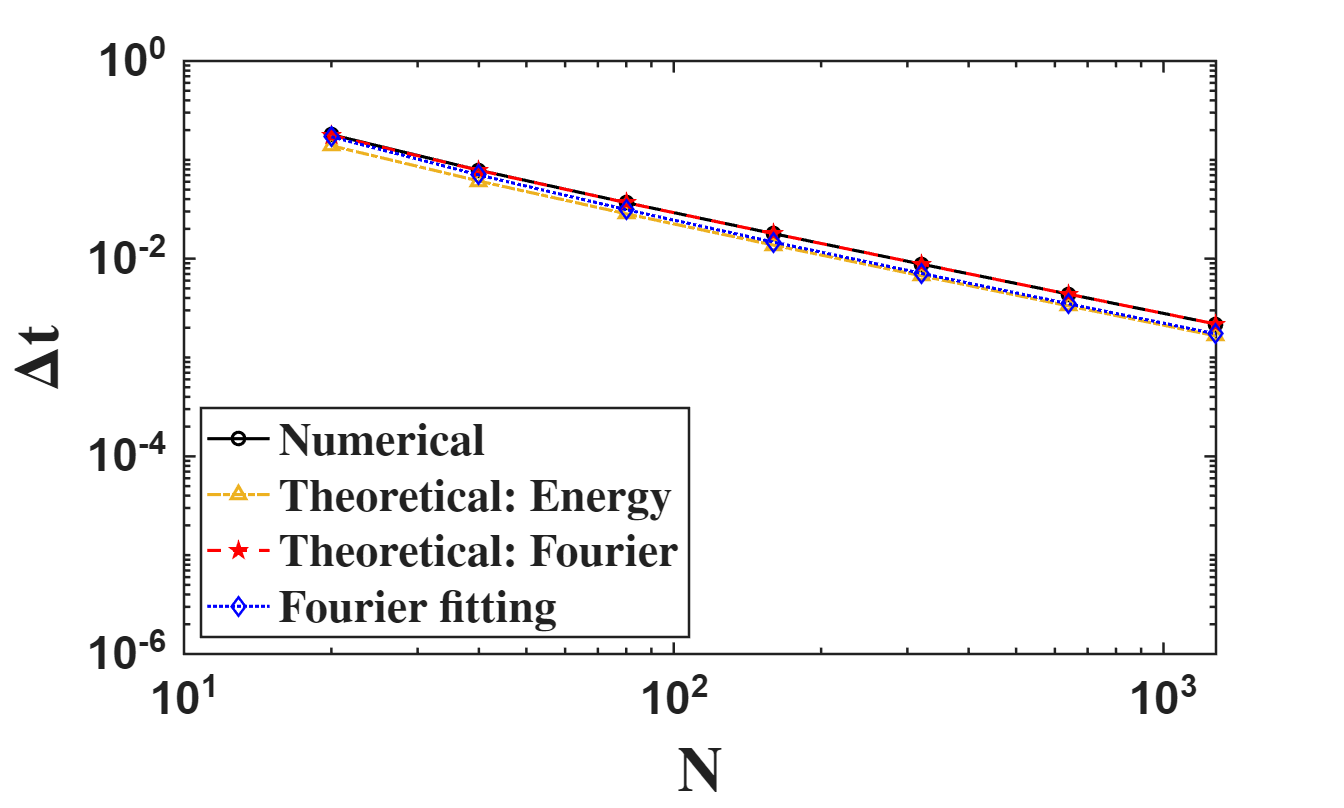}
\caption{$\vareps = 0.5$}
\end{subfigure}
\begin{subfigure}{0.32\textwidth}
\includegraphics[width=\linewidth,height=1.8in]{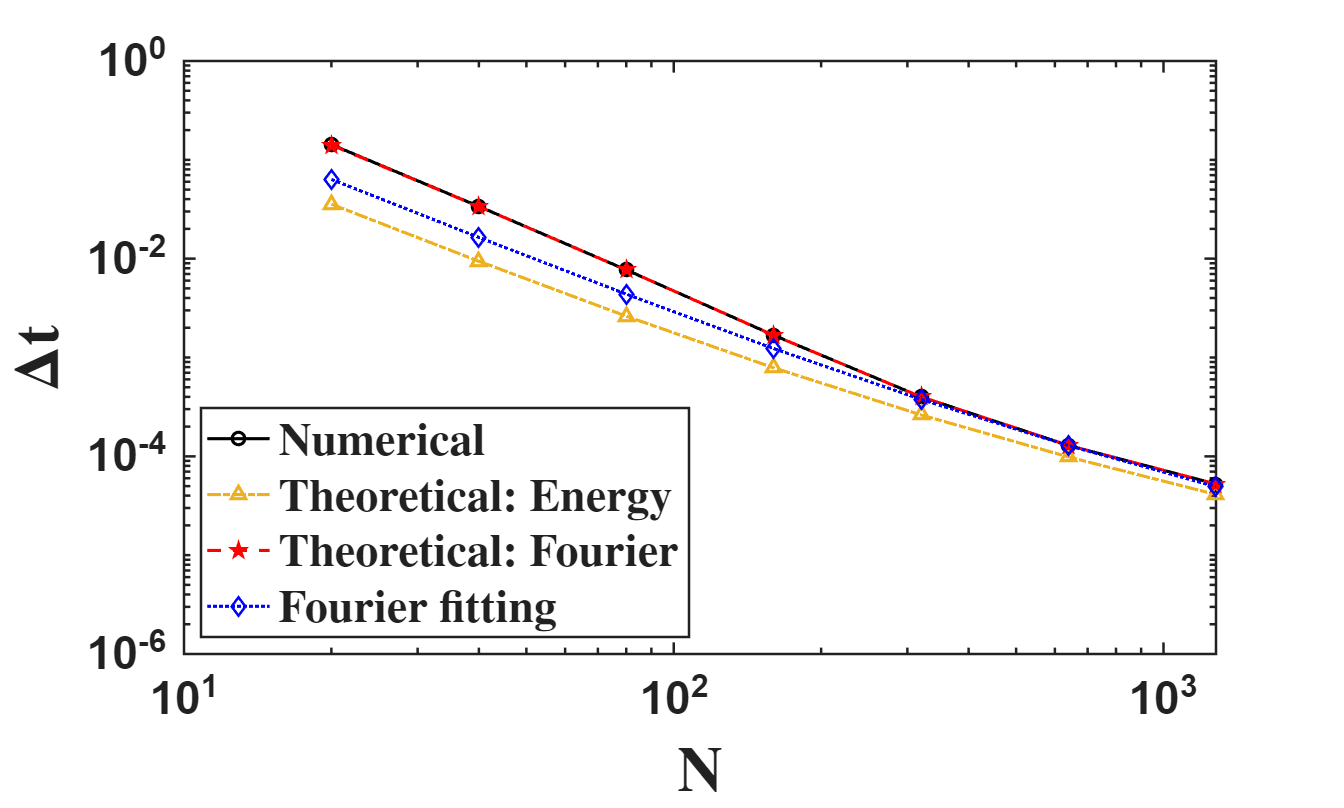}
\caption{$\vareps = 10^{-2}$}
\end{subfigure}
\begin{subfigure}{0.32\textwidth}
\includegraphics[width=\linewidth,height=1.8in]{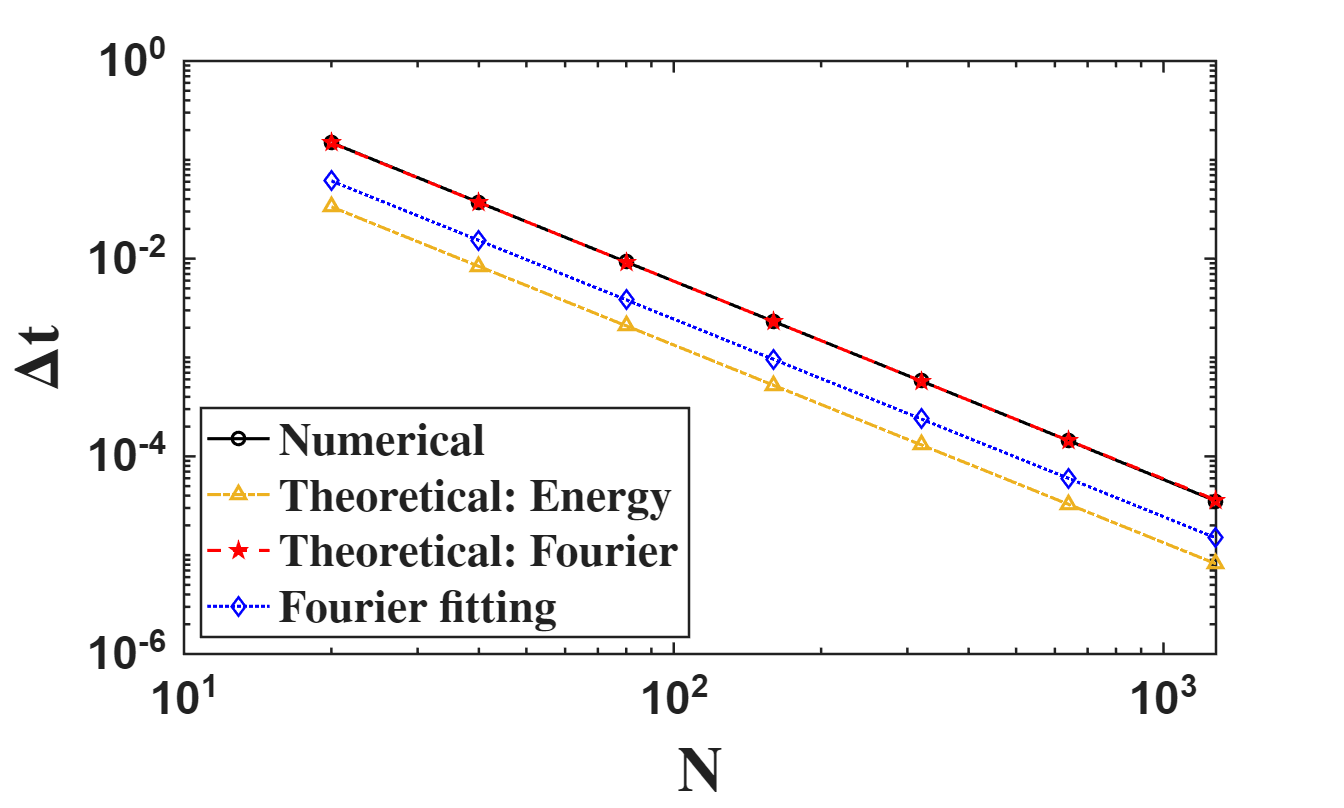}
\caption{$\vareps = 10^{-6}$}
\end{subfigure}
\caption{Maximum $\Delta t$ allowed for stability of IMEX-BDF1-DG1-$\mathcal{S}1$.}
\label{fig:stability:IMEXBDF1DG1:st1}
\end{figure}

\begin{figure}[ht]
\begin{subfigure}{0.32\textwidth}
\includegraphics[width=\linewidth,height=1.8in]
{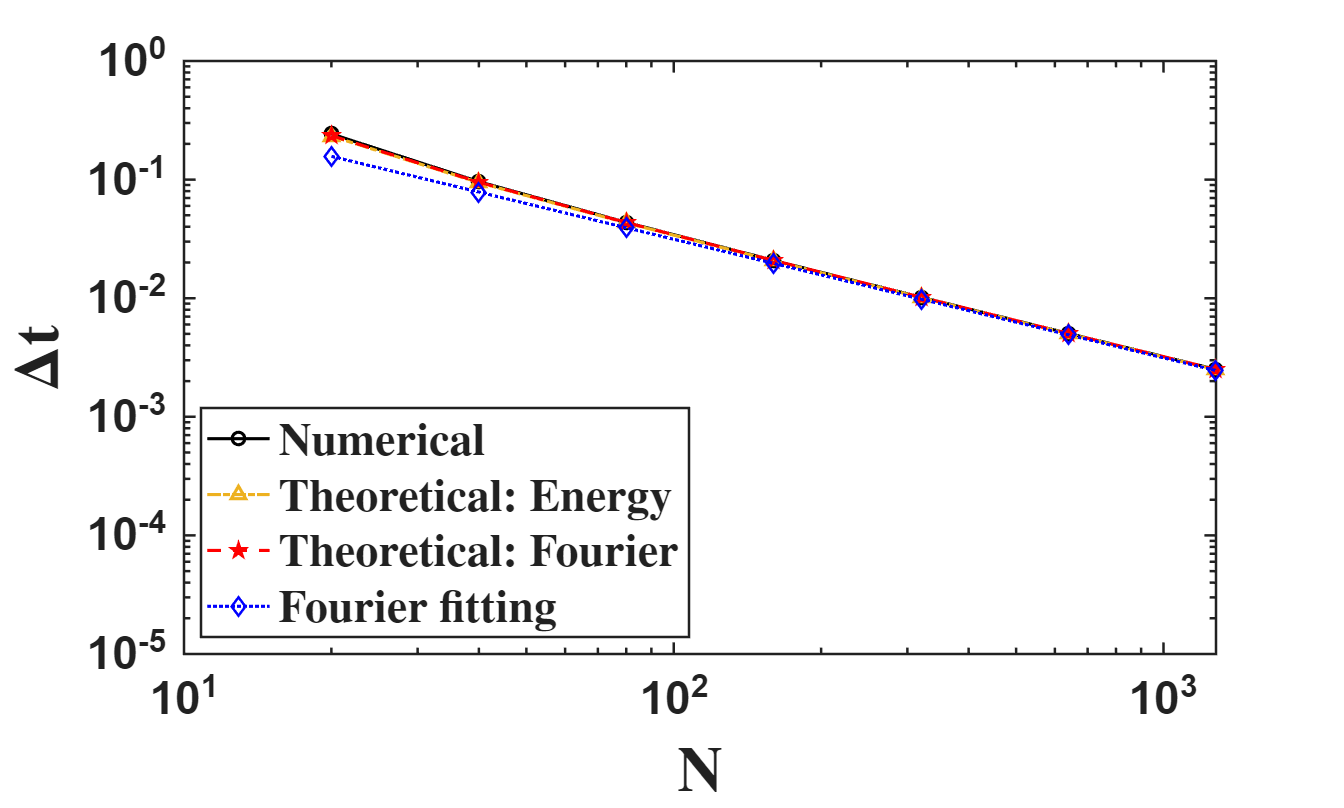}
\caption{$\vareps = 0.5$}
\end{subfigure}
\begin{subfigure}{0.32\textwidth}
\includegraphics[width=\linewidth,height=1.8in]{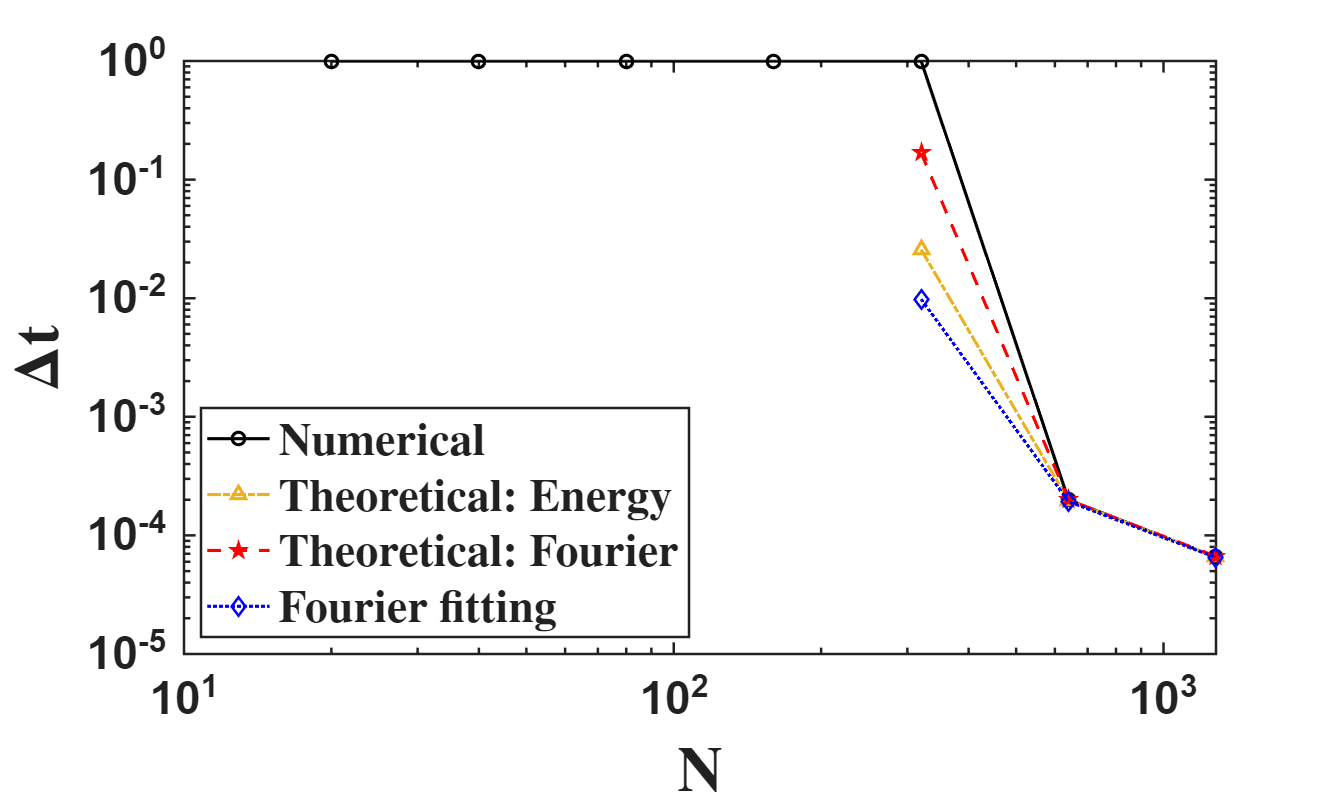}
\caption{$\vareps = 10^{-2}$}
\end{subfigure}
\begin{subfigure}{0.32\textwidth}
\includegraphics[width=\linewidth,height=1.8in]{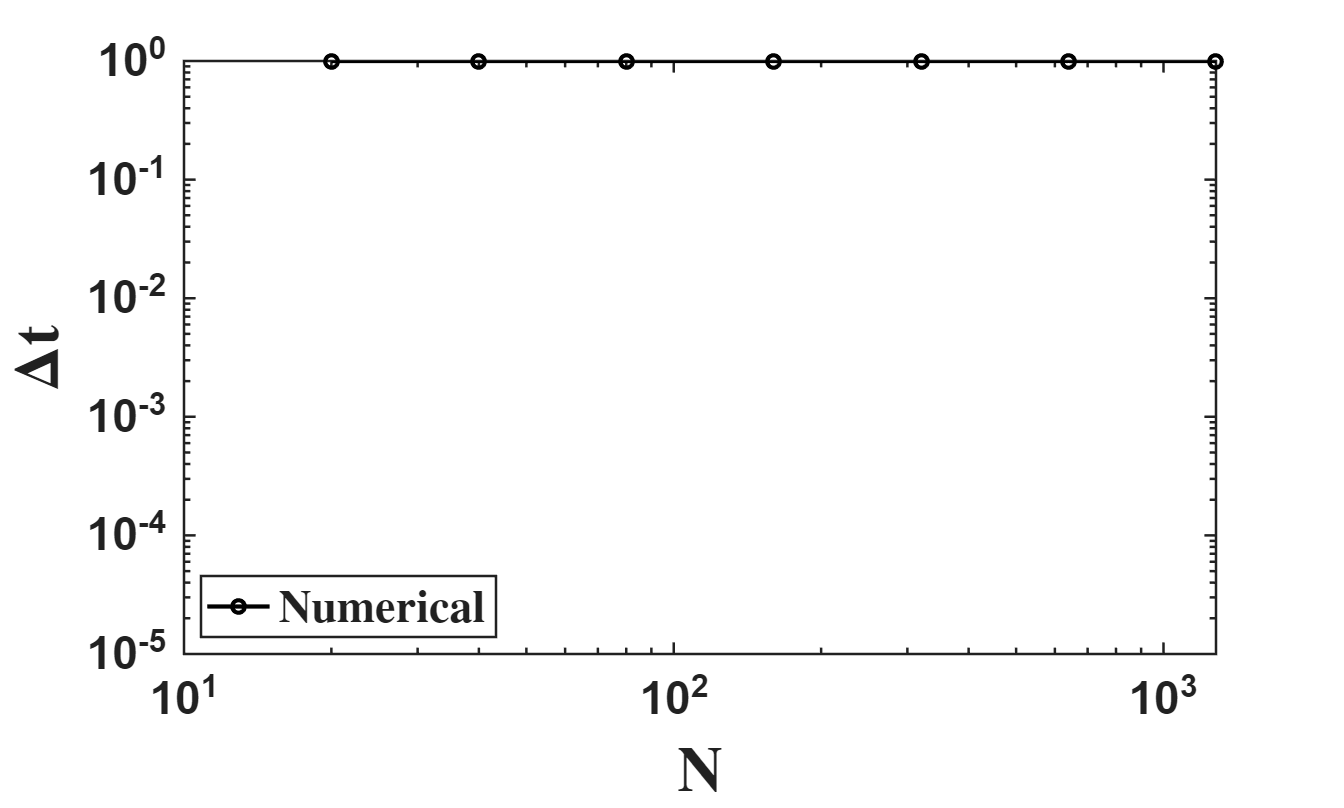}
\caption{$\vareps = 10^{-6}$}
\end{subfigure}

\caption{Maximum $\Delta t$ allowed for stability of IMEX-BDF1-DG1-$\mathcal{S}3$.}
\label{fig:stability:IMEXBDF1DG1:st3}
\end{figure}

{In  Figures \ref{fig:stability:IMEXBDF1DG1:st1}-\ref{fig:stability:IMEXBDF1DG1:st3},
we present the computed maximal time step   
$\Delta t_\star$, labeled as {\it ``Numerical"},  and compare it with the maximal time step sizes predicted by energy-based stability ({\it Theoretical: Energy}), Fourier-based stability ({\it Theoretical: Fourier}), and the fitted formulas in \eqref{eq:fourier:dt:st12}-\eqref{eq:fourier:dt:st3} based on Fourier analysis ({\it Fourier fitting}).  They are for IMEX-BDF1-DG1-$\mathcal{S}1$ and IMEX-BDF1-DG1-$\mathcal{S}3$, respectively. }

{Let us start with  Figure \ref{fig:stability:IMEXBDF1DG1:st1} for IMEX-BDF1-DG1-$\mathcal{S}1$. In all three regimes considered here with $\vareps=0.5, 10^{-2}, 10^{-6}$,  it is observed that the computed $\Delta t_\star$ ({\sf black $\circ$}) is always on the top as the largest , and it is almost indistinguishable from that by Fourier analysis ({\sf red $\star$}).  Followed next are the maximal time steps by fitted formulas based on Fourier analysis ({\sf blue $\diamond$}) and finally that by energy-based stability analysis ({\sf yellow $\triangle$}). Symbolically we write 
\begin{center}
{\sf (black $\circ$)$\approx$ (red $\star$)$\ge$(blue $\diamond$)$\ge$(yellow $\triangle$)}.
\end{center}
Let us further comment on the observations. What is not obvious but good news is (black $\circ$)$\approx$ (red $\star$), evidencing that the time step conditions predicted by the Fourier-based analysis, though generally only necessary for (monotonicity) stability, seem to also be  sufficient (at least) in the setting  examined here. Other observations are largely expected: (1) the numerically computed $\Delta t_\star$ is the largest; (2) (red $\star$)$\ge$(blue $\diamond$), due to the numerical fitting of the time step conditions by Fourier analysis into a specific ansatz $\Delta t= c_1 \vareps h + c_2\scatlower h^2$ in {\it all} regimes; 
(3) To ensure the monotonicity of the $L^2$ energy, the time step conditions by the energy-based stability analysis are sufficient while those by the Fourier-based stability analysis in general are necessary. This explains (red $\star$)$\ge$(yellow $\triangle$). The results by IMEX-BDF1-DG1-$\mathcal{S}2$ are not presented here as they are indistinguishable to those by  IMEX-BDF1-DG1-$\mathcal{S}1$.}

{Figure \ref{fig:stability:IMEXBDF1DG1:st3} is for IMEX-BDF1-DG1-$\mathcal{S}3$. Similar observations as in Figure  \ref{fig:stability:IMEXBDF1DG1:st1} can be made when the method is conditionally stable, and this corresponds to Figure \ref{fig:stability:IMEXBDF1DG1:st3}-(a)  as well as  Figure \ref{fig:stability:IMEXBDF1DG1:st3}-(b) when $N\geq 320>2\pi/(2\|v\|_{h,\infty}\vareps)$. Related to Figure \ref{fig:stability:IMEXBDF1DG1:st3}-(c) and the rest of  Figure \ref{fig:stability:IMEXBDF1DG1:st3}-(b), the method is unconditionally stable, and $\Delta t_\star\approx 1$ due to the initial interval $(0,1)$ used in the bisection search. }

%----------------------------------------------------------------------------------------------------------------

\section{Computational complexity: algebraic form, conditioning related to Strategy 3}
\label{sec:alg}

In this section, we present the matrix-vector forms of our IMEX-BDF$\IMEXBDForder$-DG$\pd$-$\mS k$ schemes,  to examine  the computational complexity when each IMEX  partitioning is applied.  With Strategy 3 (i.e. $k=3$), we will also estimate the conditioning of the linear system that needs to be solved over each time step. 

Following the same notation in \eqref{eq:fourier:0}-\eqref{eq:fourier:1} for the local basis  and solution expansions associated with each {mesh element $I_m$ from Section \ref{sec:fourier}, we define $\Phi_l^m(x)\in U_h^{r-1}$} as $\Phi_l^m(x)|_{I_n}=\delta_{mn}\phi_l^m(x)$,
and write 
$$\mathbf{\pdspacebasis} = [\Phi_0^1,\Phi_1^1, \dots, \Phi_{\pd-1}^1,\Phi_0^2,\Phi_1^2, \dots, \Phi_{\pd-1}^2, \dots,   \Phi_0^{N_x},\Phi_1^{N_x}, \dots, \Phi_{\pd-1}^{N_x}]^T.$$
The $i$th entry of $\mathbf{\pdspacebasis}$ will also be referred to as $\pdspacebasis_i$ when needed. In addition, we write 
$\boldsymbol{\macro}^n = [\boldsymbol{\macro}^n_1, \boldsymbol{\macro}^n_2, \cdots, \boldsymbol{\macro}^n_{N_x}]^T$, 
$\mathbf{\micro}_q^n = [\mathbf{\micro}_{q,1}^n, \mathbf{\micro}_{q,2}^n, \dots, \mathbf{\micro}_{q,N_x}^n]^T$, $q=1,\dots, N_{\velangle}$. The numerical solutions can then be compactly expressed as 
\begin{equation}
\macro_h^n(x) =(\boldsymbol{\macro}^n)^T \mathbf{\pdspacebasis}, \qquad \micro_{h,q}^n(x)  = (\mathbf{\micro}_q^n)^T \mathbf{\pdspacebasis}, \;\; q=1,\dots, N_{\velangle}.
\end{equation}

We further introduce matrices of size $N_{rx}\times N_{rx}$ (with $N_{rx}=\pd N_x$) associated with different terms in the numerical methods. They are the mass matrix $\M=(\M_{ij})$ with $\M_{ij} = (\pdspacebasis_j,\pdspacebasis_i)$, advection matrices $D^\pm=(D^\pm_{ij})$ with  $\Dplus_{ij} = (\Dplusop \pdspacebasis_j,\pdspacebasis_i)$, $\Dminus_{ij} = (\Dminusop \pdspacebasis_j,\pdspacebasis_i)$, the scattering matrix
$\scatmat=\big((\scatmat)_{ij}\big)$ with $(\scatmat)_{ij} = (\scat \pdspacebasis_j,\pdspacebasis_i)$, and the absorption matrix $\absorpmat=\big((\absorpmat)_{ij}\big)$ with   $(\absorpmat)_{ij} = (\absorp \pdspacebasis_j,\pdspacebasis_i)$.

For the methods with Strategy $k$, $k = 1, 2, 3$, we have the following matrix-vector form 
\begin{equation}
\mathcal{L}^{(k)}(\boldsymbol{\macro}^{n+s}, \mathbf{\micro}_1^{n+s}, \mathbf{\micro}_2^{n+s}, \dots, \mathbf{\micro}_{N_{\velangle}}^{n+s})^T =  (\mathbf{b}^{n,(k)}_{0}, \mathbf{b}_{1}^{n,(k)}, \mathbf{b}_{2}^{n,(k)}, \dots, \mathbf{b}_{N_{\velangle}}^{n,(k)})^T
\end{equation}
where
%----------------------------------------------------------------------------------------------------------------
% Strategy 1, previous 2
%
\begin{equation} 
\mathcal{L}^{(1)}
 =
\begin{bmatrix}
\Theta^{(1)} & \dt c_\IMEXBDForder\vweight_1 \velangle_1 \Dplus & \dt c_\IMEXBDForder \vweight_2 \velangle_2 \Dplus & \dots & \dt c_\IMEXBDForder \vweight_{N_{\velangle}} \velangle_{N_{\velangle}} \Dplus \\
0 & \Theta^{(2)} & 0 & \dots & 0 \\
0 & 0 & \Theta^{(2)} & \dots & 0 \\
\vdots & \vdots & \vdots & \ddots & \vdots \\
0 & 0 & 0 & \dots & \Theta^{(2)}
\end{bmatrix},
\label{eq:alg-S1}
\end{equation}

%

% Strategy 2, previous 1
%----------------------------------------------------------------------------------------------------------------
\begin{equation}
\mathcal{L}^{(2)} =
\begin{bmatrix}
\Theta^{(1)} & 0 & 0 & \dots & 0 \\
\velangle_1 \dt c_\IMEXBDForder\Dminus & \Theta^{(2)} & 0 & \dots & 0 \\
\velangle_2 \dt c_\IMEXBDForder\Dminus & 0 & \Theta^{(2)} & \dots & 0 \\
\vdots & \vdots & \vdots & \ddots & \vdots \\
\velangle_{N_{\velangle}} \dt c_\IMEXBDForder\Dminus & 0 & 0 & \dots & \Theta^{(2)}
\end{bmatrix},
\label{eq:alg-S2}
\end{equation}
%----------------------------------------------------------------------------------------------------------------
% Strategy 3
%----------------------------------------------------------------------------------------------------------------
\begin{equation}
\mathcal{L}^{(3)}
 =
\begin{bmatrix}
\M + \dt c_\IMEXBDForder\absorpmat & \dt c_\IMEXBDForder\vweight_1 \velangle_1 \Dplus & \dt c_\IMEXBDForder\vweight_2 \velangle_2 \Dplus & \dots & \dt c_\IMEXBDForder\vweight_{N_{\velangle}} \velangle_{N_{\velangle}} \Dplus \\
\velangle_1 \dt c_\IMEXBDForder\Dminus & \Theta^{(3)} & 0 & \dots & 0 \\
\velangle_2 \dt c_\IMEXBDForder\Dminus & 0 & \Theta^{(3)} & \dots & 0 \\
\vdots & \vdots & \vdots & \ddots & \vdots \\
\velangle_{N_{\velangle}} \dt c_\IMEXBDForder\Dminus & 0 & 0 & \dots & \Theta^{(3)}
\end{bmatrix}
\label{eq:alg-S3}
\end{equation}
and
\begin{equation}
\Theta^{(1)} = \M, \qquad \Theta^{(2)}  = \vareps^2 \M + \dt c_\IMEXBDForder \scatmat, \qquad 
\Theta^{(3)} = \vareps^2 \M + \dt c_\IMEXBDForder (\vareps^2 \absorpmat+\scatmat). \label{eq:alg-4}
\end{equation}
The right hand side terms $\{\mathbf{b}_{q}^{n,(k)}\}_{q=0}^{N_{\velangle}}$  depend on the numerical solutions available from $t^{n+j}, j=0,\dots, s-1$, the source terms (and possibly boundary data in the case of non-periodic boundary conditions).

With $\scat, \absorp$ being nonnegative along with the local nature of the space $\pdspace^{\pd}$ and the chosen basis $\{\Phi_l^m: l=0,\dots, \pd-1, m=1,\dots, N_x\}$, one can easily see  $\Theta^{(k)}$, $k=1,2,3$,  is symmetric positive definite (SPD), thus invertible, and it is also block-diagonal. We are now ready to discuss the implementation of the methods with IMEX Strategy 1 and 2 and get some idea of the computational complexity.
\begin{itemize}
    \item Using Strategy 1, we first solve for  $\{\mathbf{\micro}_q^{n+s}\}_{q=1}^{N_{\velangle}}$ by inverting $\Theta^{(2)}$ locally (i.e. element by element), in a parallel fashion with respect to  $q$ if one wants to. We then solve for $\boldsymbol{\macro}^{n+s}$ by inverting $\Theta^{(1)}$ locally.  Using Strategy 2, a similar procedure is followed, except that we switch the order to solve for $\boldsymbol{\macro}^{n+s}$ first {and} then {solve} for  $\{\mathbf{\micro}_q^{n+\IMEXBDForder}\}_{q=1}^{N_{\velangle}}$ (parallelly in $q$, if one wants to).
    \item When $\dt$ stays unchanged,  $\Theta^{(1)}$ and $\Theta^{(2)}$ are independent of time. Their inverses can be computed once, at a cost of $O(r^3N_x)$ due to their block-diagonal structures, and used  throughout the time marching process. In other words, the methods with  Strategy 1 and 2 have similar computational complexity as {\it fully explicit methods} per time step. 
\end{itemize}

We now turn to the methods with Strategy 3 (i.e. $k=3$).  Recall the methods of  first order accuracy in this family are the same as those in \cite{peng2021asymptotic}. Following the Schur complement procedure there, one can first express each $\mathbf{\micro}_q^{n+s}$ in terms of $\boldsymbol{\macro}^{n+s}$, namely,
\begin{equation}
\mathbf{\micro}_q^{n+s}=(\Theta^{(3)})^{-1} \big(\mathbf{b}_{q}^{n,(k)} - \velangle_q \dt c_\IMEXBDForder\Dminus\boldsymbol{\macro}^{n+s} \big),\quad  q=1,\dots, N_{\velangle},
\label{eq:alg-5}
\end{equation}
then substitute them into the equation for $\macro$ (corresponding to the first row block in \eqref{eq:alg-S3}), resulting in the following linear system for $\boldsymbol{\macro}^{n+s}$,
\begin{equation} \label{eq:alg-6}
\mathcal{H} \boldsymbol{\macro}^{n+s} = \mathbf{b}^{n,(k)},
\end{equation}
where the term $\mathbf{b}^{n,(k)}$ depends on the numerical solutions available from the previous $\IMEXBDForder$ steps and the source terms (and possibly boundary data in the case of non-periodic boundary conditions). The matrix $\mathcal{H}$ is given as
\begin{equation} \mathcal{H} = \M + \dt c_\IMEXBDForder\absorpmat - \lvg \velangle^2\rvg \dt^2 c_\IMEXBDForder^2 \Dplus (\Theta^{(3)})^{-1} \Dminus \in {\mathbb R}^{N_{rx}\times N_{rx}}. 
\label{eq:alg-7}
\end{equation}
Here we have used  $\lvg \velangle^2\rvg = \lvg \velangle^2\rvg _h$ in \eqref{eq:v square integral exact}. When the boundary conditions are periodic, we know $\Dplusop = -(\Dminusop)^T$ in \eqref{eq:adjoint}, and this implies  $\Dplus = -(\Dminus)^T$ and hence  $\mathcal{H}$ is SPD.
Once  $\boldsymbol{\macro}^{n+\IMEXBDForder}$ is computed from \eqref{eq:alg-6}, one can solve $\{\mathbf{\micro}_q^{n+\IMEXBDForder}\}_{q=1}^{N_{\velangle}}$ from \eqref{eq:alg-5}, in a parallel fashion in $q$ if one wants to.  This can be carried out very locally, as $\Dminus$ is sparse and block-structured.  

The next theorem examines the conditioning of $\mathcal{H}$ in terms of model and discretization parameters,  and the proof is given in Appendix \ref{sec:proof:CondH}.

\begin{theorem}[{\bf Conditioning of $\mH$}] With  periodic boundary conditions in space, and under the assumptions of a uniform mesh in space $\dx_m=\dx, \forall m$, as well as  constant scattering and absorption cross sections $\scat$ and $\absorp$, the following estimate holds for the condition number (in {the} 2-norm) of $\mH$:
\begin{equation}
\textrm{cond}_2(\mH)\leq 
(2r-1)(1+\dt c_\IMEXBDForder\absorp)+\frac{(2r-1)^2C_r \lvg \velangle^2\rvg \dt^2c_\IMEXBDForder^2 }{h^2\big(\vareps^2+\dt c_\IMEXBDForder(\vareps^2\absorp+\scat)\big)},
\label{eq:condH}
\end{equation}
with $C_r=4(\sqrt{3}+1)^2\left(\sum_{l=0}^{\pd-1} \frac{1}{2l+1}\right)r^4$.
\label{thm:condH}
\end{theorem}

\begin{remark}
    The estimate of the condition number of $\mH$ can be  extended to the case with a general spatial mesh and spatially varying scattering and absorption cross sections under some mild assumptions.
\end{remark}
%%%%%%%%%%%
\begin{remark}
\label{rem:condH:2} One can gain more insight on the conditioning of $\mH$ by combining with the stability analysis in Section \ref{sec:stability}. Specifically when the problem 
is in the transport dominated regime with $\vareps=O(1)$
and 
the time step condition $\dt=O(\vareps \dx)$, then $\textrm{cond}_2{\mH}=O(1)$. When the problem  
is in the scattering dominated regime with $\vareps\ll 1$ 
 we will have $\textrm{cond}_2{\mH}=O(1+\dt/h^2).$  In actual implementations, unconditional stability allows one to take $\dt=O(\dx)$ in this regime, rendering  $\textrm{cond}_2{\mH}=O(1/h)$. If one uses iterative methods e.g. Conjugate Gradient Method to solve the linear systems arising from Strategy 3, the iteration number needed to reach a fixed error tolerance will stay about the same in the transport dominated regime, and will increase by a factor of $\sqrt{2}$ when the number of mesh elements (i.e. $N_x$) doubles in the scattering dominated regime. {This is demonstrated numerically in \cite{Matsuda2025RPI} (see Section 2.7.1.5 and Tables 2.9-2.10). In practical simulations, pre-conditioning techniques (e.g., multi-grid pre-conditioners) will be needed for the ultimate efficiency if Conjugate Gradient or other Krylov-space based iterative solvers are used.}

\end{remark}

\section{Initialization: two approaches}
\label{sec:init}

At $t=t^0=0$, the numerical solution $\macro_h^0$, $\{\micro_{h,q}^0\}_{q=1}^{N_\velangle}$ are given e.g. via  the $L^2$ projection.  When the IMEX-BDF$\IMEXBDForder$-DG$\IMEXBDForder$-$\mS$$k$ scheme with $\IMEXBDForder=2,3$ is used, with the time integrators being multi-step, an additional initialization strategy is needed to provide  accurate numerical solutions at $t^1=t^0+\dt_1$, and also at $t^2=t^1+\dt_2$ when $\IMEXBDForder=3$.  Two approaches are considered in this work, and they will not affect the stability of the overall algorithm.  We use $\dt$ to denote the time step size predicted by the stability analysis for IMEX-BDF$\IMEXBDForder$-DG$\IMEXBDForder$-$\mS$$k$, e.g. in  \eqref{eq:fourier:dt:st12}-\eqref{eq:fourier:dt:st3} when $\absorp=0$.

\medskip
\noindent{\bf The First Approach.} The first one is to apply the one-step IMEX-RK$(\IMEXBDForder-1)$-DG$\IMEXBDForder$-$\mS$$k$ method with $\dt_n=\dt$ to compute the numerical solution at $t^n$, $n=1,\cdots,s-1, s=2,3$. Here
the IMEX-RK$(\IMEXBDForder-1)$-DG$\IMEXBDForder$-$\mS$$k$ scheme uses the same spatial and angular discretizations along with the same IMEX  partitionings  as the IMEX-BDF$\IMEXBDForder$-DG$\IMEXBDForder$-$\mS$$k$ scheme, yet adopt in time the IMEX-RK scheme of order $\IMEXBDForder-1$ used in \cite{jang2015high,peng2021asymptotic}.  Though such IMEX-RK scheme is an $(s-1)$th order time integrator,  its local errors at $t^n$, $n=1,\cdots,s-1$, will still be $s$-th order. As a variation, one can alternatively use the IMEX-RK$\IMEXBDForder$-DG$\IMEXBDForder$-$\mS$$k$ scheme for initialization.

\medskip
\noindent{\bf The Second Approach.} The second approach is to utilize that IMEX-BDF schemes come as a family, and one will apply the schemes  in the same family that have gradually growing yet lower order temporal accuracy (as mentioned in many classical texts, yet with little detail provided). We will present the approach using the IMEX-BDF3-DG3-$\mS$$k$ schemes as an example, that involve a three-step third order BDF method in time.  This initialization strategy first applies the one-step IMEX-BDF$1$-DG3-$\mS$$k$ scheme to compute the solution at $t^1=t^0+\dt_1$, and then applies the two-step IMEX-BDF$2$-DG3-$\mS$$k$ scheme (with variable time steps) to compute the solution at $t^2=t^1+\dt_2$. The idea is natural yet it needs to be carried out with great care to ensure the designed accuracy.

\begin{figure}[h]
\centering
\includegraphics[width=0.9\textwidth, height=4cm]{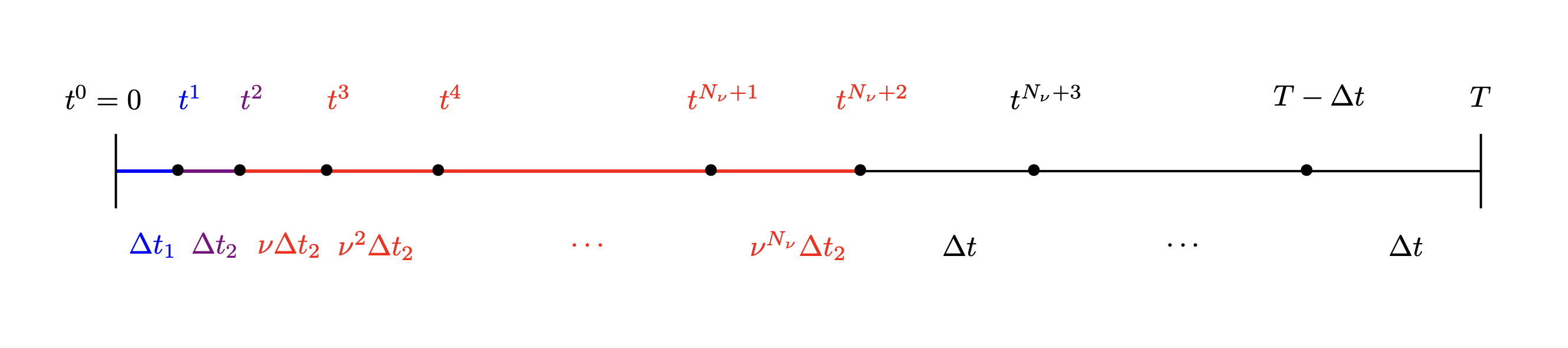}
\vspace{-0.2in}
\caption{Initialization for IMEX-BDF3-DG3-$\mS$$k$ by methods of the same family of lower order temporal accuracy. Adaptive time-stepping region: by IMEX-BDF1-DG3-$\mS$$k$ (blue), IMEX-BDF2-DG3-$\mS$$k$ (purple); transitional region (red) and  uniform time-stepping region (black) by IMEX-BDF3-DG3-$\mS$$k$ (black).}
\label{fig:init2}
\end{figure}

Based on our extensive numerical experiments {(see Tables 2.4-2.6 in \cite{Matsuda2025RPI})},  the following strategies and guidelines are identified and will be followed in our implementation. 
Firstly, to control the initial errors, the first two time steps $\dt_1$ and $\dt_2$ will be selected using  the adaptive time-stepping strategy in \cite{yan2021adaptive} with a prescribed error tolerance $err_{tol.init}$ for both $\macro$ and $\micro$. Secondly, to avoid the deterioration of the accuracy order  due to the drastic changes in time step sizes, namely $\dt_{n+1}/\dt_n$ being too large or too small, we require $\dt_{n+1}/\dt_n\in [1/\nu, \nu]$, with a hyper-parameter $\nu\in(1,2]$. Particularly,  $\dt_{2}/\dt_1\in [1/\nu, 1]$ is imposed, and  a transitional region, where $\dt_{n+1}/\dt_n=\nu$ is imposed,  is introduced between $t=t^{\IMEXBDForder-1}= t^2$ and the bulk part of the computational domain using the constant time step $\dt$ (that is determined by stability).  Thirdly, when accuracy order is examined through mesh refinements (e.g. $h\rightarrow h/2$), the initial error tolerance $err_{tol.init}$ will be reduced accordingly, namely $err_{tol.init}\rightarrow\frac{err_{tol.init}}{2^s}$ (with $\IMEXBDForder=3$). With this, much  larger $err_{tol.init}$ can be used on coarser meshes. Following the strategies described above, the computational domain can be seen to consist  of three regions: the adaptive time-stepping region, the transitional region, and the uniform time-stepping region, as illustrated in Figure \ref{fig:init2}. Without the transitional region, or using too large or too small $\nu$ in the transitional region, {or using a fixed yet not sufficiently small $err_{tol.init}$ during mesh refinement,} order reduction is observed numerically \cite{Matsuda2025RPI}.  Before entering the uniform time-stepping  region, IMEX-BDF schemes with variable time steps \cite{wang2008variable} will be applied in time.

\begin{remark}
It is known that the drastic change in step sizes in variable time-step multi-step time integrators can violate zero stability \cite{hairer1987solving, wang2008variable}. In the second initialization approach, the stability of the overall algorithm is governed by the scheme in the bulk part of the computational domain (i.e. the uniform time-stepping region).  Our choice of the hyper-parameter $\nu\in (1,2]$ is for the consideration of accuracy and efficiency.
\end{remark}
\begin{remark}
Let us also examine the size of the transitional region. The number of steps in this region, denoted as ${N_\nu}$,   is chosen as the smallest integer $N_\nu$ satisfying  $1/\nu\leq\frac{\dt}{\nu^{N_\nu}\dt_2}\leq \nu$,  namely, $N_\nu\in[\log_\nu(\dt/\dt_2)-1, \log_\nu(\dt/\dt_2)+1]$. From Figure \ref{fig:init2}, one can see that the size of the transitional region is $\dt_2\sum_{m=1}^{N_\nu} \nu^m$, that can be bounded as below,
$$\dt_2\sum_{m=1}^{N_\nu} \nu^m=\dt_2\frac{\nu(\nu^{N_\nu}-1)}{\nu-1}=\frac{\dt_2}{\dt}\frac{\nu(\nu^{N_\nu}-1)}{\nu-1} {\dt}\leq \frac{\nu(\nu^{N_\nu}-1)}{\nu^{N_\nu-1}(\nu-1)}\dt \leq  \frac{\nu^2}{\nu-1}\dt.$$
In our actual simulation,  $\nu=2$ is used, and the respective transitional region is no more than $4\dt$. 
\end{remark}

\begin{remark}  
Given the error estimators used in the adaptive time-stepping strategy \cite{yan2021adaptive} are based on Taylor's formula, the second initialization approach does not suit well the situation when initial layers are present in the solutions with {a
non-well-prepared initial condition.}
 In this case, we will use the first initialization approach, with the first one or two time step sizes modified as in \cite{peng2020stability} (see Section 5.3). This will be discussed more in the next  section and numerically illustrated in Section \ref{sec:num:non-well-prepared}.
\label{rem:init:3}
\end{remark} 

%%%%%%%%%%%%%%%%%%%%%%%
\section{AP property: a formal analysis}
\label{sec:AP}
%%%%%%%%%%%%%%%%%%%%%%%

In this section, we assume $\vareps\ll 1$ and perform a {\it formal} asymptotic analysis for the proposed three families of methods, namely \eqref{eq:fully:st1}-\eqref{eq:fully:st3},  with  constant cross sections $\scat$ and $\absorp$ (see Remark \ref{rem:ap:varyingcrossS} for  more  general case) and  fixed mesh parameters $\dx$ and $\dt$. Recall a numerical method for the equation \eqref{eq:micro-macro} is AP if its limiting scheme as $\vareps\rightarrow 0$ is a consistent scheme for the limiting diffusive equation \eqref{eq:diff-limit}. Throughout this section, a notation $A=\mO(\vareps^r)$ means: $|A|\leq {\mathcal C} |\vareps^r|$, where the constant ${\mathcal C}$ is independent of $\vareps$, yet can depend on model parameters $\scat$ and $\absorp$, mesh parameters $\dx$ and $\dt$, {and possibly also the numerical solutions at previous time steps. }

\begin{lemma}  Assume  $\macro^{n+k}_{h}=\mO(1)$, $\forall k=0,\dots, \IMEXBDForder-1$. Assume $\micro_{h,q}^{n+k}=\mO(1)$, $\forall k=0,\dots,\IMEXBDForder-1$, $\forall q=1,\dots,N_{\velangle}$, or $\micro_{h,q}^{n+k}=\mO(1/\vareps)$ for some integer $k\in[0,\IMEXBDForder-1]$ and some integer $q\in[1,N_{\velangle}]$.
\begin{itemize}
    \item [a.)]  If
$\micro_{h,q}^{n+k}=\mO(1)$, $\forall k=0,\dots, \IMEXBDForder-1$, $\forall q=1,\dots, N_{\velangle}$, the numerical solution by the proposed methods with Strategy 1 in \eqref{eq:fully:st1} satisfies 
\begin{equation}
\micro_{h,q}^{n+\IMEXBDForder} = - \frac{\velangle_q}{ \IMEXBDFimpcoef_{\IMEXBDForder} \scat}   \sum_{k=0}^{\IMEXBDForder-1} \IMEXBDFexpcoef_k\Dminusop \macro_h^{n+k}+\mO(\vareps)
, \quad q = 1, \dots, N_{\velangle},\quad \macro_{h}^{n+\IMEXBDForder}=\mO(1),
\label{eq:ap:01}
\end{equation}
while that with either Strategy 2 in \eqref{eq:fully:st2}  or Strategy 3 in \eqref{eq:fully:st3}  satisfies 
 \begin{equation}
\macro_{h}^{n+\IMEXBDForder}=\mO(1),\quad \micro_{h,q}^{n+\IMEXBDForder} =- \frac{\velangle_q}{\scat}   \Dminusop \macro_h^{n+s}+\mO(\vareps) 
, \quad q = 1, \dots, N_{\velangle}.
\label{eq:ap:02}
\end{equation}

    \item [b.)] If $ \micro_{h,q}^{n+k}=\mO(1/\vareps)$ for some integer $k\in [0, \IMEXBDForder-1]$ and some integer $q\in [1, N_\velangle]$, then with Strategy 1 in \eqref{eq:fully:st1}, the numerical solution satisfies 
\begin{equation}
\micro_{h,q}^{n+\IMEXBDForder} =- \frac{\velangle_q}{ \IMEXBDFimpcoef_{\IMEXBDForder} \scat}   \sum_{k=0}^{\IMEXBDForder-1} \IMEXBDFexpcoef_k\Dminusop \macro_h^{n+k} 
+\mO(1) +\mO(\vareps),  \quad q = 1, \dots, N_{\velangle}, \quad \macro_{h}^{n+\IMEXBDForder}=\mO(1),
\label{eq:ap:03}
\end{equation}
while  with Strategy 3 in \eqref{eq:fully:st3}, the numerical solution satisfies 
\begin{equation}
\macro_{h}^{n+\IMEXBDForder}=\mO(1),\quad \micro_{h,q}^{n+\IMEXBDForder} =- \frac{\velangle_q}{\scat}   \Dminusop \macro_h^{n+\IMEXBDForder} 
+\mO(1) +\mO(\vareps), \quad q = 1, \dots, N_{\velangle}.
\label{eq:ap:04}
\end{equation}
In both cases, 
\begin{equation}
    \micro_{h,q}^{n+\IMEXBDForder} =\mO(1),\quad \forall q = 1, \dots, N_{\velangle}.
    \label{eq:ap:05}
\end{equation}
\end{itemize} 
\label{lem:AP}
\end{lemma}

{The proof of Lemma \ref{lem:AP} is provided in Appendix \ref{sec:proof:AP}.} With this lemma and mathematical induction, the next theorem will follow.
%%%%%% THM %%%%%%%%%%%%
\begin{theorem}
\label{thm:AP}
{Assume  $\macro^{k}_{h}=\mO(1)$, $\forall k=0,\dots, \IMEXBDForder-1$. Assume $\micro_{h,q}^{k}=\mO(1)$, $\forall k=0,\dots,\IMEXBDForder-1$, $\forall q=1,\dots,N_{\velangle}$, or $\micro_{h,q}^{k}=\mO(1/\vareps)$ for some integer $k\in[0,\IMEXBDForder-1]$ and some integer $q\in[1,N_{\velangle}]$.}
\begin{itemize}
    \item [a.)]  If
$ \micro_{h,q}^{k}=\mO(1)$, $\forall k=0,\dots, \IMEXBDForder-1$, $\forall q=1,\dots, N_{\velangle}$, then in the limit of $\vareps\rightarrow 0$,  the proposed methods with Strategy 1 in \eqref{eq:fully:st1} or Strategy 2  in \eqref{eq:fully:st2}
give the following, 
\begin{equation}
\macro_h^{n+\IMEXBDForder} = 
\sum_{k=0}^{\IMEXBDForder-1} \IMEXBDFprevcoef_k \macro_h^{n+k} + \dt  
\sum_{k=0}^{\IMEXBDForder-1} \IMEXBDFexpcoef_k \left(\lvg\velangle^2\rvg\Dplusop\left(\frac{\Dminusop \macro_h^{n+k}}{\scat}\right) -\absorp \macro_h^{n+k}
+G_h\right),
\label{eq:ap:01A}
\end{equation} 
for any $n\geq n_\star$, with $n_\star=0$ for Strategy 1 and 
{$n_\star=\IMEXBDForder$} for Strategy 2,  while the proposed methods with Strategy 3 in \eqref{eq:fully:st3} lead to 
\begin{equation}
\macro_h^{n+\IMEXBDForder} = 
\sum_{k=0}^{\IMEXBDForder-1} \IMEXBDFprevcoef_k \macro_h^{n+k} + \dt  \IMEXBDFimpcoef_\IMEXBDForder 
\left(\lvg\velangle^2\rvg{\Dplusop\left(\frac{\Dminusop \macro_h^{n+\IMEXBDForder}}{\scat}\right)} -\absorp \macro_h^{n+\IMEXBDForder}
+G_h\right), 
\label{eq:ap:02A}
\end{equation} 
for any $n\ge n_\star$, with $n_\star=0$.
Both \eqref{eq:ap:01A} and \eqref{eq:ap:02A} are consistent discretizations of the limiting diffusive equation \eqref{eq:diff-limit}, with the former using the explicit part of the $\IMEXBDForder$-step IMEX-BDF time integrator in \eqref{eq:imex:BDF}, and the latter using the implicit part of the $\IMEXBDForder$-step IMEX-BDF time integrator.

    \item [b.)] If $ \micro_{h,q}^{k}=\mO(1/\vareps)$ for some integer $k\in [0, \IMEXBDForder-1]$ and some integer $q\in [1, N_\velangle]$, then in the limit of $\vareps\rightarrow 0$,  the proposed methods with Strategy 1 in \eqref{eq:fully:st1} give \eqref{eq:ap:01A} for $n\geq n_\star$, while those  with Strategy 3 in \eqref{eq:fully:st1} give \eqref{eq:ap:02A} for $n\geq n_\star$. Here
    {$n_\star=\IMEXBDForder$}.
    \end{itemize}
\end{theorem}

Based on Theorem \ref{thm:AP}, the type of AP property, captured by the different values of $n_\star$ in different scenarios, depends on the IMEX  partitioning, and it also depends on the numerical initialization when $\IMEXBDForder>1$.  With Strategy 1 and Strategy 3, the macroscopic component  $\macro_h$ drives the discrete dynamics and ensures the solution correctly exits the initial layer if there is any, while with Strategy 2, it is the microscopic component  $\micro_h$ {that} drives the dynamics instead. In general, {the initial condition at $t^0=0$ produces $\macro^0_h=\mO(1)$, and $\micro^0_{h,q}=\mO(1)$ or $\mO(1/\vareps)$ for a given $q$.}
We refer to the initial condition with $\micro^0_{h,q}=\mO(1), \forall q=1,\dots, N_\velangle$ as being {\it well-prepared}.  
\begin{itemize}
    \item 
With a well-prepared initial condition, both the initialization approaches in Section \ref{sec:init} will lead to bounded solutions, namely, $\macro^k_h=\mO(1), \micro^k_{h,q}=\mO(1), \forall k=1,\dots, \IMEXBDForder-1, \forall q=1, \dots, N_\velangle$ during the initialization stage. This will lead to the limiting scheme as in  scenario a.) in the theorem. Furthermore, the value of $n_\star$ for Strategy 2 can be reduced if the initialization relaxes the solution to its local equilibrium  $\micro_{h,q}^{k} =- \frac{\velangle_q}{\scat}   \Dminusop \macro_h^{k}+\mO(\vareps)$ sooner for some $k\leq \IMEXBDForder-1$.

\item
{When the initial condition} is not well-prepared, only the first initialization approach  in  the previous section shall be used due to the presence of an initial layer (also see Remark \ref{rem:init:3}). Moreover, with Strategy 1 or 3, the IMEX-RK-DG methods {used in the initialization}, like the IMEX-BDF-DG methods in \eqref{eq:fully:st1} and \eqref{eq:fully:st3}, are driven by the macroscopic component $\macro_h$ and this ensures  the numerical solutions correctly exit the initial layer with a bounded $\micro_h$ after one time step.

\item In Lemma \ref{lem:AP} or Theorem \ref{thm:AP}, we exclude the case with Strategy 2 when {the initial condition}  is not well-prepared, as the respective methods will lead to a nonphysical solution with $\macro_{h}=\mO(1/\vareps)$. This situation can be avoided via suitably chosen initialization, e.g. by using the first  initialization approach combined with IMEX Strategy 1 or 3. The resulting IMEX-BDF-DG-$\mS$2 {scheme} will then be AP after some initial {time} steps.
\end{itemize}

\begin{remark}
\label{rem:non-well-prep-reduction} {When the initial condition} is not well-prepared, the $\mO(1)$ term in \eqref{eq:ap:03}-\eqref{eq:ap:04} will render a first order temporal error  during the initialization regardless of the value of $\IMEXBDForder$ in the methods. This can be showed similarly as in Section 5 of \cite{peng2020stability}. The remedy for the associated order reduction, just as suggested by Remark 5.2 in \cite{peng2020stability}, is to apply smaller time  step sizes during the initialization stage. {This will also  be illustrated numerically in Section \ref{sec:num:non-well-prepared}.}
\end{remark}
\begin{remark}
\label{rem:ap:varyingcrossS}
  The AP property can be analyzed for the proposed methods with  general scattering and absorption cross sections $\scat=\scat(x)$ and $\absorp=\absorp(x)$, especially with the nodal treatment in \eqref{eq:varying-crossS2}. One can refer to Section 4 of \cite{peng2020asymptotic} for such  analysis.
\end{remark}

%----------------------------------------------------------------------------------------------------------------

\section{Numerical Examples}
\label{sec:num}

In this section, we demonstrate the performance of the IMEX-BDF$s$-DG$s$-$\mS$$k$ methods, $s, k=1,2,3$, when they are applied to examples with smooth or low-regularity solutions.  Unless stated otherwise, our methods of second and third  order use the first approach in Section \ref{sec:init} for the initialization. When numerical solutions are visualized (e.g. in Figure \ref{fig:two-material epsilon 1 T 15e-1 rho}),  their element-wise averages are plotted;  when point-wise errors are visualized (e.g. in Figure \ref{fig:two-material epsilon 1 T 15e-1 rho error}), several points are sampled and plotted per element. All methods are implemented by our own codes, written in MATLAB, and the linear system \eqref{eq:alg-6} is solved by {\it mldivide}, also known as {\it backslash}.

\subsection{Example 1: {Smooth example with constant material properties}} 
\label{subsec:smooth example}

In this example, we consider the one-group transport equation in slab geometry on $\spaceset = [0,2\pi]$, with {constant material properties} $\scat= 1$ and $\absorp= 0$ and zero source. The solution is smooth, with the well-prepared initial data $\macro(x,0) = \sin(x)$ and $\micro(x,\velangle,0) = -\velangle \cos(x)$, and periodic boundary conditions. The final time is $T=1$, unless otherwise specified.  In velocity, 16-point normalized Gauss-Legendre quadrature  is used with $N_v=16$. Uniform meshes are used in space with  $N_x=N$ elements and $h=2\pi/N$.
The $L^1$ errors and convergence orders for $\macro$ and $\lvg \velangle \micro \rvg$ are computed as follows. 
\begin{subequations} \label{eq:L1 errors and convergence orders}
\begin{align}
E_N^{\macro} & = \| \macro_{\dx}(x,T) - \macro_{\frac{\dx}{2}}(x,T) \|_{L^1(\spaceset)}, & O_N^{\macro} & = \log_2(E_N^{\macro}/E_{2N}^{\macro}), \\
E_N^{\lvg \velangle \micro \rvg} & = \| \lvg \velangle \micro_{\dx}(x,\velangle,T) \rvg_h - \lvg \velangle \micro_{\frac{\dx}{2}}(x,\velangle,T) \rvg_h \|_{L^1(\spaceset)}, & O_N^{\lvg \velangle \micro \rvg} & = \log_2(E_N^{\lvg \velangle \micro \rvg}/E_{2N}^{\lvg \velangle \micro \rvg}).
\end{align}
\end{subequations}

\subsubsection{Accuracy and AP property}

\begin{figure}[h]
\centering{

\includegraphics[width=0.32\textwidth, height=4cm]  
{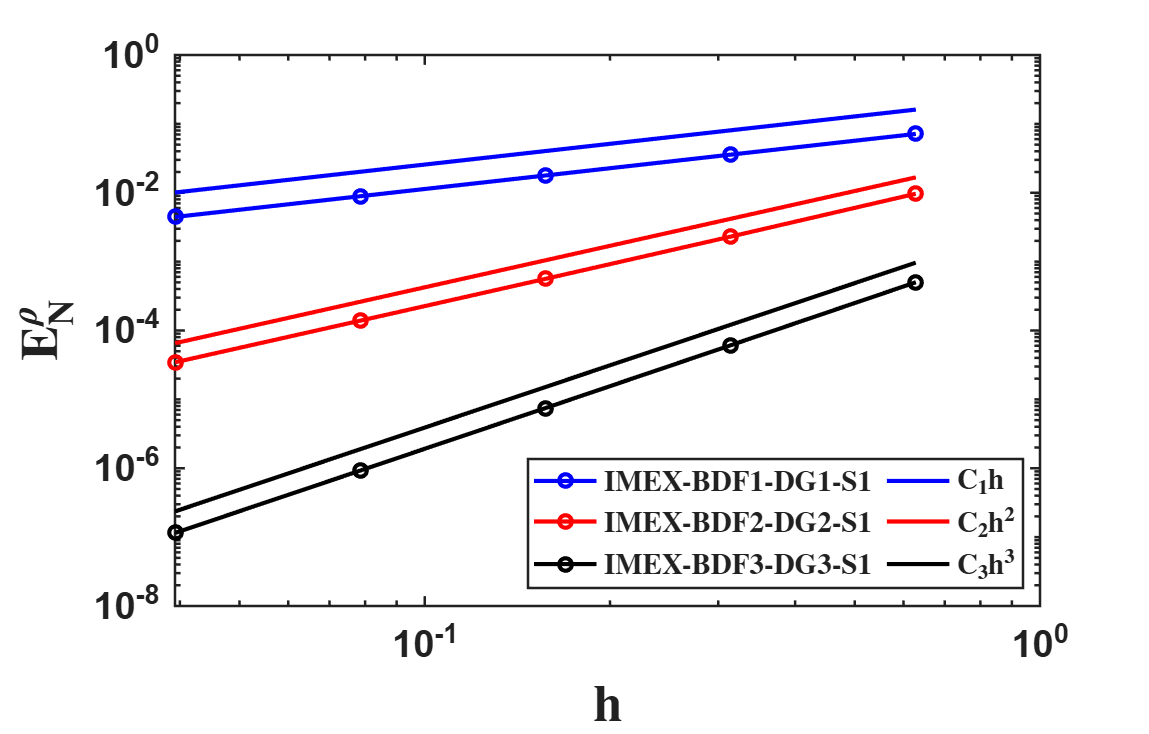}
\includegraphics[width=0.32\textwidth, height=4cm]   
{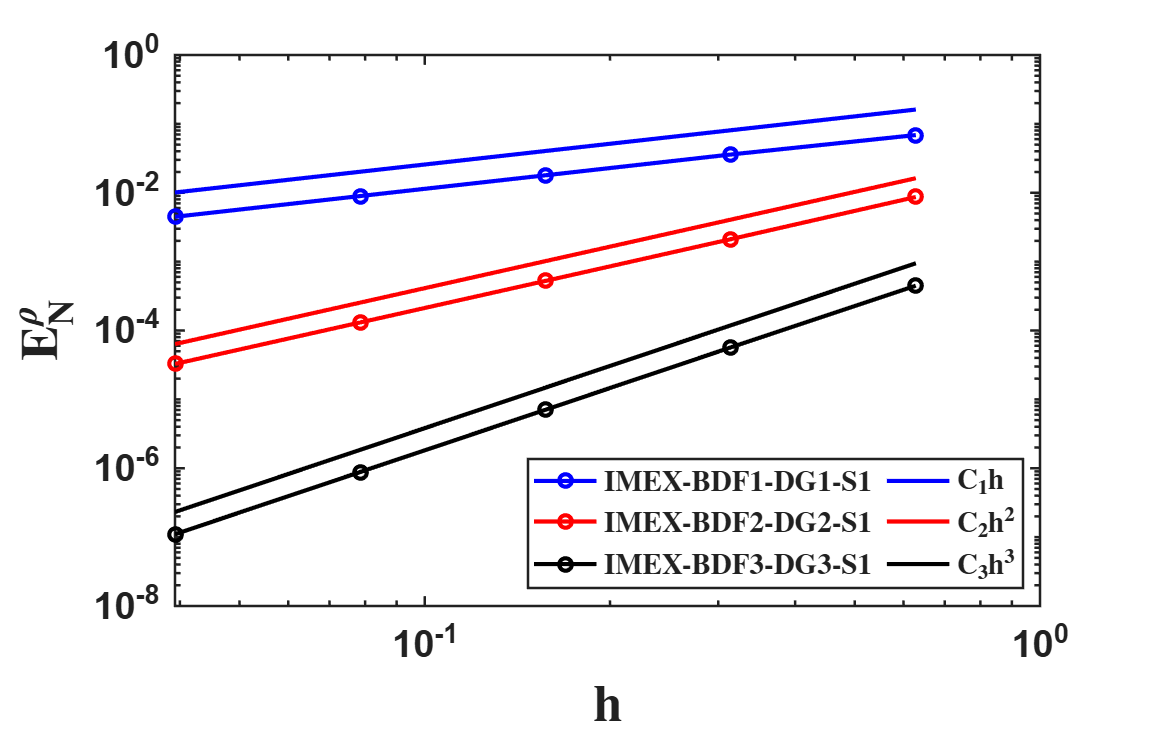}
\includegraphics[width=0.32\textwidth, height=4cm]  
{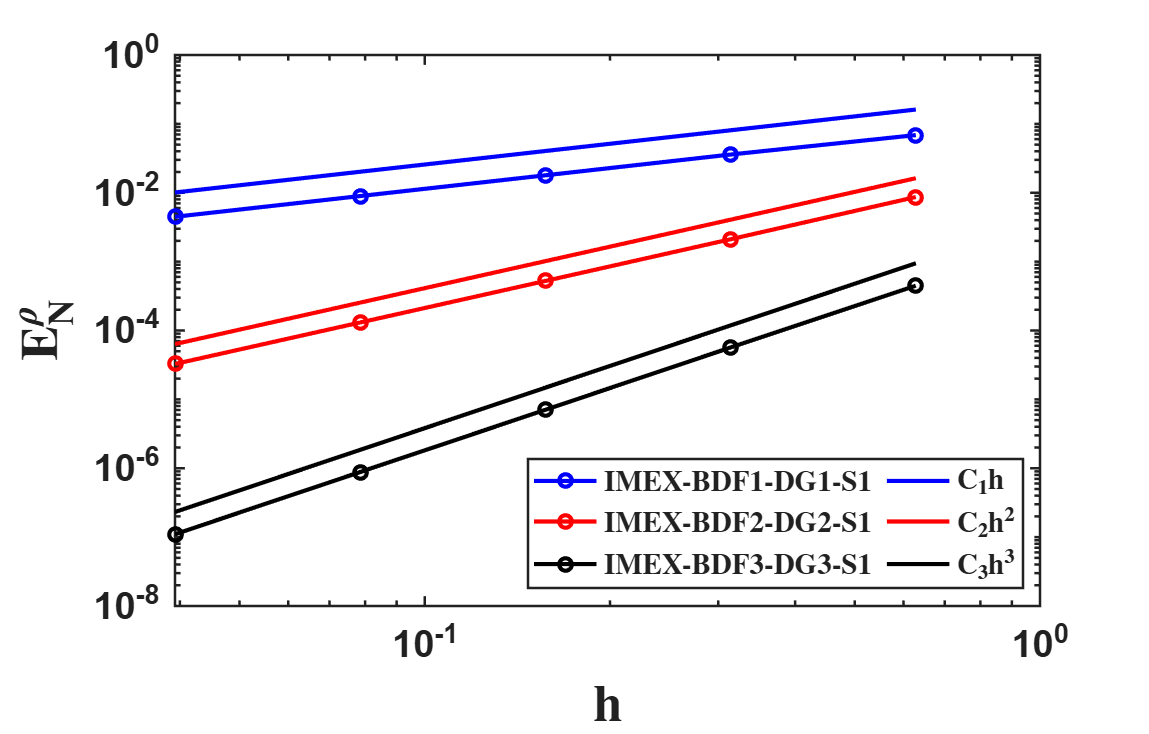}
\includegraphics[width=0.32\textwidth, height=4cm] 
{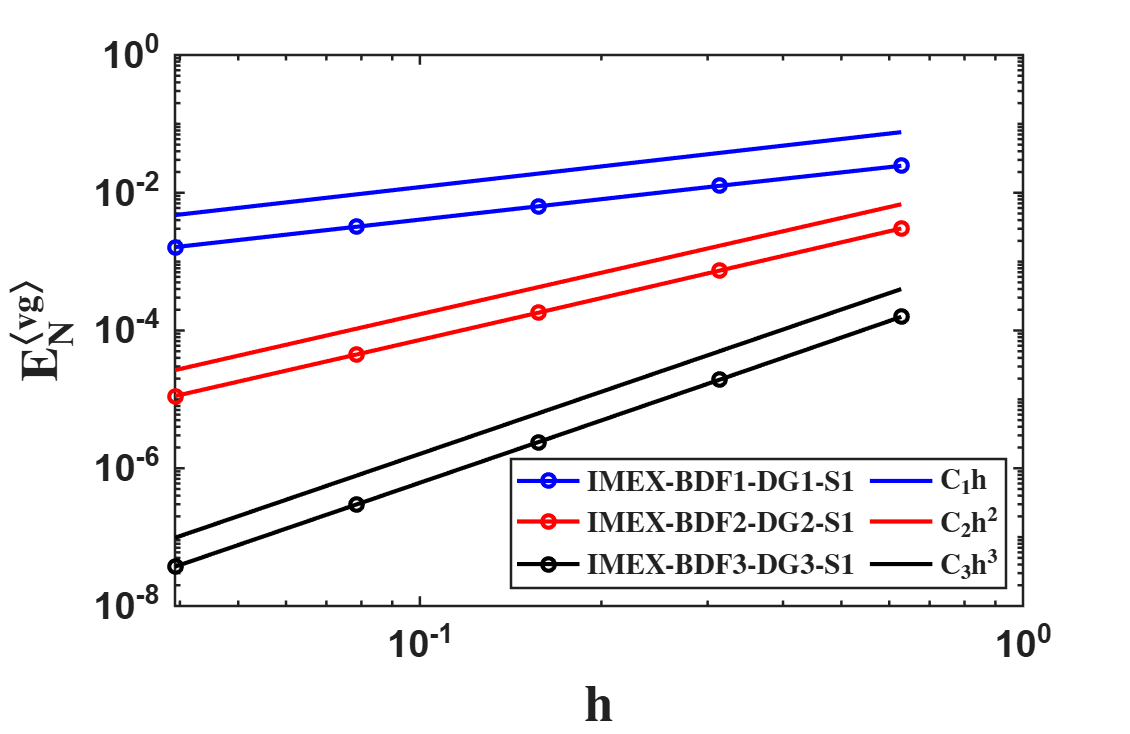}
\includegraphics[width=0.32\textwidth, height=4cm] 
{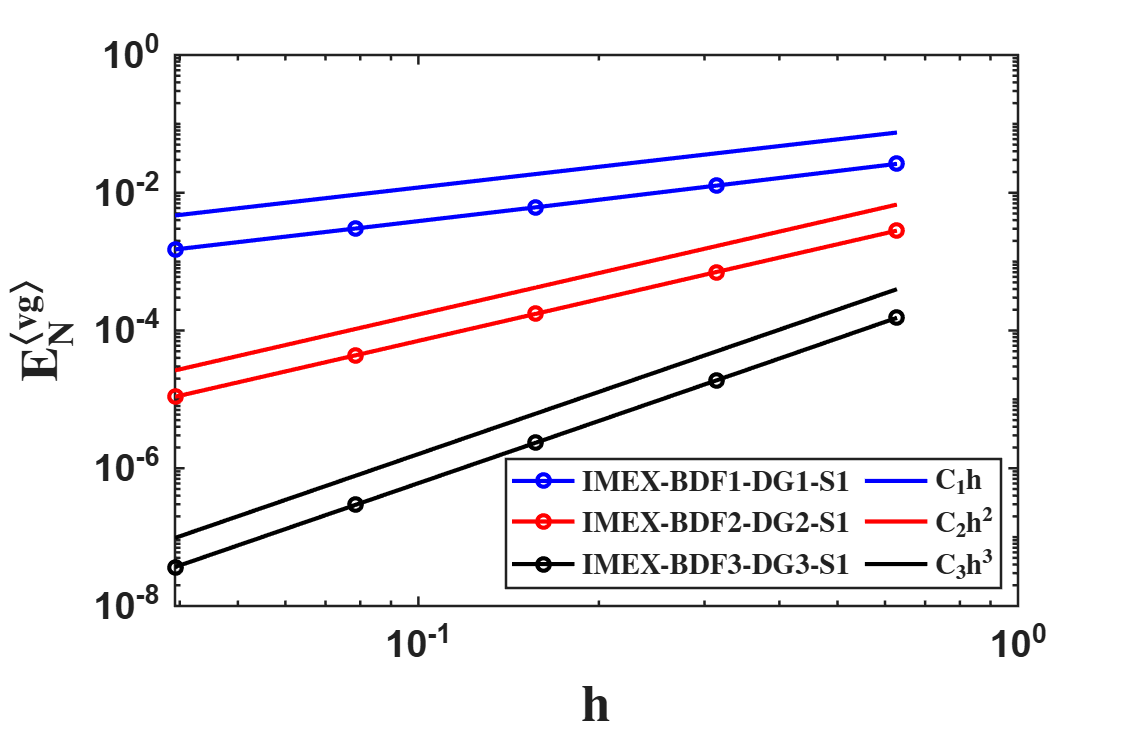}
\includegraphics[width=0.32\textwidth,height=4cm]  
{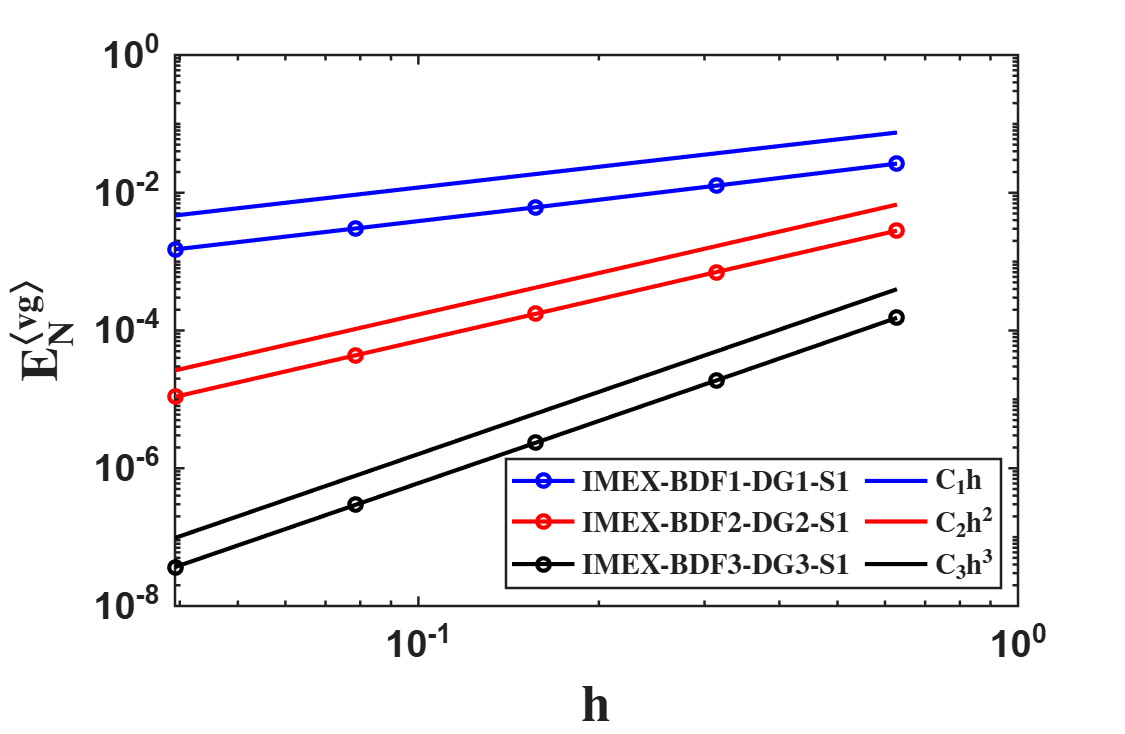}
}
\caption{Example 1: Smooth example with constant material properties. Convergence study for  IMEX-BDF$s$-DG$s$-$\mS$1, $s=1,2,3$. $L^1$ errors in  $\macro$ (top) and in $\lvg \velangle \micro \rvg$ (bottom); From left to right: $\vareps=0.5, 10^{-2}, 10^{-6}$.}
\label{fig:Smooth Example Periodic BC L1 Error Convergence}
\end{figure}

%%%%%%%%%%%%%%%%%%%%%%%%%%%%%%%%%
We first demonstrate the order of  convergence and the
AP property of the proposed IMEX-BDF$s$-DG$s$-$\mS$$k$ schemes, with $s, k=1,2, 3$, and $\vareps = 0.5, 10^{-2}, 10^{-6}$. In Figure \ref{fig:Smooth Example Periodic BC L1 Error Convergence}, we plot $h$ versus $L^1$ errors in $\macro$ (top row) and $\lvg \velangle \micro \rvg$ (bottom row) for {the methods with Strategy 1.} The methods  show their designed accuracy order of $s$. The uniform convergence orders for different values of $\vareps$ (especially when it is small) further support the AP property of the methods.  {For this smooth example,  results by the IMEX-BDF$s$-DG$s$-$\mS$$k$ schemes, with $k=2,3$, are not plotted, as they are visually indistinguishable from those by IMEX-BDF$s$-DG$s$-$\mS$$1$ schemes. To get some idea, the errors by IMEX-BDF$s$-DG$s$-$\mS$$k$ schemes, with $k=2,3$, are the same as those by IMEX-BDF$s$-DG$s$-$\mS$1 in their leading two or more significant digits (after the decimal point) on all $\vareps$ and meshes considered.}

\subsubsection{Comparison: Two approaches for initialization}
%%%%%%%%%%%%%%%%%%%%%%%%%%%%%%%%%%%%%%%%%%%
In Section \ref{sec:init}, two approaches are proposed for initialization when the proposed methods are  second and third order accurate. For this smooth example with a well-prepared initial condition, there is little difference in their performance, as illustrated by Table \ref{tab:Smooth Example Initialization Comparison IMEX-BDF2-DG2} for the second order IMEX-BDF2-DG2  method with Strategy 1. Here Approach 2 starts with $err_{tol.init}=10^{-4}$ as the adaptive time-stepping tolerance when $N=10$. One should note that the IMEX-BDF1 time integrator and the IMEX-RK1 time integrator being used are the same. Results on more refined meshes or by the third order methods are not shown as the two initialization approaches produce the same errors, at least up to the first four significant digits.

\begin{table}[ht] 
\centering 
\caption{Example 1:  Smooth example with constant material properties. To study the two initialization approaches  in Section \ref{sec:init} by comparing $L^1$ errors in  $\macro$ and in  $\lvg \velangle \micro \rvg$ at $T=1$ by IMEX-BDF2-DG2-$\mS$1.  Approach 1: IMEX-RK1-DG2-$\mS$1 with $\dt$ determined by \eqref{eq:fourier:dt:st12}; Approach 2: IMEX-BDF1-DG2-$\mS$1 with adaptive time-stepping + transitional region before the uniform time-stepping region with time step $\dt$.} 
\begin{tabular}{c|c|c|c|c|c}\hline 
& & \multicolumn{2}{|c|}{Approach 1} & \multicolumn{2}{|c}{Approach 2} \\ 
\hline 
$\varepsilon$ & $N$ & $E_N^{\rho}$ & $E_N^{\langle v g \rangle}$ & $E_N^{\rho}$ & $E_N^{\langle v g \rangle}$ \\ \hline 
 \multirow{4}{*}{$0.5$}  & 20 & 2.307E-03 & 7.417E-04 & 2.298E-03 & 7.383E-04 \\ 
& 40 & 5.620E-04 & 1.811E-04 & 5.599E-04 & 1.803E-04 \\ 
& 80 & 1.389E-04 & 4.483E-05 & 1.384E-04 & 4.463E-05 \\ 
& 160 & 3.451E-05 & 1.116E-05 & 3.438E-05 & 1.111E-05 \\ 
\hline 
 \multirow{2}{*}{$10^{-2}$} & 10 & 8.683E-03 & 2.820E-03 & 8.682E-03 & 2.820E-03 \\ 
& 20 & 2.116E-03 & 7.065E-04 & 2.116E-03 & 7.065E-04 \\ 
\hline 
\multirow{2}{*}{$10^{-6}$}& 10 & 8.670E-03 & 2.818E-03 & 8.669E-03 & 2.818E-03 \\ 
& 20 & 2.114E-03 & 7.062E-04 & 2.114E-03 & 7.062E-04 \\
\hline 
\end{tabular}
\label{tab:Smooth Example Initialization Comparison IMEX-BDF2-DG2}
\end{table}

\subsubsection{Comparison in cost efficiency and accuracy: BDF vs RK in time}
%%%%%%%%%%%%%%%%%%%%%%%%%%%%%%%%%%%%%%%%%%%
The focus of the present work is on using IMEX-BDF time integrators to achieve AP methods with high order temporal accuracy for the model \eqref{eq:kinetic transport equation}. As a family of multi-step methods, they have the potential to be more efficient compared with the IMEX-RK type time integrators, that are multi-stage, when the methods are higher than first order accurate, hence involve more function evaluations. On the other hand, {methods with a lower cost} may not produce more accurate numerical solutions. The {next} set of tests  is designed and performed to shed some light regarding the cost efficiency versus accuracy of the IMEX-BDF-DG methods and IMEX-RK-DG  methods when the algorithmic difference  solely comes from the use of time integrators, BDF vs RK. Particularly, we consider the second and third order IMEX-BDF$\IMEXBDForder$-DG$\IMEXBDForder$-$\mS3$  methods proposed here ($\IMEXBDForder=2,3$), and the methods in \cite{peng2021asymptotic} which differ only in using IMEX-RK methods in time.  We pick  IMEX Strategy 3 for the study, since the methods are  unconditionally stable in the diffusive regime with $\vareps\ll1$, and the methods with RK in time do not {seem to} suffer from order reduction (see next study). 
The resolution of the computed solutions will be measured by the $L^1$ errors in $\macro$ and $\lvg \velangle \micro \rvg$ defined in \eqref{eq:L1 errors and convergence orders}. With the similarity in the observations, only results in  $\macro$ are presented. For the computational time, we measure  the  average wall-clock time  for each scheme marching from $t=0$ to $T = 5$ in three repeated runs. To reach a more fair comparison, the initialization is done using a variation of the first approach, namely by using  the IMEX-RK$\IMEXBDForder$-DG$\IMEXBDForder$-$\mS3$ scheme to compute the  solutions at $t^j, j=1, \dots, s-1$ and  with  the same time step size as that for the IMEX-BDF$\IMEXBDForder$-DG$\IMEXBDForder$-$\mS3$ scheme.

\underline{In the first test}, we consider $\vareps = 0.5$. In this transport regime with  relatively weak scattering, both IMEX-BDF-DG and IMEX-RK-DG methods are conditionally stable. With BDF in time, the time step is set according to \eqref{eq:fourier:dt:st3}; 
with RK in time, 
the time step conditions of similar kind
can be derived by Fourier-based stability analysis, and the actual time step sizes are provided by equation (6.1) in  \cite{peng2021asymptotic}. That is, each scheme uses the ``largest'' time step size predicted by the Fourier stability analysis. 
\underline{In the second test}, we consider $\vareps = 10^{-6}$. In this diffusive regime  with strong scattering, both IMEX-BDF-DG and IMEX-RK-DG methods are unconditionally stable. We take the same time step size $\dt=0.5h$ in all simulations.
\underline{In the third test}, we still consider $\vareps = 10^{-6}$, except that we now  take into account RK time integrators are multi-stage. Particularly,   for the second order and the third order IMEX-RK methods  in \cite{peng2021asymptotic}, they are 2-stage and 4-stage, respectively. With these in  mind, while $\dt=0.5h$  is still used for IMEX-BDF-DG methods, we use larger time step sizes for IMEX-RK-DG methods, specifically, $\dt=2\times 0.5h =h$ for the second order IMEX-RK2-DG2 method, and $\dt=4\times 0.5h =2h$ for the third order IMEX-RK3-DG3 method. In other words, the effective time step sizes (i.e. the ratio of $\dt$ and the number of inner stages) are the same as   $0.5h$. 

\begin{figure}[ht]
\centering
{
\includegraphics[width=0.7\linewidth]
{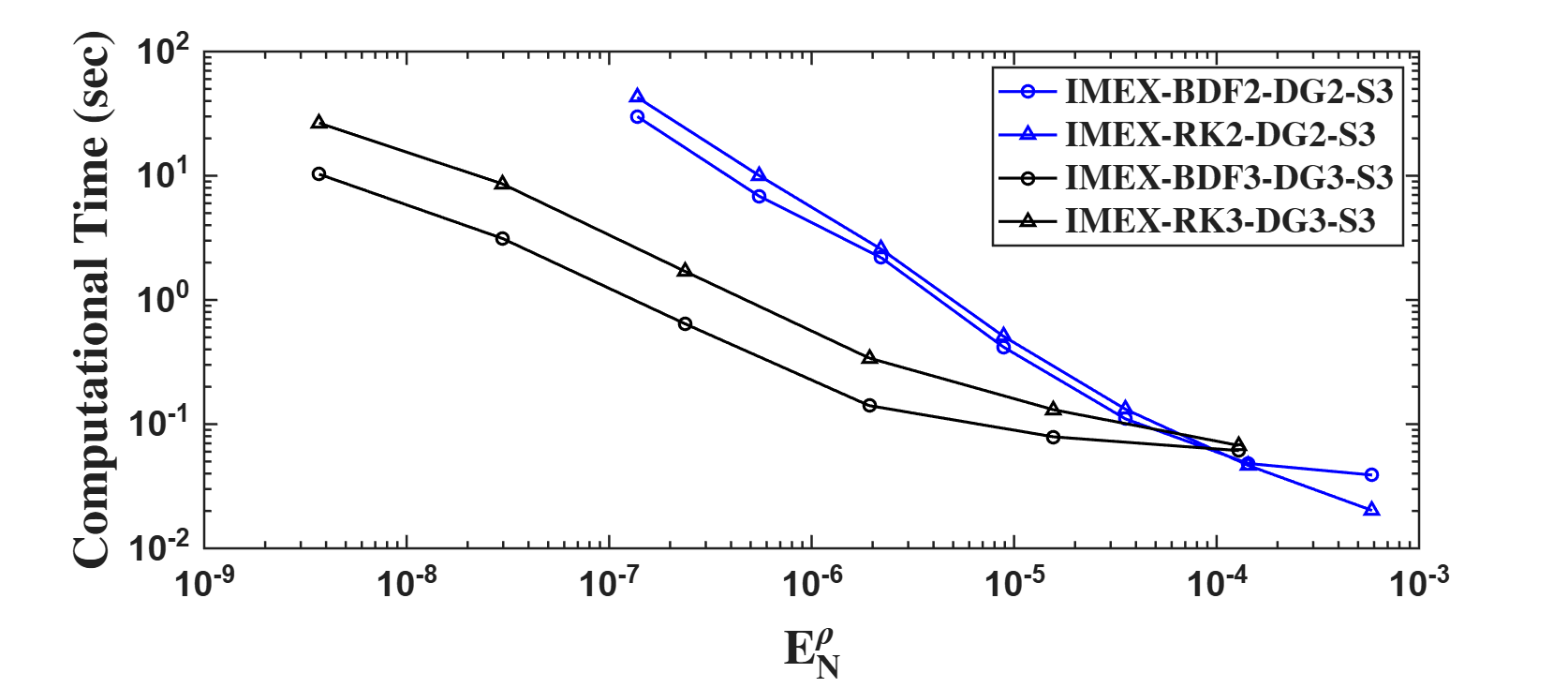}
\includegraphics[width=0.7\linewidth]
{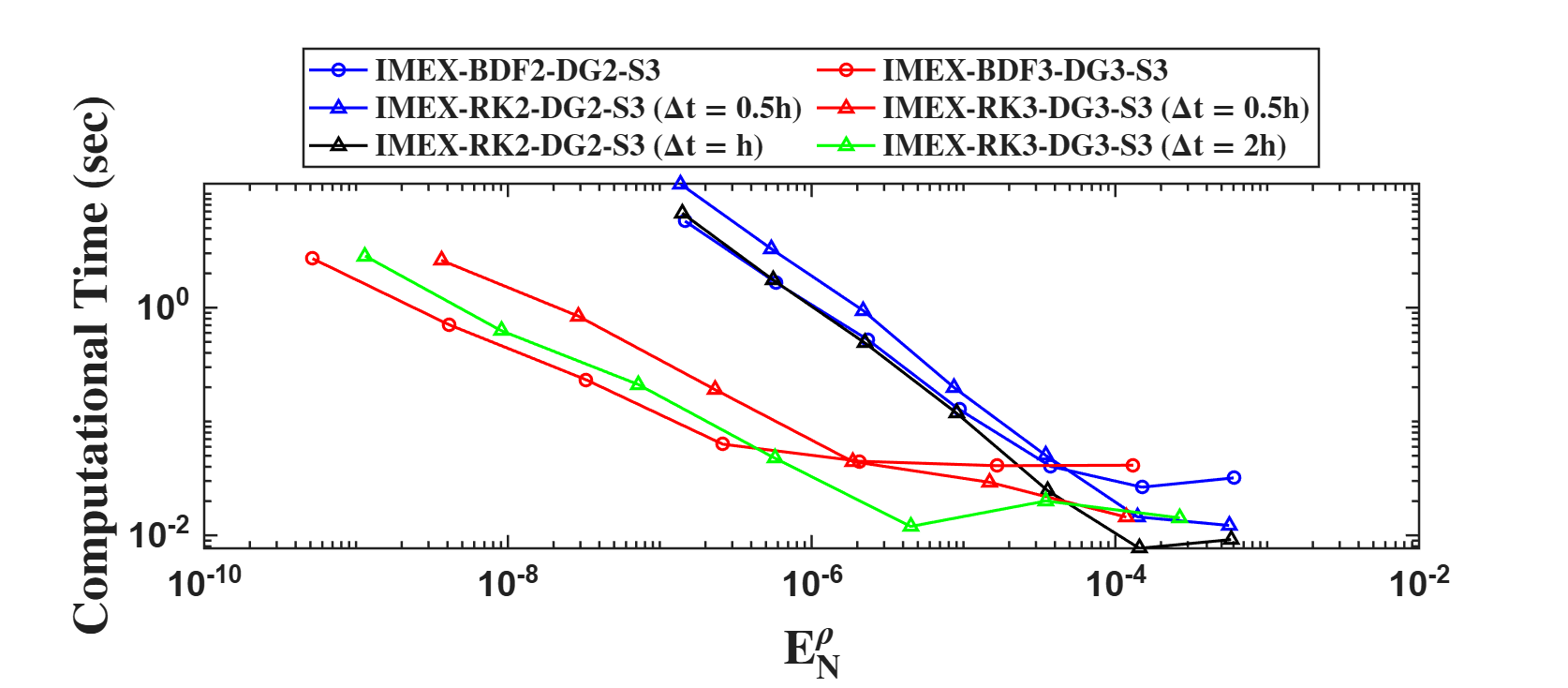}
}
\caption{Example 1:  Smooth example with constant material properties. To compare  IMEX-BDF$\IMEXBDForder$-DG$\IMEXBDForder$-$\mS$3  and IMEX-RK$\IMEXBDForder$-DG$\IMEXBDForder$-$\mS$3 by examining  $L^1$ errors in $\macro$ versus the computational times. Top:  $\vareps = 0.5$; Below:  $\vareps = 10^{-6}$.}
\label{fig:BDF-RK}
\end{figure}

In Figure \ref{fig:BDF-RK}, the $L^1$ errors in $\macro$ versus the computational times are reported, in subfigures separately for test 1 with $\vareps=0.5$ and for  test 2 \& 3 with $\vareps=10^{-6}$. They are based on the results computed from  a sequence of meshes, with a starting mesh of $N=20$ for the second order methods and of $N=10$ for the third order methods, and  a refinement with the factor 2.

For the methods of the same order of accuracy, it is observed that  in order to achieve the same but relatively higher resolution with smaller errors, methods with BDF {take less computational time than those with RK in time}, even with the number of inner stages being taken into account in the RK setting. {This difference is particularly meaningful for third order methods.} 
Methods with RK in time take less computational time to achieve the same but relatively lower resolution. Furthermore, to achieve the same level of errors especially of higher resolution, third order methods are much more cost effective than the second order counterparts,  as widely known now {for sufficiently smooth solutions. }

\subsubsection{A detour: order reduction  of IMEX-RK-DG-$\mS$1}
%%%%%%%%%%%%%%%%%%%%%
One reason that we do not compare the cost efficiency between the methods using BDF and RK in time with {\it all} IMEX partitionings is that the IMEX-RK-DG methods, particularly with the IMEX Strategy 1, can experience order reduction in $g$ 
when $\vareps$ is of moderate to small size, hence lose the uniform in $\vareps$ order of accuracy. The root of this phenomenon is the $O(\dt)$ error in approximating  $\micro+v\partial_x\macro=O(\vareps)$ by the {\it multi-stage} RK methods using IMEX  Strategy 1. This can be explained by  a similar formal analysis as in Appendix A of  \cite{peng2020stability} (also see its Remark 3.4), and illustrated numerically by the results in Table \ref{tab:order reduction} (panels on the left). Even though the $O(\dt)$-error is expected for $\lvg v\micro \rvg$ computed by the third order IMEX-RK3-DG3-$\mS$1 method, with the use of the time step size $\dt=O(\dx^2)$ in the considered regimes\footnote{The time step size for the IMEX-RK3-DG3-$\mS$$k$ {scheme}, with $k=1,2$, is determined by $\dt = 0.13 \vareps \dx + 0.015 {\scatlower} \dx^2$, a formula found through a similar Fourier analysis as in Section \ref{sec:fourier}.}, we have $O(\dt)=O(\dx^2)$ and observe second order accuracy in $\lvg v\micro \rvg$, confirming the order reduction.

Order reduction {seems to be avoided} if one instead  uses IMEX-BDF methods in time  with any IMEX partitioning (see Figure \ref{fig:Smooth Example Periodic BC L1 Error Convergence}), highlighting one {potential} advantage of working with BDF-type methods in time for solving multi-scale problems. Alternatively, if one uses the IMEX-RK methods in time along with  IMEX  Strategy 2 (see panels on the right in Table \ref{tab:order reduction}) or Strategy 3 (see \cite{peng2021asymptotic}), {order reduction is not observed.}

\begin{table}[htbp]
  \centering
  \caption{Example 1: Smooth example with constant material properties.  $L^1$ errors and orders at $T=1$ by IMEX-RK3-DG3-$\mathcal{S}$1 and IMEX-RK3-DG3-$\mathcal{S}$2.}
   % \medskip
    \begin{tabular}{c|c|c|c|c|c|c|c|c|c}
    \hline        
    \multirow{2}{*} {$\vareps$}& \multirow{2}{*} {$N$}&
       \multicolumn{4}{c|}{IMEX-RK3-DG3-$\mathcal{S}$1}& \multicolumn{4}{c}{IMEX-RK3-DG3-$\mathcal{S}$2}\\
       \cline{3-10}
	& & $E_N^{\macro}$ & order  &$E_N^{\lvg \velangle \micro \rvg}$&order&$E_N^{\macro}$ & order  &$E_N^{\lvg \velangle \micro \rvg}$ &order\\
    \hline
     \multirow{4}{*}{$10^{-3}$}  
    & 20 & 5.668E-05 & - & 3.062E-05 & - & 5.668E-05 & - & 1.889E-05 & -  \\
&40 & 7.061E-06 & 3.00 & 6.290E-06 & {\bf 2.28} & 7.061E-06 & 3.00 & 2.354E-06 & 3.00 \\ 
&80 & 8.820E-07 & 3.00 & 1.482E-06 & {\bf 2.09} & 8.820E-07 & 3.00 & 2.940E-07 & 3.00 \\ 
&160 & 1.102E-07 & 3.00 & 3.435E-07 & {\bf 2.11} & 8.820E-07 & 3.00 & 2.940E-07 & 3.00 \\ \hline
     \multirow{4}{*}{$10^{-6}$} 
&20 & 5.668E-05 & - & 3.031E-05 & - & 5.668E-05 & - & 1.889E-05 & - \\ 
&40 & 7.061E-06 & 3.00 & 6.117E-06 & {\bf 2.31} & 7.061E-06 & 3.00 & 2.354E-06 & 3.00 \\ 
&80 & 8.819E-07 & 3.00 & 1.412E-06 & {\bf 2.12} & 8.819E-07 & 3.00 & 2.940E-07 & 3.00 \\ 
&160 & 1.102E-07 & 3.00 & 3.448E-07 & {\bf 2.03} & 1.102E-07 & 3.00 & 3.674E-08 & 3.00 \\ \hline	      
          \end{tabular}%
           \label{tab:order reduction}
\end{table}%

%----------------------------------------------------------------------------------------------------------------
\subsection{Example 2: {Smooth example with spatially varying scattering cross section and nonzero constant source}}
\label{subsec:varying scattering}

{In this section, we consider another smooth example, set up similarly as Example 1 in Section \ref{subsec:smooth example}, except that the scattering cross section profile is spatially dependent,  specifically, 
\[\scat(x) = 1 + 10 \sin^2(x),\]  and $\source(x)=1$.  For this example with the spatially varying scattering cross section and nonzero source, just as in Section \ref{subsec:smooth example}, we examine and confirm the designed accuracy order $s$ for the proposed IMEX-BDF$s$-DG$s$-$\mS$$k$ schemes, with $s, k=1,2, 3$, and $\vareps = 0.5, 10^{-2}, 10^{-6}$.    To save space, we only present    in Table \ref{tab:Ex2:errOrderS3-rev-r2} the $L^1$ errors and orders of convergence for $\macro$ at the final time $T=1$ by the methods with IMEX Strategy 3. The time step is  determined by \eqref{eq:fourier:dt:st3} with $\scatlower =\min_{x\in\spaceset}\scat(x)=1$.  Once again, the uniform convergence orders for different values of $\vareps$ (especially when it is small) further support the AP property of the methods. }

\begin{table}[ht] 
\vspace{0.1in} 
\centering 
\caption{Example 2: Smooth example with spatially varying scattering cross section and nonzero constant source. $L^1$ errors and orders in $\macro$ at $T = 1$ by IMEX-BDF-DG-$\mathcal{S}$3.}
\begin{tabular}{c|c|c|c|c|c|c|c}\hline 
& & \multicolumn{2}{c|}{IMEX-BDF1-DG1-$\mathcal{S}$3} & \multicolumn{2}{c|}{IMEX-BDF2-DG2-$\mathcal{S}$3} & \multicolumn{2}{c}{IMEX-BDF3-DG3-$\mathcal{S}$3}\\ \hline 
$\varepsilon$ & $N$ & $E_N^{\rho}$ & order & $E_N^{\rho}$ & order & $E_N^{\rho}$ & order \\ \hline 
\multirow{4}{*}{$0.5$} & 40 & 2.496E-02 & - & 1.339E-03 & - & 5.167E-05 & - \\ 
& 80 & 1.252E-02 & 1.00 & 3.371E-04 & 1.99 & 6.584E-06 & 2.97 \\ 
& 160 & 6.273E-03 & 1.00 & 8.382E-05 & 2.01 & 8.179E-07 & 3.01 \\ 
& 320 & 3.140E-03 & 1.00 & 2.092E-05 & 2.00 & 1.017E-07 & 3.01 \\ 
\hline 
\multirow{4}{*}{$10^{-2}$} & 40 & 2.477E-02 & - & 1.265E-03 & - & 4.881E-05 & - \\ 
& 80 & 1.241E-02 & 1.00 & 3.189E-04 & 1.99 & 6.177E-06 & 2.98 \\ 
& 160 & 6.208E-03 & 1.00 & 7.960E-05 & 2.00 & 7.714E-07 & 3.00 \\ 
& 320 & 3.106E-03 & 1.00 & 1.990E-05 & 2.00 & 9.653E-08 & 3.00 \\ 
\hline 
\multirow{4}{*}{$10^{-6}$} & 40 & 2.477E-02 & - & 1.282E-03 & - & 5.338E-05 & - \\ 
& 80 & 1.241E-02 & 1.00 & 3.210E-04 & 2.00 & 6.575E-06 & 3.02 \\ 
& 160 & 6.208E-03 & 1.00 & 8.022E-05 & 2.00 & 8.241E-07 & 3.00 \\ 
& 320 & 3.105E-03 & 1.00 & 2.006E-05 & 2.00 & 1.027E-07 & 3.00 \\ 
\hline 
\end{tabular} 
\label{tab:Ex2:errOrderS3-rev-r2}
\end{table}

\subsection{Example 3: Two-material problem}
\label{sec:two-material}

In this example, we consider the one-group transport equation in slab geometry  that involves two materials on $\spaceset = [0,11]$, modeled by 
\begin{equation}
\scat = 0,\; \absorp = 1, \; \text{on} \; x \in [0,1]; \quad 
\scat = 100, \; \absorp = 0, \; \ \text{on} \; x \in [1,11], 
\end{equation}
along with zero source and $\vareps = 1$.  The inflow Dirichlet boundary conditions are 
\begin{equation*}
\angflux(0,\velangle,t)  = 5 \;\;\text{for} \; \;v\geq 0,\quad \angflux(11,\velangle,t)  = 0 \;\;\text{for}\;\; v\leq 0.
\end{equation*}
  The material over $[0,1]$ is purely absorbing, and the material over $[1,11]$ contributes to strong scattering. An isotropic inflow enters the domain from the left boundary and  becomes anisotropic, and this is followed by the formation of an interior layer between {the} two materials.  The initial condition is zero and the final time is $T=1.5$.

{To impose the  inflow boundary conditions,  the inflow-outflow close-loop numerical boundary treatment as in \cite{peng2021asymptotic} is applied, and the details can be found in its Section 6.1 with two specifications:  (1) $(\breve{\macro_h})_{1/2}=\macro_L(t)$ and $(\breve{\macro_h})_{N+1/2}=\macro_R(t)$ are treated explicitly to ensure $\Theta^{(3)}$ in \eqref{eq:alg-4} stays block-diagonal; (2)  $C_R=-1$  is used due to a typo in \cite{peng2021asymptotic}. In this treatment, terms such as $\int_{-1}^{0} \angflux(x,\velangle,t) d\velangle$ and $\int_{0}^{1} \angflux(x,\velangle,t) d\velangle$ need to be evaluated. Related, 2-panel 16-point normalized Gauss-Legendre quadrature  is used  for the velocity discretization with $N_v=32$.}

 With the presence of an interior layer in the solution,  a non-uniform spatial mesh is used, with  $\dx=1/20$ on $[0,1.5]$ and $\dx=19/40$ on $[1.5,11]$.  We also write  $\dx_{\text{min}}=1/20$ for this mesh.  The reference solution is computed using the IMEX-BDF3-DG3-$\mS$3 scheme with $\dx=1/320$ on $[0,1.5]$ and $\dx=19/640$ on $[1.5,11]$, with $\dx_{\text{min}}=1/320$. To determine  the time step size,  we first note that  $\min_{x\in\spaceset}\scat(x)=0$ and  $\max_{x\in\spaceset}\absorp(x)=1$.
Based on the insight we get from the stability analysis for the contribution of the absorption term to the time step conditions (see Theorem 
\ref{thm:stab:energy},
especially  Eqn. \eqref{eq:st1:dt} for IMEX Strategy 1 and 2, as well as Figure \ref{fig:Fourier:absorp:3rdOrder}), the following  time steps are used: 
\begin{itemize}
    \item 
With Strategy 1 and 2 (i.e. $k=1,2$), we set $\dt=\frac{\dt^{(k)}_{\text{CFL}\IMEXBDForder}}{1+0.5\dt^{(k)}_{\text{CFL}\IMEXBDForder}}$ with $\scatlower=0$ for $\IMEXBDForder=1,2,3$, along with $\dt^{(k)}_{\text{CFL}\IMEXBDForder}$ defined in \eqref{eq:fourier:dt:st12}. {That is, we set $\dt= \frac{0.7\vareps \dx_{\text{min}}}{1+0.35\vareps\dx_{\text{min}}}$ for the IMEX-BDF1-DG1 scheme, $\dt = \frac{0.151\vareps \dx_{\text{min}}}{1+0.0755\vareps\dx_{\text{min}}}$ for the IMEX-BDF2-DG2 scheme, and $\dt=\frac{0.062\vareps \dx_{\text{min}}}{1+0.031\vareps\dx_{\text{min}}}$ for the IMEX-BDF3-DG3 scheme.}
\item With Strategy 3, $\dt=\dt^{(3)}_{\text{CFL}\IMEXBDForder}$ is used with $\scatlower=0$ along with $\dt^{(3)}_{\text{CFL}\IMEXBDForder}$ defined in \eqref{eq:fourier:dt:st3}. 
{That is, we set $\dt=\min(\vareps, 0.5) \dx_{\text{min}}=0.5\dx_{\text{min}}$ for the IMEX-BDF1-DG1 scheme,  and $\dt= 0.2\vareps \dx_{\text{min}}, 0.07\vareps \dx_{\text{min}}$ for the IMEX-BDF$\IMEXBDForder$-DG$\IMEXBDForder$ scheme, with $\IMEXBDForder= 2, 3$ respectively.}
\end{itemize}

\begin{figure}[ht]

\begin{subfigure}{0.33\textwidth}
\includegraphics[width=\linewidth, height=4.5cm]{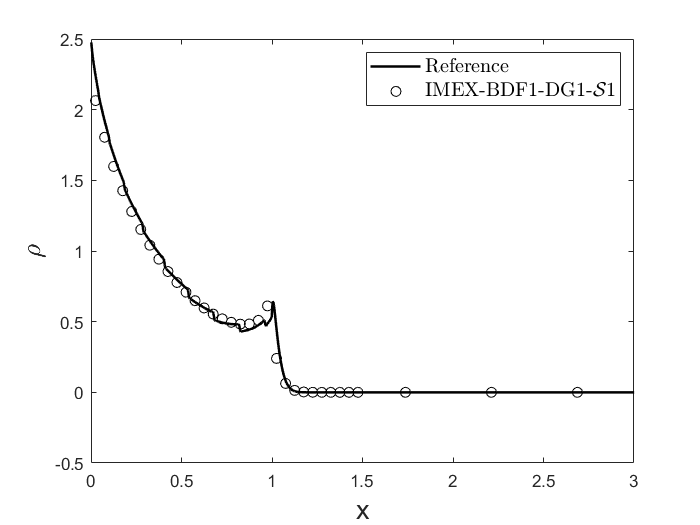}
\end{subfigure}
\begin{subfigure}{0.33\textwidth}
\includegraphics[width=\linewidth, height=4.5cm]{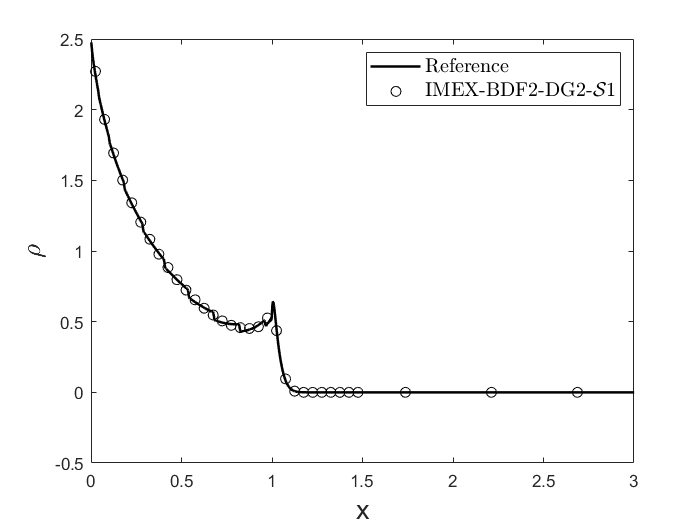}
\end{subfigure}
\begin{subfigure}{0.33\textwidth}
\includegraphics[width=\linewidth, height=4.5cm]{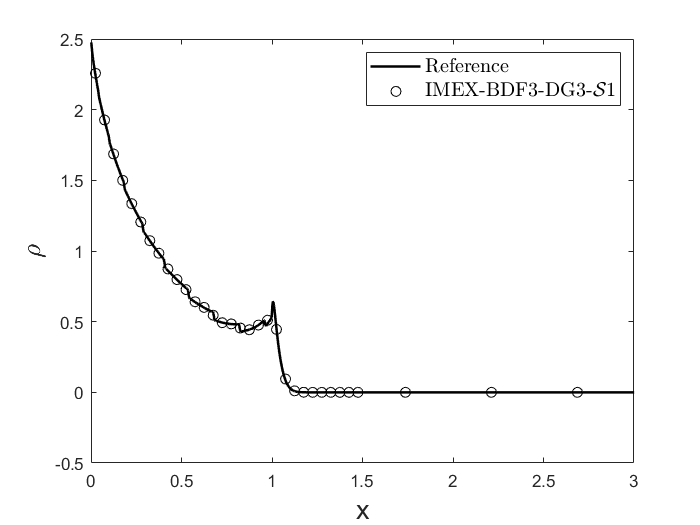}
\end{subfigure}
\begin{subfigure}{0.33\textwidth}
\includegraphics[width=\linewidth, height=4.5cm]{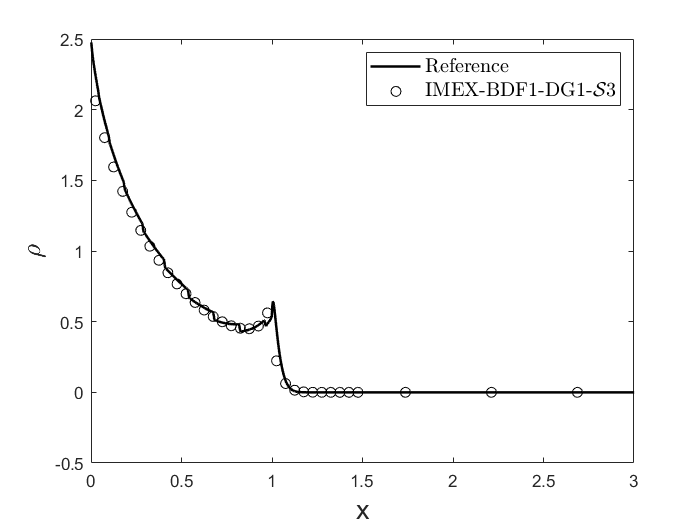}
\end{subfigure}
\begin{subfigure}{0.33\textwidth}
\includegraphics[width=\linewidth, height=4.5cm]{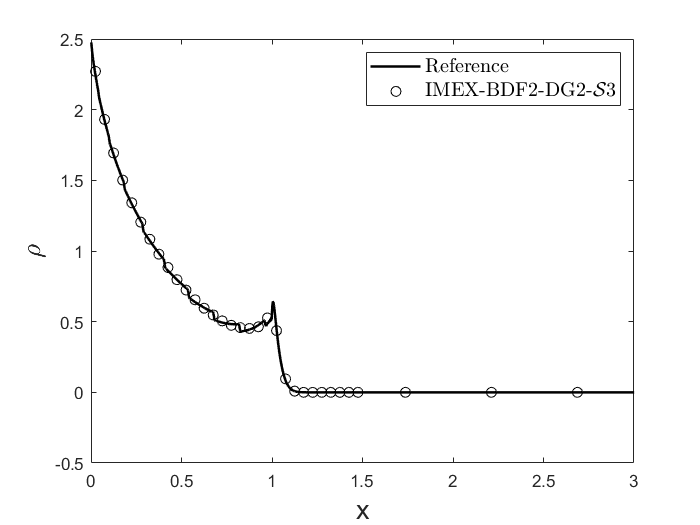}
\end{subfigure}
\begin{subfigure}{0.33\textwidth}
\includegraphics[width=\linewidth, height=4.5cm]
{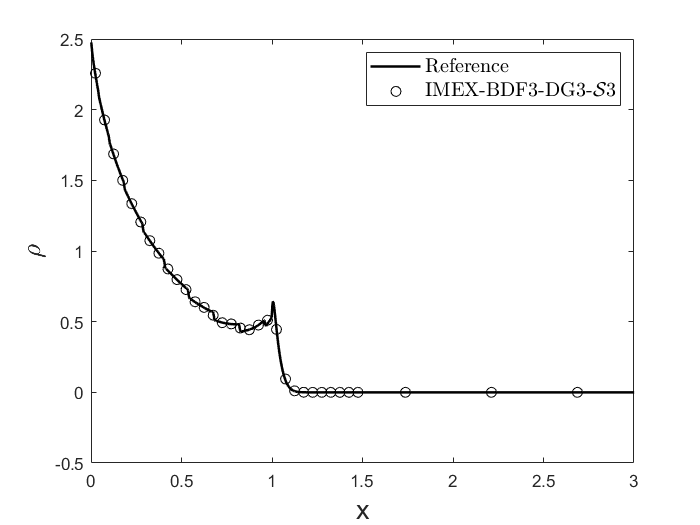}
\end{subfigure}
\caption{Example 3: Two-material problem. The computed $\macro$ at $T = 1.5$.}
\label{fig:two-material epsilon 1 T 15e-1 rho}
\end{figure}

\begin{figure}[ht]
\begin{subfigure}{0.45\textwidth}
\includegraphics[width=\linewidth]
{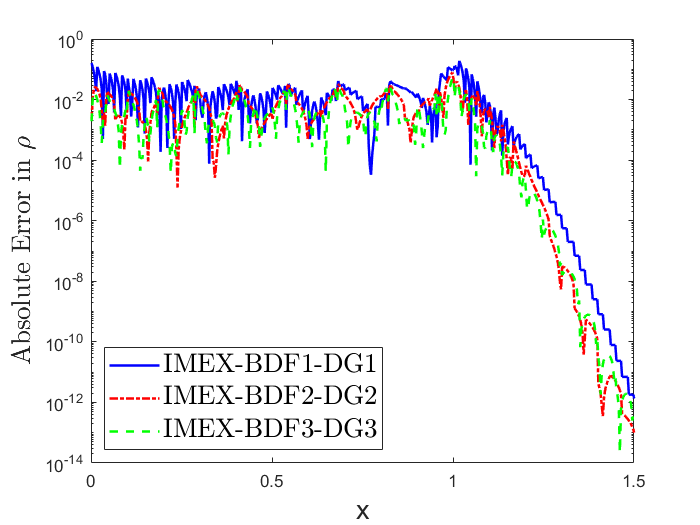}
\end{subfigure}
\begin{subfigure}{0.45\textwidth}
\includegraphics[width=\linewidth]{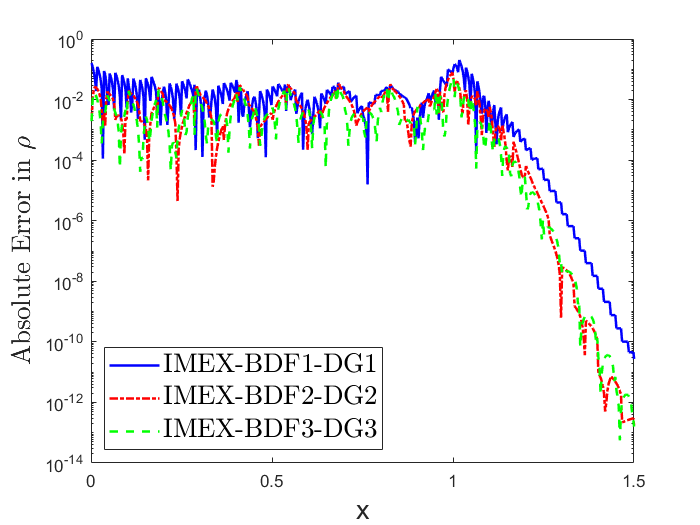}
\end{subfigure}
\caption{{Example 3: Two-material problem. Absolute error in $\macro$ at $T = 1.5$ by methods with Strategy 1 (Left) and Strategy 3 (Right), using  $\dx = 1/60$ on $[0,1.5]$ and $\dx = 19/120$ on $[1.5,11]$ for IMEX-BDF1-DG1, $\dx = 1/30$ on $[0,1.5]$ and $\dx = 19/60$ on $[1.5,11]$ for IMEX-BDF2-DG2, and $\dx = 1/20$ on $[0,1.5]$ and $\dx = 19/40$ on $[1.5,11]$ for IMEX-BDF3-DG3, so all methods have the same number of spatial degrees of freedom. Reference solution uses $\dx = 1/320$ on $[0,1.5]$ and $\dx = 19/640$ on $[1.5,11]$.}}
\label{fig:two-material epsilon 1 T 15e-1 rho error}
\end{figure}

In Figure  \ref{fig:two-material epsilon 1 T 15e-1 rho}, we present the computed $\macro$ by the IMEX-BDF$\IMEXBDForder$-DG$\IMEXBDForder$-$\mS$$k$ scheme with $\IMEXBDForder=1,2,3, k=1,3$ for this multi-scale problem. All results match  the reference solution well, and the second and third order methods produce more accurate solutions than the first order ones as expected. The results by the methods with Strategy 2 are indistinguishable from those by Strategy 1, and they are not included. 
{Finally, in Figure \ref{fig:two-material epsilon 1 T 15e-1 rho error}, the point-wise absolute errors of $\macro$ on the subregion $[0,1.5]$ are plotted for Strategy 1 and 3, when the same 150 total number of spatial degrees of freedom are used in each IMEX-BDF$s$-DG$s$, $s=1, 2, 3$. The results support that it is relatively beneficial to work with higher than first order methods to produce numerical solutions with better resolution, even when the solution has limited regularity.}

\subsection{Example 4: Riemann problem for telegraph equation}
\label{subsec:riemann}
%%%%%%%%%%%%%%%%%%%%%%%

In this example,  a Riemann problem  is considered for the telegraph equation\footnote{{Note that the telegraph equation is different from the one-group transport equation with 2-point Gauss-Legendre quadrature.}} on $\spaceset = [-1,1]$, with $\scat(x)=1$, $\absorp(x)=0$,  zero source, and the following 
initial conditions
\begin{align*}
\macro(x,0)
& =
\begin{cases}
2, & x \leq 0 \\
1, & x > 0
\end{cases}
,
&
\micro(x,\velangle,0) 
& = 
\begin{cases}
0, & x \leq 0 \\
0, & x > 0
\end{cases}.
\end{align*}

{We will focus on the transport regime with $\vareps = 0.7$, particularly to report results on effectively controlling numerical oscillations for formally second and third order methods.  Interested readers can refer to \cite{Matsuda2025RPI} for numerical results when the model is in the diffusive regime and our methods are confirmed to be AP. } There is no visible difference between the computed solutions by the methods using IMEX Strategy 1 and 2 (except some comments below regarding the use of the nonlinear limiter), hence we only present the results with Strategy 2 and 3. 
All the tests are done with a uniform spatial mesh of  $\dx = 1/40$ and a final time $T=0.15$. 
The time step $\dt$ is chosen as follows for the IMEX-BDF$s$-DG$s$ method. With   Strategy 2, we take $\dt=0.5 \vareps \dx + 0.25 \scatlower \dx^2$ when $s=1$, $\dt=0.13 \vareps \dx + 0.013 \scatlower \dx^2$ when $s=2$, and $\dt=0.053 \vareps \dx + 0.0014 \scatlower \dx^2$ when $s=3$; with  Strategy 3, $\dt$ is determined by  \eqref{eq:fourier:dt:st3}.  % 
The reference solution is computed using the first order forward Euler upwind finite difference scheme with $\dx = 1/1000$ 
and $\dt=0.05 \vareps \dx$. {The numerical treatment for the inflow boundary condition follows the same approach as for the example in Section \ref{sec:two-material}, except that no numerical integration is needed to evaluate any integration in $\velangle$.}

\begin{figure}[ht]
\centering{

\includegraphics[width=0.4\textwidth, height=1.7in]{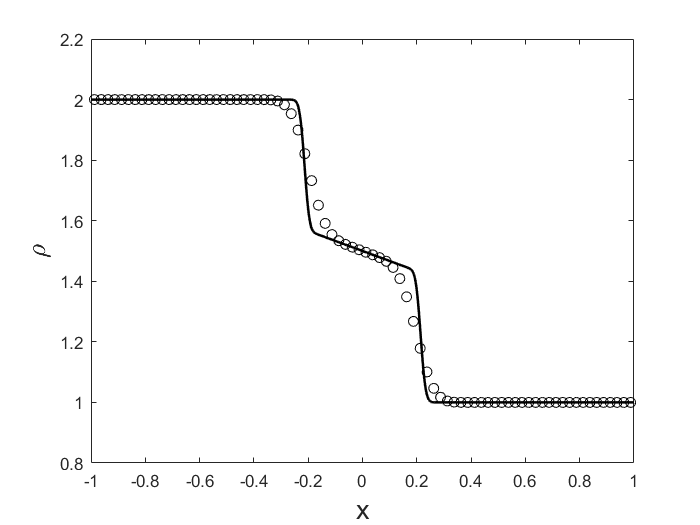}
\includegraphics[width=0.4\textwidth, height=1.7in]{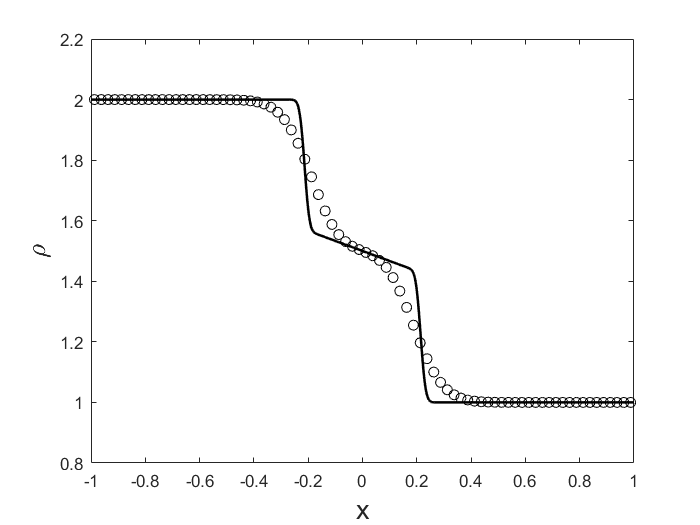}

\includegraphics[width=0.4\textwidth, height=1.7in]{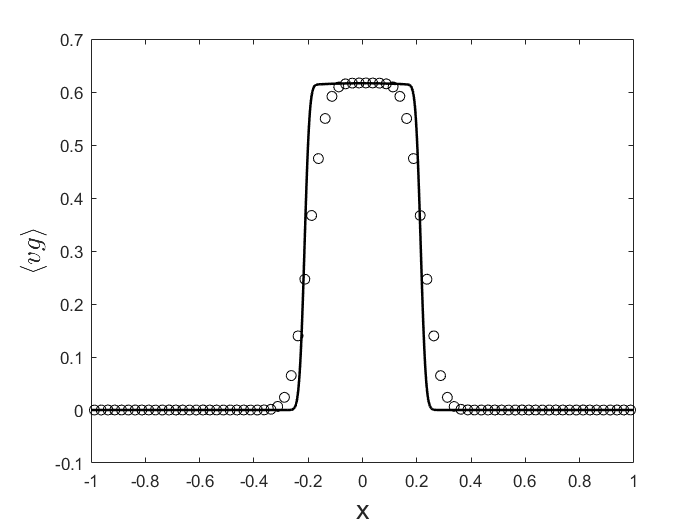}
\includegraphics[width=0.4\textwidth, height=1.7in]{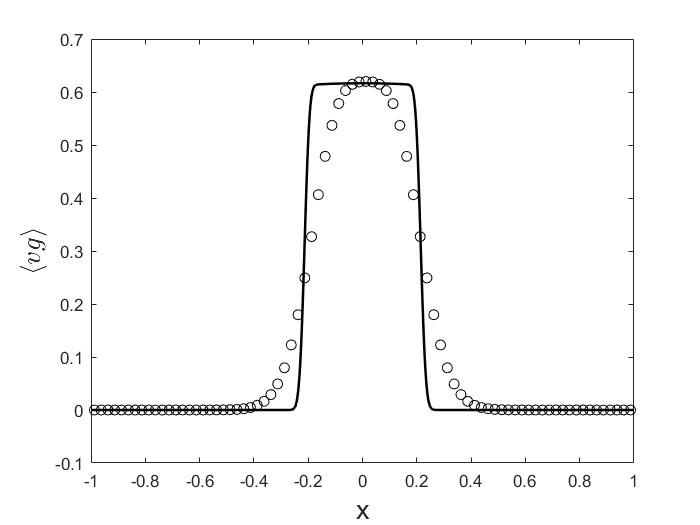}
}

\caption{Example 4: Riemann problem for telegraph equation. $\macro$ (top row) and $\lvg \velangle \micro\rvg$ (bottom row) by the first order IMEX-BDF1-DG1 schemes, with $\vareps = 0.7$, and  Strategy 2 (left column) and 3 (right column).}
\label{fig:Riemann problem-1stOrder-eps07}
\end{figure}

In Figure \ref{fig:Riemann problem-1stOrder-eps07}, both $\macro$ and $\lvg \velangle\micro \rvg$  computed by the first order methods with Strategy 2 and 3  are reported. The computed solutions capture the overall profile of the solutions. They are smeared  around  the discontinuities yet without any oscillation. The non-sharpness is more prominent with IMEX Strategy 3 in  Figure \ref{fig:Riemann problem-1stOrder-eps07} (right column).  

\begin{figure}[h!]
\centering{
\begin{subfigure}{0.4\textwidth}
\includegraphics[width=\linewidth, height=1.7in]{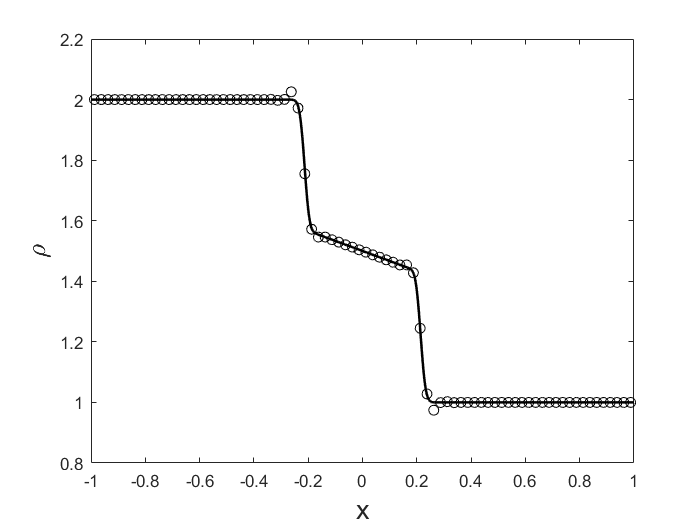}
\caption{$\macro$, IMEX-BDF2-DG2}
\end{subfigure}
\begin{subfigure}{0.4\textwidth}
\includegraphics[width=\linewidth, height=1.7in]{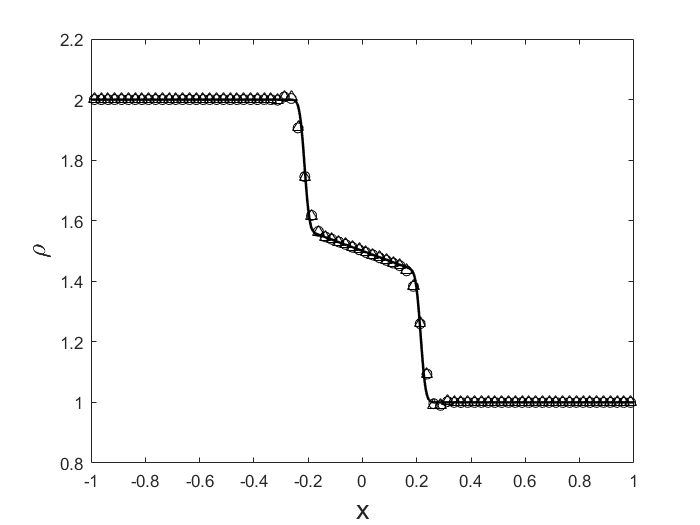}
\caption{$\macro$, IMEX-BDF2-DG2, OE-limiter}
\end{subfigure}

\begin{subfigure}{0.4\textwidth}
\includegraphics[width=\linewidth, height=1.7in]{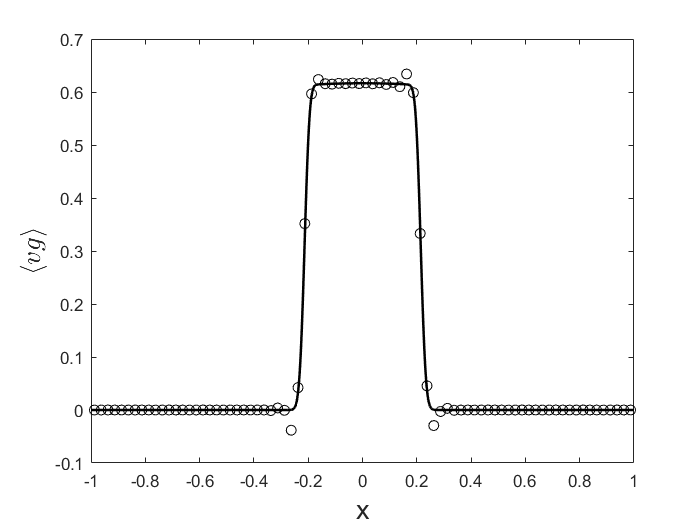}
\caption{$\lvg \velangle \micro \rvg$, IMEX-BDF2-DG2}
\end{subfigure}
\begin{subfigure}{0.4\textwidth}
\includegraphics[width=\linewidth, height=1.7in]{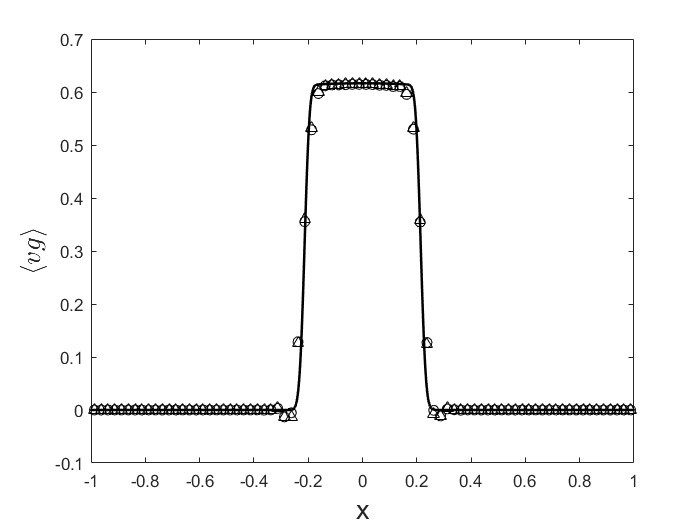}
\caption{$\lvg \velangle \micro \rvg$, IMEX-BDF2-DG2, OE-limiter}
\end{subfigure}

\begin{subfigure}{0.4\textwidth}
\includegraphics[width=\linewidth, height=1.7in]{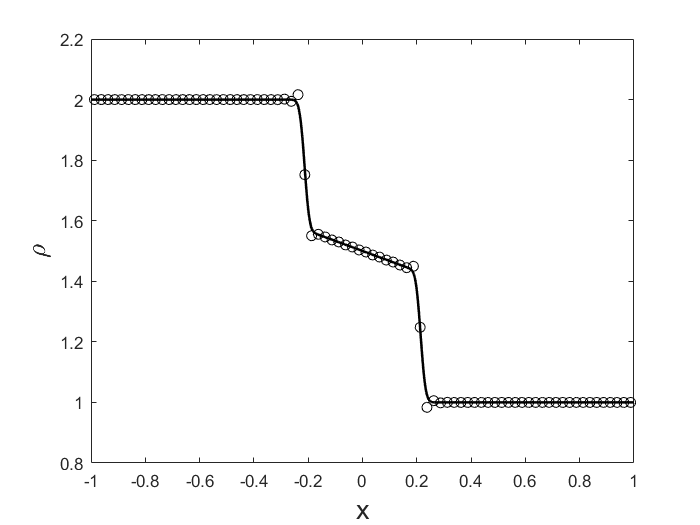}
\caption{$\macro$, IMEX-BDF3-DG3}
\end{subfigure}
\begin{subfigure}{0.4\textwidth}
\includegraphics[width=\linewidth, height=1.7in]{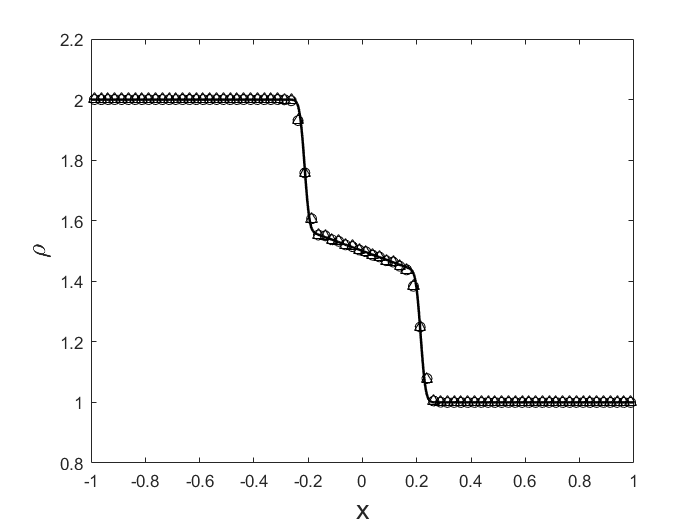}
\caption{$\macro$, IMEX-BDF3-DG3, OE-limiter}
\end{subfigure}

\begin{subfigure}{0.4\textwidth}
\includegraphics[width=\linewidth, height=1.7in]{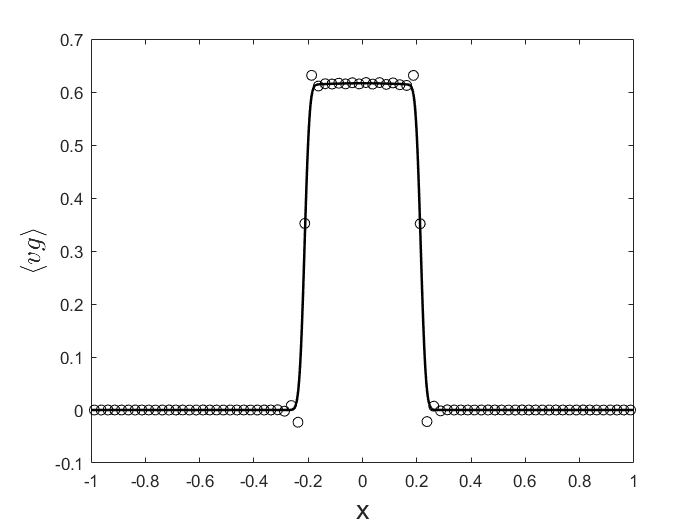}
\caption{$\lvg \velangle \micro \rvg$, IMEX-BDF3-DG3}
\end{subfigure}
\begin{subfigure}{0.4\textwidth}
\includegraphics[width=\linewidth, height=1.7in]{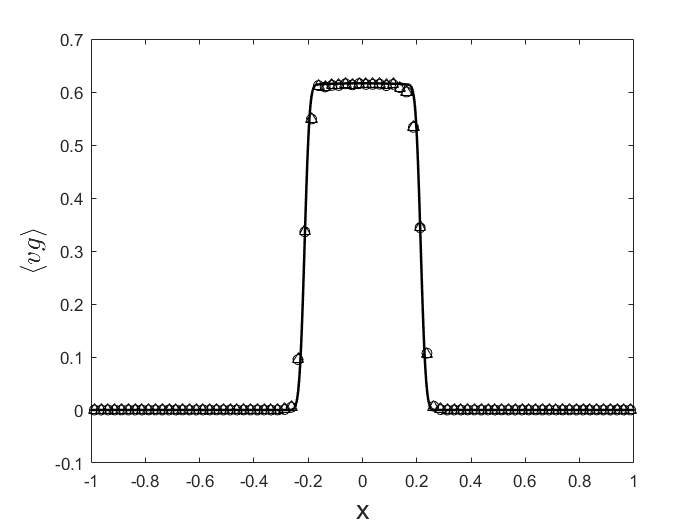}
\caption{$\lvg \velangle \micro \rvg$, IMEX-BDF3-DG3, OE-limiter}
\end{subfigure}
}

\caption{Example 4: Riemann problem for telegraph equation.  $\macro$ and $\lvg \velangle \micro\rvg$  by IMEX-BDF2-DG2 (top two rows) and IMEX-BDF3-DG3 (bottom two rows) with Strategy 2 with $\vareps = 0.7$. Left: no limiter; right: with OE-limiter, applied to $\macro$ before (marked with $\circ$) and after (marked with $\Delta$) solving for $\micro$. 
}
\label{fig:Riemann problem-2nd3rdOrder-eps07-st2}
\end{figure}

\begin{figure}[h!]
\centering{
\begin{subfigure}{0.4\textwidth}
\includegraphics[width=\linewidth, height=1.7in]{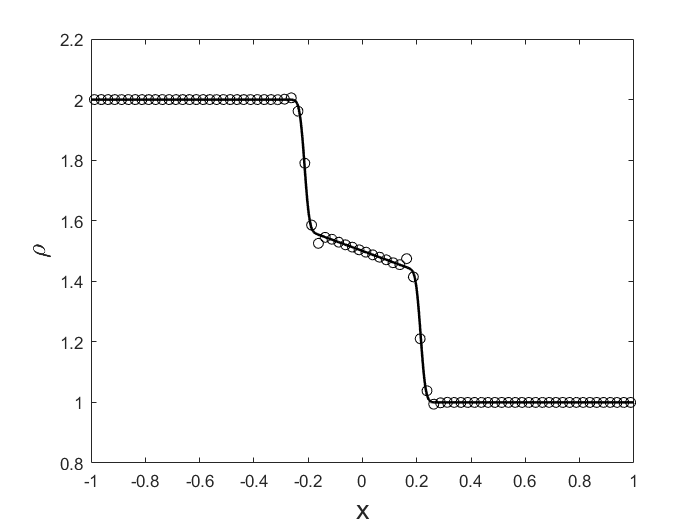}
\caption{$\macro$, IMEX-BDF2-DG2}
\end{subfigure}
\begin{subfigure}{0.4\textwidth}
\includegraphics[width=\linewidth, height=1.7in]{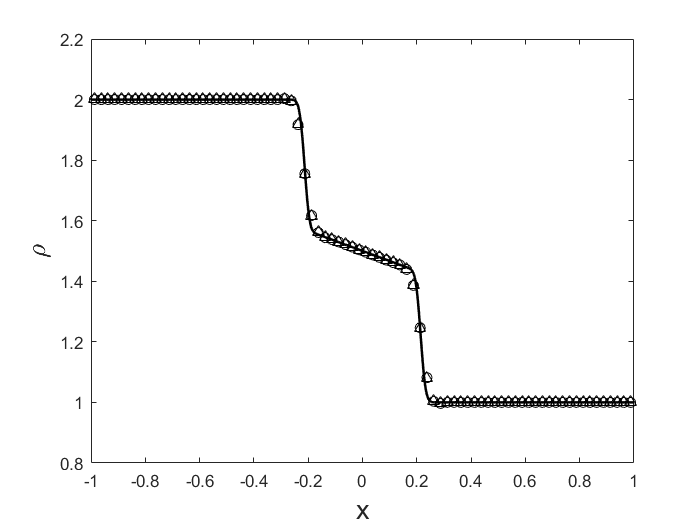}
\caption{$\macro$, IMEX-BDF2-DG2, OE-limiter}
\end{subfigure}

\begin{subfigure}{0.4\textwidth}
\includegraphics[width=\linewidth, height=1.7in]{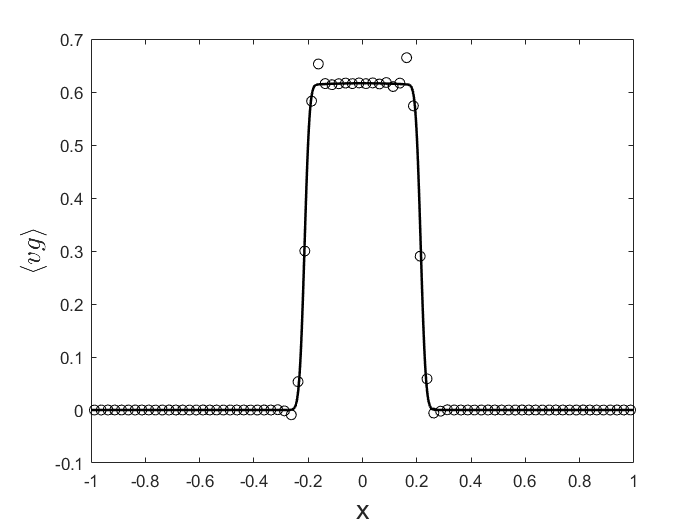}
\caption{$\lvg \velangle \micro \rvg$, IMEX-BDF2-DG2}
\end{subfigure}
\begin{subfigure}{0.4\textwidth}
\includegraphics[width=\linewidth, height=1.7in]{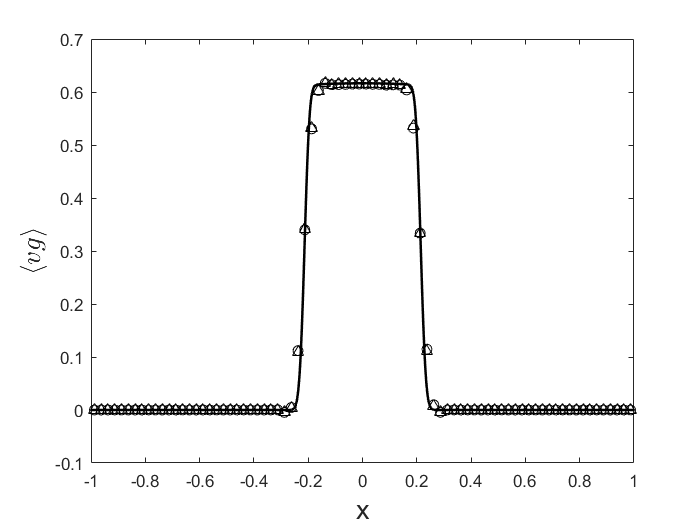}
\caption{$\lvg \velangle \micro \rvg$, IMEX-BDF2-DG2, OE-limiter}
\end{subfigure}

\begin{subfigure}{0.4\textwidth}
\includegraphics[width=\linewidth, height=1.7in]{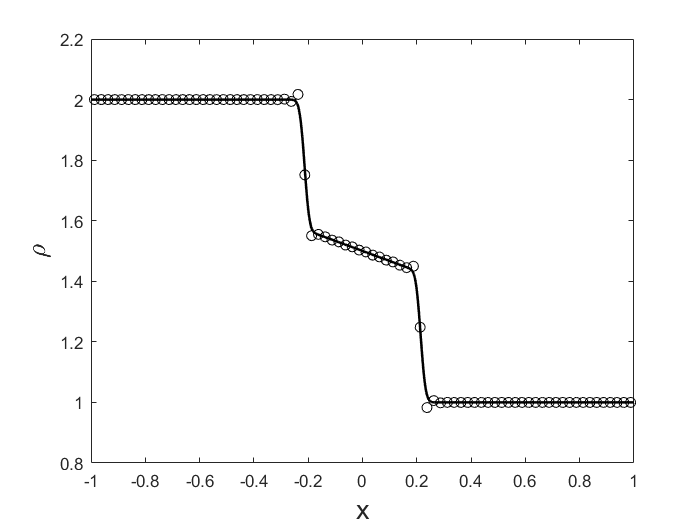}
\caption{$\macro$, IMEX-BDF3-DG3}
\end{subfigure}
\begin{subfigure}{0.4\textwidth}
\includegraphics[width=\linewidth, height=1.7in]{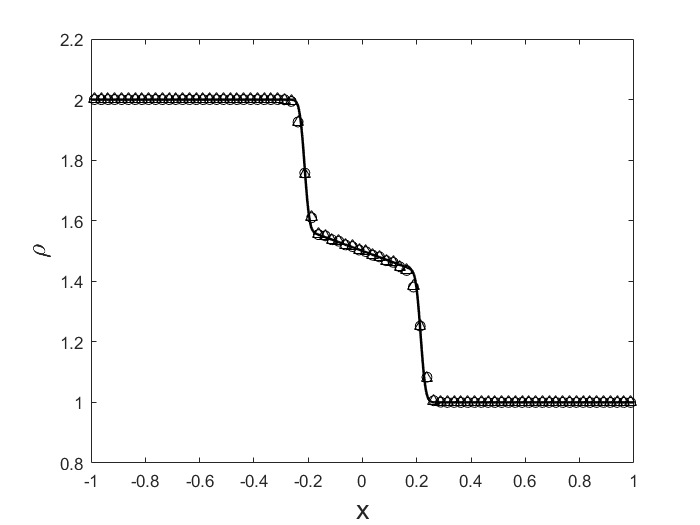}
\caption{$\macro$, IMEX-BDF3-DG3, OE-limiter}
\end{subfigure}

\begin{subfigure}{0.4\textwidth}
\includegraphics[width=\linewidth, height=1.7in]{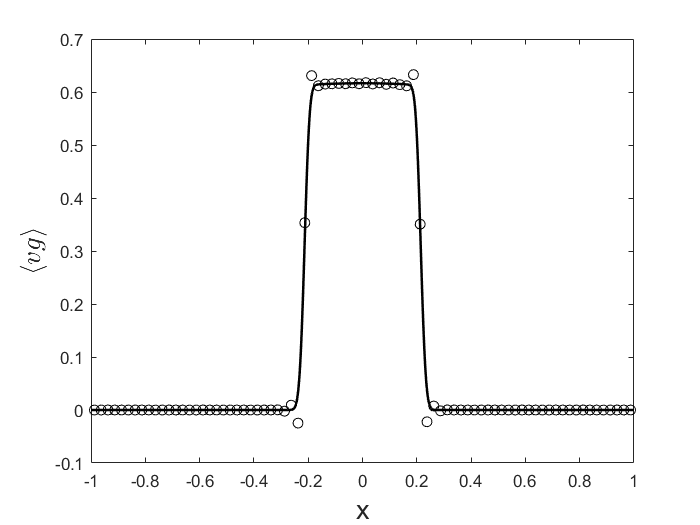}
\caption{$\lvg \velangle \micro \rvg$, IMEX-BDF3-DG3}
\end{subfigure}
\begin{subfigure}{0.4\textwidth}
\includegraphics[width=\linewidth, height=1.7in]{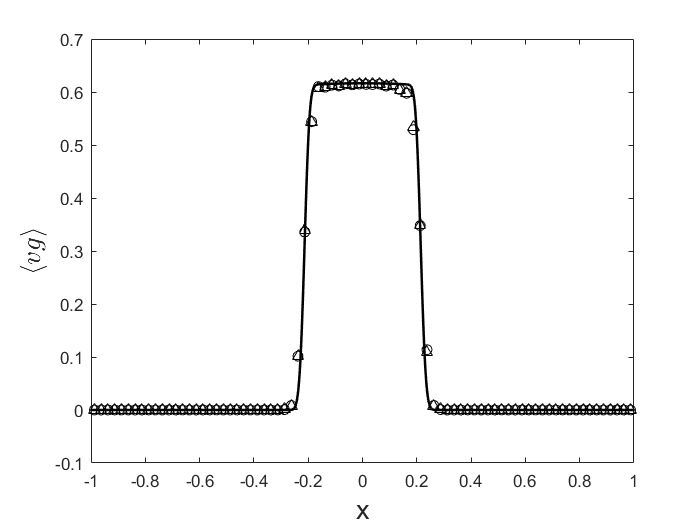}
\caption{$\lvg \velangle \micro \rvg$, IMEX-BDF3-DG3, OE-limiter}
\end{subfigure}
}
\caption{Example 4: Riemann problem for telegraph equation.  $\macro$ and $\lvg \velangle \micro\rvg$  by IMEX-BDF2-DG2 (top two rows) and IMEX-BDF3-DG3 (bottom two rows) with Strategy 3 with $\vareps = 0.7$. Left: no limiter; right: with OE-limiter, applied to $\macro$ before (marked with $\circ$) and after (marked with $\Delta$) solving for $\micro$. 
}
\label{fig:Riemann problem-2nd3rdOrder-eps07-st3}
\end{figure}

In Figure \ref{fig:Riemann problem-2nd3rdOrder-eps07-st2} (left column), we present the results by the second and third order methods with Strategy 2. Compared with the first order methods, the computed solutions better capture the sharp transitions at the discontinuities but display numerical oscillations. 
To control the oscillations, we apply the nonlinear limiting strategy in \cite{peng2024oedg}, developed in the setting of the  oscillation-eliminating discontinuous Galerkin methods. In our simulation, the limiting strategy, referred to as OE-limiter, is implemented  as a post-processing step\footnote{There is a parameter $\beta_j$ in the definition of the strategy in Eqn (2.5)  of \cite{peng2024oedg}. It is the characteristic velocity in the free-streaming operator and is taken as $\beta_j=1/\vareps=10/7$.},
and it is applied (only) to  $\macro$ after each time step of our schemes. 
The results are reported in Figure \ref{fig:Riemann problem-2nd3rdOrder-eps07-st2} (right column).
Note that with IMEX Strategy 2, we solve for $\macro$ first and then solve for $\micro$  over each time step. This gives two possible ways to apply the OE-limiter: to apply it to $\macro$ before or after we solve for  $\micro$. The results by these two  implementations, marked with circles and triangles in Figure \ref{fig:Riemann problem-2nd3rdOrder-eps07-st2}, show little difference, and the numerical oscillations are effectively reduced and nearly eliminated completely, not only in $\macro$ but also in $\lvg \velangle\micro \rvg$.  The tests above are further carried out for the second and third order methods with Strategy 3, see Figure \ref{fig:Riemann problem-2nd3rdOrder-eps07-st3}, and very similar observations can be made. We will end this example with a few comments.
\begin{itemize}
    \item With Strategy 1 in the second and third order methods, we solve for  $\micro$ before solving for $\macro$ over each time step, therefore we only look into one implementation of the OE-limiter, that is to apply it to $\macro$ once it is available. The resulting solutions are comparable with those using Strategy 2.
\item With Strategy 3 and without the OE-limiter, the overshoots around discontinuities of $\lvg \velangle\micro \rvg$ are larger than those by Strategy 2. 
\end{itemize}

\subsection{Example 5: An example with non-well-prepared initial data}
\label{sec:num:non-well-prepared}

The last numerical test is to examine an example with a not well-prepared initial configuration,  namely with $g(x,\velangle,0)=O(1/\vareps)$. We here consider the one-group transport equation in slab geometry on $\spaceset = [0,2\pi]$ with the initial data $f(x,\velangle, 0)=(1+(\velangle-0.5)^2) (1 + 0.05 \cos(x))$, periodic boundary conditions and zero source. It is easy to derive $\macro(x,0)=\frac{19}{12}(1 + 0.05 \cos(x))$. Moreover, $\scat=1$ and $\absorp=0$, and the focus is on the diffusive regime with  $\vareps = 10^{-6}$.  

With the unconditional stability of the proposed methods with Strategy 3, we consider the IMEX-BDF$s$-DG$s$-$\mS$3 methods on uniform spatial meshes with $\dt=0.5\dx$, and present the numerical results in Table \ref{tab:initiallayer} by the (formally) second order and third order methods with $s=2,3$ at the final time $T=1$. The results in the {\bf left-half} panel of Table \ref{tab:initiallayer} are computed when a variation of the first initialization approach is applied. More specifically, the IMEX-RK$s$-DG$s$-$\mS$$k$ scheme with the time step $\dt$ is applied to compute the solutions at $t=t^1, \cdots, t^{s-1}$.  According to the AP analysis in Section \ref{sec:AP} especially Remark \ref{rem:non-well-prep-reduction}, with the non-well-prepared initial condition,  this initialization will result in a temporal error of $O(\dt)$  and lead to the overall error $O(\dt)+O(\dx^s)=O(\dx)+O(\dx^s)$.  This order reduction is observed numerically.

The results in the {\bf right-half} panel of Table \ref{tab:initiallayer} are computed with a modified initialization strategy:  the IMEX-RK$s$-DG$s$-$\mS$3 scheme with a reduced time step, $h^s$, is applied to compute the solutions at $t=t^1, \cdots, t^{s-1}$ for the IMEX-BDF$s$-DG$s$-$\mS$3 scheme.  To avoid using the BDF time integrators with variable time step sizes, the IMEX-RK$s$-DG$s$-$\mS$3 method with the time step $\dt$ is also used to compute the solution at $t^{s},\cdots, t^{2s-2}$ as a transitional phase,  before the IMEX-BDF$s$-DG$s$-$\mS$3 scheme (with $\dt$) is applied to the rest of the time period. This simple initial fix is in the similar spirit of Remark 5.2 in \cite{peng2020stability}. With this, the designed order of accuracy is recovered.

\begin{table}[ht] 
\vspace{0.1in} 
\centering 
\caption{Example 5. Example with a non-well-prepared initial condition.  $L^1$ errors and orders at $T=1$ of IMEX-BDF$s$-DG$s$-$\mS$$3$,  initialized by IMEX-RK$s$-DG$s$-$\mS$3,  $s=2,3$, $\vareps = 10^{-6}$. Left-half panel: uniform time-stepping; right-half  panel: uniform time-stepping with an initial fix.}
\begin{tabular}{c|c|c|c|c|c|c|c|c|c}\hline 
\multirow{2}{*} {} & \multirow{2}{*} {$N$} & \multicolumn{4}{c|}{uniform time-stepping} & \multicolumn{4}{c}{with an initial fix} \\ \cline{3-10} 
& & $E_N^{\macro}$ & order & $E_N^{\lvg \velangle \micro \rvg}$ & order & $E_N^{\macro}$ & order & $E_N^{\lvg \velangle \micro \rvg}$ & order \\ \hline 
\multirow{5}{*}{$s=2$} & 20 & 3.09E-04 & - & 1.03E-04 & - & 7.65E-04 & - & 1.18E-03 \\ 
& 40 & 1.71E-04 & 0.86 & 5.69E-05 & 0.86 & 1.67E-04 & 2.20 & 2.77E-04 & 2.09 \\ 
& 80 & 8.11E-05 & 1.07 & 2.70E-05 & 1.07 & 4.57E-05 & 1.87 & 7.08E-05 & 1.97 \\ 
& 160 & 4.16E-05 & 0.96 & 1.39E-05 & 0.96 & 1.19E-05 & 1.94 & 1.79E-05 & 1.99 \\ 
& 320 & 2.06E-05 & 1.01 & 6.88E-06 & 1.01 & 3.04E-06 & 1.97 & 4.49E-06 & 1.99 \\ 
\hline 
\multirow{5}{*}{$s=3$} & 20 & 1.91E-04 & - & 6.35E-05 & - & 6.30E-04 & - & 1.46E-03 \\ 
& 40 & 9.93E-05 & 0.94 & 3.31E-05 & 0.94 & 6.64E-05 & 3.25 & 1.04E-04 & 3.81 \\ 
& 80 & 4.66E-05 & 1.09 & 1.55E-05 & 1.09 & 9.36E-06 & 2.83 & 1.33E-05 & 2.96 \\ 
& 160 & 2.37E-05 & 0.97 & 7.90E-06 & 0.97 & 1.23E-06 & 2.92 & 1.69E-06 & 2.98 \\ 
& 320 & 1.17E-05 & 1.02 & 3.91E-06 & 1.02 & 1.59E-07 & 2.96 & 2.12E-07 & 2.99 \\ 
\hline 
\end{tabular} 
\label{tab:initiallayer}
\end{table}

%===============
\section{Concluding remarks}
\label{sec:conclusions}

This work presents our recent development  in the design and analysis of high order asymptotic preserving (AP) methods for multi-scale simulation of kinetic transport models, especially in the framework of using discontinuous  Galerkin spatial discretizations based on the micro-macro decomposition of the governing model.  Linear multi-step methods, specifically, the implicit-explicit (IMEX) BDF methods are used as the time integrators, coupled with three different {partitioned strategies}  to delineate which terms are treated explicitly or implicitly. The impact of these strategies on the stability, accuracy, computational complexity,  and AP property is systematically investigated.
The main contributions and findings are summarized below.
\begin{itemize}
\item[-] The stability with respect to the model scales and discretization parameters is established for the proposed methods by  energy-based and Fourier-based analysis, with the former relying on a suitably defined discrete energy for each family of the methods,  and the latter providing the time-step conditions for actual simulations. The Fourier analysis has utilized some intrinsic structure {analytically} identified for the amplification matrix.  Stability analysis also improves the understanding of the role of the absorption term and its numerical treatments.
\item[-] {In the transport regime with $\vareps=O(1)$, all the proposed methods require a hyperbolic type time step condition for stability. When coming to the diffusive regime,  methods with Strategy 1 and 2 require a more stringent diffusion type time step condition, while methods with  Strategy 3 are unconditionally stable and allow much larger time step sizes. The relative efficiency of these methods would depend  on the specific examples (e.g., whether temporal errors dominate) as well as on how the linear system \eqref{eq:alg-6} encountered by Strategy 3 is numerically solved. }
\item[-] With IMEX Strategy 3, the condition number is analyzed for the discrete diffusive matrix, to be inverted per time step, in terms of the model and discretization parameters. Such estimate informs about the computational complexity of the respective methods.
\item[-] A numerical  comparison is made to understand the cost efficiency of using IMEX-BDF and  (certain type)  IMEX Runge-Kutta time integrators, while other ingredients to define the methods stay the same. It is shown that the methods with IMEX-BDF in time are more cost efficient to generate numerical solutions with {higher resolution. } 
\item[-] Unlike the IMEX Runge-Kutta time integrators  coupled with certain IMEX partitionings (i.e. Strategy 1), 
the methods with the IMEX-BDF in time {are not observed numerically} to suffer from order reduction in the stiff case  with all three IMEX partitionings.  {More theoretical analysis would be needed to fully understand this.}

\item[-] Two initialization approaches are examined.  One is to use one-step time integrators, and the other  is based on an idea outlined in many textbooks, namely using the linear multi-step methods of lower order accuracy in the same family with increasing order. We find the former is more straightforward to use, while the latter requires careful implementation to work well especially in consideration of the accuracy.

\item[-] With well-prepared initial conditions, the proposed methods are AP. Particularly, as $\vareps\rightarrow 0$, the limiting schemes with Strategy 1 and 2 are explicit discretizations of the limiting diffusion equation \eqref{eq:diff-limit}, and those with Strategy 3 are  implicit discretizations of \eqref{eq:diff-limit}.  When the initial condition is not well-prepared, the methods with Strategy 1 and 3 are still AP, however the methods with Strategy 2  alone are not, due to their incapability  to correctly exit the initial layer. Fortunately, methods with  Strategy 2 will become AP if the initialization is done by utilizing methods with Strategy 1 or 3. 

\item[-] In the presence of discontinuities in the solution, the oscillation-eliminating limiting strategy \cite{peng2024oedg}, applied only to the macroscopic $\macro$ as a post-processing step {after each time step,}  is effective to control spurious  oscillations in kinetic simulations. 
\end{itemize}

{What deserves further investigation is the numerical boundary treatment of methods formulated based on the micro–macro decomposition of the model, for example, to preserve the designed order of accuracy for smooth solutions with inflow Dirichlet boundary conditions, or to ensure the correct behavior of numerical solutions near domain boundaries under anisotropic boundary conditions in the intermediate and diffusive regimes \cite{lemou2012micro}.}   The methodologies along with the mathematical understanding established here can be further used or developed for multi-scale simulation of full radiative transfer equations and other kinetic transport models.

\bibliographystyle{plain}
\bibliography{references.bib}

\begin{thebibliography}{10}

\bibitem{adams2001discontinuous}
Marvin~L Adams.
\newblock Discontinuous finite element transport solutions in thick diffusive
  problems.
\newblock {\em Nuclear science and engineering}, 137(3):298--333, 2001.

\bibitem{adams2002fast}
Marvin~L Adams and Edward~W Larsen.
\newblock Fast iterative methods for discrete-ordinates particle transport
  calculations.
\newblock {\em Progress in nuclear energy}, 40(1):3--159, 2002.

\bibitem{ascher1995implicit}
Uri~M Ascher, Steven~J Ruuth, and Brian~TR Wetton.
\newblock Implicit-explicit methods for time-dependent partial differential
  equations.
\newblock {\em SIAM Journal on Numerical Analysis}, 32(3):797--823, 1995.

\bibitem{boscarino2013implicit}
Sebastiano Boscarino, Lorenzo Pareschi, and Giovanni Russo.
\newblock Implicit-explicit {R}unge-{K}utta schemes for hyperbolic systems and
  kinetic equations in the diffusion limit.
\newblock {\em SIAM Journal on Scientific Computing}, 35(1):A22--A51, 2013.

\bibitem{bourgat1994coupling}
JF~Bourgat, Patrick Le~Tallec, B~Perthame, and Y~Qiu.
\newblock Coupling {B}oltzmann and {E}uler equations without overlapping.
\newblock {\em Contemporary Mathematics}, 157:377--398, 1994.

\bibitem{crestetto2019asymptotically}
Ana{\"\i}s Crestetto, Nicolas Crouseilles, Giacomo Dimarco, and Mohammed Lemou.
\newblock Asymptotically complexity diminishing schemes ({ACDS}) for kinetic
  equations in the diffusive scaling.
\newblock {\em Journal of Computational Physics}, 394:243--262, 2019.

\bibitem{degond2005smooth}
Pierre Degond and Shi Jin.
\newblock A smooth transition model between kinetic and diffusion equations.
\newblock {\em SIAM Journal on Numerical Analysis}, 42(6):2671--2687, 2005.

\bibitem{dimarco2017implicit}
Giacomo Dimarco and Lorenzo Pareschi.
\newblock Implicit-explicit linear multistep methods for stiff kinetic
  equations.
\newblock {\em SIAM Journal on Numerical Analysis}, 55(2):664--690, 2017.

\bibitem{dimarco2018asymptotic}
Giacomo Dimarco, Lorenzo Pareschi, and Giovanni Samaey.
\newblock Asymptotic-preserving {M}onte {C}arlo methods for transport equations
  in the diffusive limit.
\newblock {\em SIAM Journal on Scientific Computing}, 40(1):A504--A528, 2018.

\bibitem{einkemmer2021asymptotic}
Lukas Einkemmer, Jingwei Hu, and Yubo Wang.
\newblock An asymptotic-preserving dynamical low-rank method for the
  multi-scale multi-dimensional linear transport equation.
\newblock {\em Journal of Computational Physics}, 439:110353, 2021.

\bibitem{filbet2015hierarchy}
Francis Filbet and Thomas Rey.
\newblock A hierarchy of hybrid numerical methods for multiscale kinetic
  equations.
\newblock {\em SIAM Journal on Scientific Computing}, 37(3):A1218--A1247, 2015.

\bibitem{golse1999convergence}
Fran{\c{c}}ois Golse, Shi Jin, and C~David Levermore.
\newblock The convergence of numerical transfer schemes in diffusive regimes i:
  {D}iscrete-ordinate method.
\newblock {\em SIAM Journal on Numerical Analysis}, 36(5):1333--1369, 1999.

\bibitem{guermond2010asymptotic}
Jean-Luc Guermond and Guido Kanschat.
\newblock Asymptotic analysis of upwind discontinuous {G}alerkin approximation
  of the radiative transport equation in the diffusive limit.
\newblock {\em SIAM Journal on Numerical Analysis}, 48(1):53--78, 2010.

\bibitem{guo2022low}
Wei Guo and Jing-Mei Qiu.
\newblock A low rank tensor representation of linear transport and nonlinear
  {V}lasov solutions and their associated flow maps.
\newblock {\em Journal of Computational Physics}, 458:111089, 2022.

\bibitem{hairer1987solving}
Ernst Hairer, Syvert~Paul N{\o}rsett, and Gerhard Wanner.
\newblock {\em Solving Ordinary Differential Equations {I}: Nonstiff Problems}.
\newblock Springer-Verlag, 1987.

\bibitem{hesthaven2007nodal}
Jan~S Hesthaven and Tim Warburton.
\newblock {\em Nodal discontinuous Galerkin methods: algorithms, analysis, and
  applications}.
\newblock Springer Science \& Business Media, 2007.

\bibitem{jang2014analysis}
Juhi Jang, Fengyan Li, Jing-Mei Qiu, and Tao Xiong.
\newblock Analysis of asymptotic preserving {DG-IMEX} schemes for linear
  kinetic transport equations in a diffusive scaling.
\newblock {\em SIAM Journal on Numerical Analysis}, 52(4):2048--2072, 2014.

\bibitem{jang2015high}
Juhi Jang, Fengyan Li, Jing-Mei Qiu, and Tao Xiong.
\newblock High order asymptotic preserving {DG-IMEX} schemes for
  discrete-velocity kinetic equations in a diffusive scaling.
\newblock {\em Journal of Computational Physics}, 281:199--224, 2015.

\bibitem{jin2010asymptotic}
Shi Jin.
\newblock Asymptotic preserving ({AP}) schemes for multiscale kinetic and
  hyperbolic equations: a review.
\newblock {\em Lecture notes for summer school on methods and models of kinetic
  theory (M\&MKT), Porto Ercole (Grosseto, Italy)}, pages 177--216, 2010.

\bibitem{jin2000uniformly}
Shi Jin, Lorenzo Pareschi, and Giuseppe Toscani.
\newblock Uniformly accurate diffusive relaxation schemes for multiscale
  transport equations.
\newblock {\em SIAM Journal on Numerical Analysis}, 38(3):913--936, 2000.

\bibitem{larsen1983numerical}
Edward~W Larsen.
\newblock On numerical solutions of transport problems in the diffusion limit.
\newblock {\em Nuclear Science and Engineering}, 83(1):90--99, 1983.

\bibitem{larsen1987asymptotic}
Edward~W Larsen, Jim~E Morel, and Warren~F Miller~Jr.
\newblock Asymptotic solutions of numerical transport problems in optically
  thick, diffusive regimes.
\newblock {\em Journal of Computational Physics}, 69(2):283--324, 1987.

\bibitem{lemou2012micro}
Mohammed Lemou and Florian M{\'e}hats.
\newblock Micro-macro schemes for kinetic equations including boundary layers.
\newblock {\em SIAM Journal on Scientific Computing}, 34(6):B734--B760, 2012.

\bibitem{lemou2008new}
Mohammed Lemou and Luc Mieussens.
\newblock A new asymptotic preserving scheme based on micro-macro formulation
  for linear kinetic equations in the diffusion limit.
\newblock {\em SIAM Journal on Scientific Computing}, 31(1):334--368, 2008.

\bibitem{lewis1984computational}
Elmer~Eugene Lewis and Warren~F Miller.
\newblock {\em Computational Methods of Neutron Transport}.
\newblock John Wiley and Sons, Inc., New York, NY, 1984.

\bibitem{liu2004boltzmann}
Tai-Ping Liu and Shih-Hsien Yu.
\newblock Boltzmann equation: micro-macro decompositions and positivity of
  shock profiles.
\newblock {\em Communications in mathematical physics}, 246(1):133--179, 2004.

\bibitem{Matsuda2025RPI}
Kimberly Matsuda.
\newblock {\em Multi-scale simulation and model order reduction for the
  radiative transfer equation}.
\newblock PhD thesis, Rensselaer Polytechnic Institute, 2025.

\bibitem{pareschi2011efficient}
Lorenzo Pareschi and Giovanni Russo.
\newblock Efficient asymptotic preserving deterministic methods for the
  {B}oltzmann equation.
\newblock {\em Models and Computational Methods for Rarefied Flows, AVT-194 RTO
  AVT/VKI, Rhode St. Genese, Belgium}, page~34, 2011.

\bibitem{peng2024oedg}
Manting Peng, Zheng Sun, and Kailiang Wu.
\newblock {OEDG}: Oscillation-eliminating discontinuous {G}alerkin method for
  hyperbolic conservation laws.
\newblock {\em Mathematics of Computation}, 2024.

\bibitem{peng2020structure}
Zhichao Peng.
\newblock {\em Structure-preserving discontinuous Galerkin methods for
  multi-scale kinetic transport equations and nonlinear optics models}.
\newblock PhD thesis, Rensselaer Polytechnic Institute, 2020.

\bibitem{peng2020asymptotic}
Zhichao Peng, Vrushali~A Bokil, Yingda Cheng, and Fengyan Li.
\newblock Asymptotic and positivity preserving methods for {K}err-{D}ebye model
  with {L}orentz dispersion in one dimension.
\newblock {\em Journal of Computational Physics}, 402:109101, 2020.

\bibitem{peng2024micro}
Zhichao Peng, Yanlai Chen, Yingda Cheng, and Fengyan Li.
\newblock A micro-macro decomposed reduced basis method for the time-dependent
  radiative transfer equation.
\newblock {\em Multiscale Modeling \& Simulation}, 22(1):639--666, 2024.

\bibitem{peng2020stability}
Zhichao Peng, Yingda Cheng, Jing-Mei Qiu, and Fengyan Li.
\newblock Stability-enhanced {AP} {IMEX}-{LDG} schemes for linear kinetic
  transport equations under a diffusive scaling.
\newblock {\em Journal of Computational Physics}, 415:109485, 2020.

\bibitem{peng2021stability}
Zhichao Peng, Yingda Cheng, Jing-Mei Qiu, and Fengyan Li.
\newblock Stability-enhanced {AP} {IMEX}1-{LDG} method: energy-based stability
  and rigorous {AP} property.
\newblock {\em SIAM Journal on Numerical Analysis}, 59(2):925--954, 2021.

\bibitem{peng2021asymptotic}
Zhichao Peng and Fengyan Li.
\newblock Asymptotic preserving {IMEX-DG-S} schemes for linear kinetic
  transport equations based on {S}chur complement.
\newblock {\em SIAM Journal on Scientific Computing}, 43(2):A1194--A1220, 2021.

\bibitem{reyna2015operator}
Matthew~A Reyna and Fengyan Li.
\newblock Operator bounds and time step conditions for the {DG} and central
  {DG} methods.
\newblock {\em Journal of Scientific Computing}, 62(2):532--554, 2015.

\bibitem{tencer2016reduced}
John Tencer, Kevin Carlberg, Roy Hogan, and Marvin Larsen.
\newblock Reduced order modeling applied to the discrete ordinates method for
  radiation heat transfer in participating media.
\newblock In {\em Heat Transfer Summer Conference}, volume 50336, page
  V002T15A011. American Society of Mechanical Engineers, 2016.

\bibitem{wang2008variable}
Dong Wang and Steven~J Ruuth.
\newblock Variable step-size implicit-explicit linear multistep methods for
  time-dependent partial differential equations.
\newblock {\em Journal of Computational Mathematics}, pages 838--855, 2008.

\bibitem{xiong2015high}
Tao Xiong, Juhi Jang, Fengyan Li, and Jing-Mei Qiu.
\newblock High order asymptotic preserving nodal discontinuous {G}alerkin
  {IMEX} schemes for the {BGK} equation.
\newblock {\em Journal of Computational Physics}, 284:70--94, 2015.

\bibitem{xu20192}
Yuan Xu, Qiang Zhang, Chi-wang Shu, and Haijin Wang.
\newblock The {L}2-norm stability analysis of {R}unge--{K}utta discontinuous
  {G}alerkin methods for linear hyperbolic equations.
\newblock {\em SIAM Journal on Numerical Analysis}, 57(4):1574--1601, 2019.

\bibitem{yan2021adaptive}
David Yan, Mary~C Pugh, and Francis~P Dawson.
\newblock Adaptive time-stepping schemes for the solution of the
  {P}oisson-{N}ernst-{P}lanck equations.
\newblock {\em Applied Numerical Mathematics}, 163:254--269, 2021.

\bibitem{zhang2023asymptotic}
Guoliang Zhang, Hongqiang Zhu, and Tao Xiong.
\newblock Asymptotic preserving and uniformly unconditionally stable finite
  difference schemes for kinetic transport equations.
\newblock {\em SIAM Journal on Scientific Computing}, 45(5):B697--B730, 2023.

\end{thebibliography}

%-------

%%%%%%%%%%%%%%%%%%%
\begin{appendices}
%%%%%%%%%%%%%%%%%%

\counterwithin{equation}{section}

\section{{An example of non-AP implicit methods}}
\label{sec:appendA}

\newcommand{\beq}{\begin{equation}}
\newcommand{\eeq}{\end{equation}}
\newcommand{\al}{{\alpha}}   
\newcommand{\be}{{\beta}}   
\newcommand{\df}{\partial}

{With $\al(x, t)=f(x, v=1, t)$, $\be(x,t)=f(x, v=-1, t)$, the telegraph equation is given as 
\beq
\label{APP:eq:1}
\eps\df_t \al+\df_x \al=\frac{1}{2\eps}(\be-\al), \quad 
\eps\df_t \be-\df_x \be=\frac{1}{2\eps}(\al-\be). 
\eeq
By introducing ${ \rho}=\frac{\al+\be}{2}=\langle f\rangle$, ${ J}=\frac{\al-\be}{2\eps}=\frac{1}{\eps} \langle vf\rangle$, the model \eqref{APP:eq:1} is further rewritten as  
$\partial_t\rho+\partial_x J=0$, $ \eps^2\partial_tJ+\partial_x\rho=-J.$
Particularly when $\eps\ll 1$,  one formally has a diffusive model for $\rho$,
\beq
\label{APP:eq:2}
\partial_t\rho=\df_{xx}\rho +O(\eps^2).
\eeq 
}

{
For \eqref{APP:eq:1}, we now apply  the first order upwind finite difference (also the $P^0$ upwind discontinuous Galerkin) method in space with the backward Euler method in time, 
\begin{subequations}
\label{APP:eq:3}
\begin{align}
		\eps\frac{\al^{n+1}_j - \al^n_j}{\Delta t} &+  \frac{ \al^{n+1}_{j}- \al^{n+1}_{j-1}} {\Delta x}=\frac{1}{2\eps}(\be^{n+1}_j-\al^{n+1}_j),\\		
		\eps\frac{\be^{n+1}_j - \be^n_j}{\Delta t} &- \frac{ \be^{n+1}_{j+1}- \be^{n+1}_{j}} {\Delta x}=\frac{1}{2\eps}(\al^{n+1}_j-\be^{n+1}_j).
\end{align}
\end{subequations}
Here $\Delta x$ and $\Delta t$ denote the meshsize in space and time, respectively, and $\al^n_j\approx \al(x_j, t^n)$, $\be^n_j\approx \be(x_j, t^n)$.  The scheme \eqref{APP:eq:3}, in terms of $\rho$ and $J$, is equivalent to 
\begin{subequations}
\label{APP:eq:4}
\begin{align}
		\frac{\rho^{n+1}_j - \rho^n_j}{\Delta t} &+ \frac{J^{n+1}_{j+1}- J^{n+1}_{j-1}} {2\Delta x} 
		 - { \frac{{\Delta x}}{2\eps} \frac{\rho^{n+1}_{j+1} - 2\rho^{n+1}_j +\rho^{n+1}_{j-1}}{{\Delta x}^2}}= 0,\\	
\eps^2\frac{J^{n+1}_j - J^n_j}{\Delta t}&+ \frac{\rho^{n+1}_{j+1} - \rho^{n+1}_{j-1}}{2\Delta x}-{ \frac{\eps {\Delta x} }{2} \frac{J^{n+1}_{j+1}-2J^{n+1}_j + J^{n+1}_{j-1}}{{\Delta x}^2}}= -J^{n+1}_j.
\end{align}
\end{subequations}
Formally when $\eps\ll 1$, one gets 
\begin{align}
\label{APP:eq:5}
\frac{\rho^{n+1}_j - \rho^n_j}{\Delta t}
=\frac{\rho^{n+1}_{j+2}-2\rho^{n+1}_j + \rho^{n+1}_{j-2}}{(2\Delta x)^2}
		 +{  \frac{{\Delta x}}{2\eps} \frac{\rho^{n+1}_{j+1} - 2\rho^{n+1}_j +\rho^{n+1}_{j-1}}{{\Delta x}^2}}+O(\eps).
\end{align}
Particularly on a fixed mesh and with  $\eps\ll 1$,   the scheme \eqref{APP:eq:5} is inconsistent to the diffusive model in \eqref{APP:eq:2}, due to the numerical dissipation related to the upwind treatment. This implies  that the first order method in \eqref{APP:eq:3}, though being implicit, is not asymptotic preserving (AP), and it cannot faithfully capture the diffusive limit on under-resolved meshes. 
}

%%%%%%%%%%%%%%%%%%%%%%%%%%%
\section{Proof of Theorem \ref{thm:stab:energy}}
\label{sec:stab:energy:proof}
%%%%%%%%%%%%%%%%%%%%%%%%%%%

Note that the first order IMEX-BDF temporal scheme in this work is also a first order IMEX-Runge Kutta (RK) scheme, and our proposed IMEX-BDF1-DG$r$-$\mS$2 scheme is essentially the same as the DG$r$-IMEX1 scheme in \cite{jang2015high,jang2014analysis}, while our IMEX-BDF1-DG$r$-$\mS$3 scheme is the same as the IMEX1-DG$r$-S scheme in \cite{peng2021asymptotic}. Hence the energy-based stability analysis for the DG1-IMEX1 scheme in \cite{jang2014analysis}, that assumes $\absorp(x)\equiv0$ and the velocity variable is continuous, 
can be adapted to the setting here and  leads to \eqref{eq:st2:stab-ene} with \eqref{eq:st1:dt}-\eqref{eq:st1:dt.a}. {With very minor modification, the energy-based stability analysis for the IMEX1-DG1-S scheme in \cite{peng2021asymptotic} will give \eqref{eq:st3:stab}, a slightly more refined result.} 
With these, we here only present the analysis for the IMEX-BDF1-DG1-$\mS$1 scheme. To this end, a similar proof procedure as in \cite{jang2014analysis} is followed, together with some technique in \cite{peng2020structure} to estimate the explicit discretization for the  absorption terms (see  \eqref{eq:stab:5.1} below).

We take $\testmacro = \macro_h^{n+1}$ in \eqref{eq:BDF1:DGr:st1:w.a} and $\testmicro =  \micro_{h,q}^{n+1}$ in \eqref{eq:BDF1:DGr:st1:w.b},  multiply \eqref{eq:BDF1:DGr:st1:w.b} by $\vweight_q$ and sum over all $q = 1, 2, \dots, N_{\velangle}$, and get
%
% Energy analysis: test functions = macro, micro equations
%----------------------------------------------------------------------------------------------------------------
\begin{subequations} \label{eq:stab:1}
\begin{align} 
& \frac{1}{\dt} (\macro_h^{n+1} - \macro_h^n, \macro_h^{n+1})
+ l_h (\lvg \velangle \micro_h^{n+1} \rvg_h, \macro_h^{n+1})
\notag \\
& = \frac{1}{2 \dt} (|| \macro_h^{n+1} ||^2 - || \macro_h^n ||^2 + || \macro_h^{n+1} - \macro_h^n ||^2) 
+ l_h (\lvg \velangle \micro_h^{n+1} \rvg_h, \macro_h^{n+1})= -(\absorp \macro_h^n, \macro_h^{n+1}),\label{eq:stab:1.a}\\
& \frac{\vareps^2}{\dt} \lvg (\micro_h^{n+1} - \micro_h^n, \micro_h^{n+1}) \rvg_h + \vareps \lvg b_{h,\velangle} (\micro_h^n, \micro_h^{n+1}) \rvg_h - \lvg \velangle d_h(\macro_h^n, \micro_h^{n+1}) \rvg_h 
\notag \\
& = \frac{\vareps^2}{2 \dt} (||| \micro_h^{n+1} |||^2 - ||| \micro_h^n |||^2 + ||| \micro_h^{n+1} - \micro_h^n |||^2)
+ \vareps \lvg b_{h,\velangle} (\micro_h^n, \micro_h^{n+1}) \rvg_h
- \lvg \velangle d_h(\macro_h^n, \micro_h^{n+1}) \rvg_h 
\notag \\
& = -||| \micro_h^{n+1} |||_s^2
- \vareps^2 \lvg (\absorp \micro_h^n, \micro_h^{n+1}) \rvg_h.
\label{eq:stab:1.b}
\end{align}
\end{subequations}
Using the definition in \eqref{eq:bilinear:a} and the property in \eqref{eq:adjoint}, one has
%----------------------------------------------------------------------------------------------------------------
% Energy analysis: l_h(m,n) = <v d_h(n,m)>
%----------------------------------------------------------------------------------------------------------------
\begin{equation} \label{eq:stab:2}
l_h(\lvg \velangle \micro_h^{n+1} \rvg_h, \macro_h^{n+1}) =  d_h(\macro_h^{n+1}, \lvg \velangle \micro_h^{n+1}\rvg_h) =\lvg \velangle d_h(\macro_h^{n+1}, \micro_h^{n+1}) \rvg_h.
\end{equation}
We now combine \eqref{eq:stab:1}-\eqref{eq:stab:2} and reach
%----------------------------------------------------------------------------------------------------------------
% Energy analysis: Added equations
%----------------------------------------------------------------------------------------------------------------
\begin{align}  \label{eq:stab:3}
& \frac{1}{2 \dt} (|| \macro_h^{n+1} ||^2 + \vareps^2 ||| \micro_h^{n+1} |||^2 - || \macro_h^n ||^2 - \vareps^2 ||| \micro_h^n |||^2) 
+ \frac{1}{2 \dt} (|| \macro_h^{n+1} - \macro_h^n ||^2  + \vareps^2 ||| \micro_h^{n+1} - \micro_h^n |||^2)
\notag \\
& + \vareps \lvg b_{h,\velangle}(\micro_h^{n+1}, \micro_h^{n+1}) \rvg_h - \vareps \lvg b_{h,\velangle}(\micro_h^{n+1} - \micro_h^{n}, \micro_h^{n+1}) \rvg_h + \lvg \velangle d_h(\macro_h^{n+1} - \macro_h^n, \micro_h^{n+1}) \rvg_h \notag \\
& = -||| \micro_h^{n+1} |||_s^2 - (\absorp \macro_h^n, \macro_h^{n+1}) - \vareps^2 \lvg (\absorp \micro_h^n, \micro_h^{n+1}) \rvg_h.    
\end{align}

Since $\pdspace^0$ is the space of all piecewise constant functions, one has $\partial_x \micro_h^{n+1}|_{I_j}=0,\forall j$.
Using this and the property \eqref{eq:micro v integral zero property}, following some similar  steps in \cite{jang2014analysis} (e.g.  those to show (3.22)-(3.25) there),  we obtain
%----------------------------------------------------------------------------------------------------------------
% Energy analysis: Bilinear form inequalities
%----------------------------------------------------------------------------------------------------------------
\begin{subequations}
\label{eq:stab:4}
\begin{align} 
\lvg b_{h,\velangle}(\micro_h^{n+1}, \micro_h^{n+1}) \rvg_h 
& = \left \lvg \frac{\left| v \right|}{2} \sum_j [\micro_h^{n+1}]_{j-\frac{1}{2}}^2 \right \rvg_h,
\label{eq:stab:4.a} \\
\left| \lvg b_{h,\velangle}(\micro_h^{n+1} - \micro_h^n, \micro_h^{n+1}) \rvg_h \right|
& \leq \eta_1 ||| \micro_h^{n+1} - \micro_h^n |||^2 + \frac{|| \velangle ||_{h,\infty}}{2 \eta_1 \dx} \left \lvg \frac{\left| v \right|}{2} \sum_j [\micro_h^{n+1}]_{j-\frac{1}{2}}^2 \right \rvg_h,
\label{eq:stab:4.b} \\
\left| \lvg \velangle d_h(\macro_h^{n+1} - \macro_h^n, \micro_h^{n+1} \rvg_h \right|
& \leq \eta_2 \| \macro_h^{n+1} - \macro_h^n ||^2 + \frac{\lvg \left| \velangle \right| \rvg_h}{2 \eta_2 \dx} \left \lvg \frac{\left| v \right|}{2} \sum_j [\micro_h^{n+1}]_{j-\frac{1}{2}}^2 \right \rvg_h, 
\label{eq:stab:4.c}
\end{align}
\end{subequations}
where the positive quantities $\eta_1$ and $\eta_2$ will be specified later. 
With these estimates,
\eqref{eq:stab:3} becomes
%----------------------------------------------------------------------------------------------------------------
% Energy analysis: Added equations with bilinear inequalities
%----------------------------------------------------------------------------------------------------------------
\begin{align} \label{eq:stab:5}
& \frac{1}{2 \dt} (|| \macro_h^{n+1} ||^2 + \vareps^2 ||| \micro_h^{n+1} |||^2 - || \macro_h^n ||^2 - \vareps^2 ||| \micro_h^n |||^2) 
+ \Big(\frac{1}{2 \dt} - \eta_2\Big) || \macro_h^{n+1} - \macro_h^n ||^2 
\notag \\
& + \Big(\frac{\vareps^2}{2 \dt} - \vareps \eta_1\Big) ||| \micro_h^{n+1} - \micro_h^n |||^2
+ \Big(\vareps - \frac{\vareps || \velangle ||_{h,\infty}}{2 \eta_1 \dx} - \frac{\lvg \left| \velangle \right| \rvg_h}{2 \eta_2 \dx}\Big) \left \lvg \frac{\left| v \right|}{2} \sum_j [\micro_h^{n+1}]_{j-\frac{1}{2}}^2 \right \rvg_h
\notag \\
& \leq -||| \micro_h^{n+1} |||_s^2 - (\absorp \macro_h^n, \macro_h^{n+1}) - \vareps^2 \lvg (\absorp \micro_h^n, \micro_h^{n+1}) \rvg_h.    
\end{align}
To estimate the absorption terms, we apply  a simple yet not so obvious estimate in Remark 3.6.2 of \cite{peng2020structure}, together with the upper bound of $\absorp(x)$ in \eqref{eq:sig:bound}, and get
\begin{subequations}
\label{eq:stab:5.1}
\begin{align}
-(\absorp\macro_h^n, \macro_h^{n+1}) &\leq \frac{1}{4}|| \macro_h^{n+1}-\macro_h^n ||_a^2\leq\frac{\absorpupper}{4}|| \macro_h^{n+1}-\macro_h^n ||^2, \label{eq:stab:5.1a}\\
-\lvg (\absorp\micro_h^n, \micro_h^{n+1}) \rvg_h &\leq \frac{1}{4} ||| \micro_h^{n+1}-\micro_h^n |||_a^2\leq \frac{\absorpupper}{4}||| \micro_h^{n+1}-\micro_h^n |||^2.
\end{align}
\end{subequations}
The first inequality in \eqref{eq:stab:5.1a} is nothing but $|| \macro_h^{n+1}+\macro_h^n ||_a^2\geq 0$. Now combining \eqref{eq:stab:5}-\eqref{eq:stab:5.1} and the lower bound of $\scat(x)$ in \eqref{eq:sig:bound}, we have
%----------------------------------------------------------------------------------------------------------------
% Energy analysis: Added equations with bilinear inequalities and scat, absorp inequalities
%----------------------------------------------------------------------------------------------------------------
\begin{align} \label{eq:stab:6} 
& \frac{1}{2\dt} E_{h,\mS 1}^{n+1} - \frac{1}{2\dt} E_{h,\mS 1}^n
+ {\Big[\frac{1}{2} \Big(\frac{1}{\dt} - \frac{\absorpupper}{2}\Big) - \eta_2\Big]}
|| \macro_h^{n+1} - \macro_h^n ||^2
\notag \\
& + \Big[\frac{\vareps^2}{2} \Big(\frac{1}{\dt} - \frac{\absorpupper}{2}\Big) - \vareps \eta_1\Big] ||| \micro_h^{n+1} - \micro_h^n |||^2
+ \Big(\vareps - \frac{\vareps || \velangle ||_{h,\infty}}{2 \eta_1 \dx} - \frac{\lvg \left| \velangle \right| \rvg_h}{2 \eta_2 \dx}\Big) \left \lvg \frac{\left| v \right|}{2} \sum_j [\micro_h^{n+1}]_{j-\frac{1}{2}}^2 \right \rvg_h
\notag \\
& \leq -\scatlower ||| \micro_h^{n+1} |||^2.
\end{align}

What remains is to identify the time step condition that  ensures $E_{h,\mS 1}^{n+1} \leq E_{h,\mS 1}^n$, hence  gives the stability.
 Firstly, we require
\begin{equation}
\label{eq:cond1:dt}
    \dt< {\frac{2}{\absorpupper}, }
    \quad \text{if}\;\;\; \absorpupper\ne 0.
\end{equation}
Next by using (3.35) in \cite{jang2014analysis}, namely, $\left \lvg \frac{\left| v \right|}{2} \sum_j [\micro_h^{n+1}]_{j-\frac{1}{2}}^2 \right \rvg_h \leq \frac{2 || \velangle ||_{h,\infty}}{\dx} ||| \micro_h^{n+1} |||^2$,
\eqref{eq:stab:6}  becomes
%----------------------------------------------------------------------------------------------------------------
% Energy analysis: Added equations with bilinear inequalities and scat, absorp, micro square jump inequalities
%----------------------------------------------------------------------------------------------------------------
\begin{align} 
 \frac{1}{2\dt} E_{h,\mS 1}^{n+1} & - \frac{1}{2\dt} E_{h,\mS 1}^n
+ \Big[\frac{1}{2} \Big(\frac{1}{\dt} - \frac{\absorpupper}{2}\Big) - \eta_2\Big] || \macro_h^{n+1} - \macro_h^n ||^2
\notag \\
& + \Big[\frac{\vareps^2}{2} \Big(\frac{1}{\dt} - \frac{\absorpupper}{2}\Big) - \vareps \eta_1\Big] ||| \micro_h^{n+1} - \micro_h^n |||^2
\notag \\
& \leq -\Big[\scatlower+\min\big(\vareps - \frac{\vareps || \velangle ||_{h,\infty}}{2 \eta_1 \dx} - \frac{\lvg \left| \velangle \right| \rvg_h}{2 \eta_2 \dx},0\big)\; \frac{2 || \velangle ||_{h,\infty}}{h} \Big] ||| \micro_h^{n+1} |||^2.
\end{align}
From here, $E_{h,\mS 1}^{n+1} \leq E_{h,\mS 1}^n$ can be ensured if one further requires 
%----------------------------------------------------------------------------------------------------------------
% Energy analysis: Nonnegative coefficient requirement
%----------------------------------------------------------------------------------------------------------------
\begin{subequations} \label{eq:stab:7}
\begin{align}
&\frac{1}{2} \Big(\frac{1}{\dt} - \frac{\absorpupper}{2}\Big) - \eta_2  \geq 0,\quad
\frac{\vareps^2}{2} \Big(\frac{1}{\dt} - \frac{\absorpupper}{2}\Big) - \vareps \eta_1  \geq 0,
\label{eq:stab:7.a} \\
&\scatlower+\min\Big(\vareps - \frac{\vareps || \velangle ||_{h,\infty}}{2 \eta_1 \dx} - \frac{\lvg \left| \velangle \right| \rvg_h}{2 \eta_2 \dx},0\Big)\; \frac{2 || \velangle ||_{h,\infty}}{h}  
\geq 0.
\label{eq:stab:7.b}
\end{align}
\end{subequations}
By setting $\eta_2 = \frac{\eta_1}{\vareps}=
 \frac{1}{2} (\frac{1}{\dt} - \frac{\absorpupper}{2})$,  
\eqref{eq:stab:7.a} is satisfied, and  \eqref{eq:stab:7.b} can be simplified to the time step condition 
\eqref{eq:st1:dt}-\eqref{eq:st1:dt.a}. As the final step, one can check and confirm these conditions are compatible with 
the requirement in \eqref{eq:cond1:dt}.

%%%%%%%%%%%%%%%%%%%%%%%%%%%%%%%%%%%%%%%%%%
\section{Proof of Theorem \ref{theorem:similar matrix}}%
\label{sec:stab:fourier:proof}
%%%%%%%%%%%%%%%%%%%%%%%%%%%%%%%%%%%%%%%%%%

%----------------------------------------------------------------------------------------------------------------
% Fourier analysis: Proof of similarity matrix theorem
%----------------------------------------------------------------------------------------------------------------

With similarity, the proof is only carried out for the IMEX-BDF$\IMEXBDForder$-DG$\pd$-$\mS$1 scheme in \eqref{eq:fully:st1}. 

\medskip
{\bf Step 1: Derivation of \eqref{eq:fourier:2}.} Taking the ansatz in \eqref{eq:fourier:ansatz} for the numerical solutions, we obtain the matrix form \eqref{eq:fourier:2} of the scheme \eqref{eq:fully:st1},

%----------------------------------------------------------------------------------------------------------------
% Fourier analysis: Strategy 2 matrix equation
%----------------------------------------------------------------------------------------------------------------
\begin{multline} \label{eq:fourier:3}
\underbrace{
\begin{bmatrix}
G_{L,\IMEXBDForder} \\
& \iden_{\pd(N_{\velangle}+1)} \\
& & \ddots \\
& & & \iden_{\pd(N_{\velangle}+1)}
\end{bmatrix}
}_{\text{$\mathcal{G}_L$}}
\underbrace{\begin{bmatrix}
\mathbf{V}^{n+\IMEXBDForder} \\
\mathbf{V}^{n+\IMEXBDForder-1} \\
\vdots \\
\mathbf{V}^{n+1}
\end{bmatrix}}_{\mathbf{W}^{n+\IMEXBDForder}}
=
\underbrace{
\begin{bmatrix}
G_{R,n+\IMEXBDForder-1} & \dots & G_{R,n+1} & G_{R,n} \\
\iden_{\pd(N_{\velangle}+1)} \\
& \ddots \\
& & \iden_{\pd(N_{\velangle}+1)}
\end{bmatrix}
}_{\text{$\mathcal{G}_R$}}
\underbrace{\begin{bmatrix}
\mathbf{V}^{n+\IMEXBDForder-1} \\
\mathbf{V}^{n+\IMEXBDForder-2} \\
\vdots \\
\mathbf{V}^{n}
\end{bmatrix}}_{\mathbf{W}^{n+\IMEXBDForder-1}},
\end{multline}
with the amplification matrix $\mathbf{G}^{(\IMEXBDForder,\pd)} = \mathbf{G}^{(\IMEXBDForder,\pd)}(\vareps,\scatlower,\dx,\dt;\xi) = \text{$\mathcal{G}_L$}^{-1} \text{$\mathcal{G}_R$}$.
Here, $\iden_m$ is the $m \times m$ identity matrix, $\mathbf{V}^n = [\widehat{\boldsymbol{\macro}}^n, \widehat{\mathbf{\micro}}_1^n, \widehat{\mathbf{\micro}}_2^n, \dots, \widehat{\mathbf{\micro}}_{N_{\velangle}}^n]^T$, and $G_{L,\IMEXBDForder}, G_{R,n+k}$, $k = 0, 1, \dots, \IMEXBDForder - 1$, are $\big(\pd(N_{\velangle}+1)\big)\times\big(\pd(N_{\velangle}+1)\big)$ matrices defined as follows: 
%----------------------------------------------------------------------------------------------------------------
% Fourier analysis: GL, GR
%----------------------------------------------------------------------------------------------------------------
\begin{subequations}
\label{eq:GL:GRk}
\begin{align}
& G_{L,\IMEXBDForder} = 
\begin{bmatrix}
\dx \Mhat & \dt \IMEXBDFimpcoef_{\IMEXBDForder} \vweight_1 \velangle_1 \Dplushat & \dt \IMEXBDFimpcoef_{\IMEXBDForder} \vweight_2 \velangle_2 \Dplushat & \dots & \dt \IMEXBDFimpcoef_{\IMEXBDForder} \vweight_{N_{\velangle}} \velangle_{N_{\velangle}} \Dplushat \\
0 & \dx (\vareps^2 + \dt \IMEXBDFimpcoef_{\IMEXBDForder} \scatlower) \Mhat & 0 & \dots & 0 \\
0 & 0 & \dx (\vareps^2 + \dt \IMEXBDFimpcoef_{\IMEXBDForder} \scatlower) \Mhat & \dots & 0 \\
\vdots & \vdots & \vdots & \ddots & \vdots \\
0 & 0 & 0 & \dots & \dx (\vareps^2 + \dt \IMEXBDFimpcoef_{\IMEXBDForder} \scatlower) \Mhat
\end{bmatrix}
\label{eq:GL}
\\
& G_{R,n+k} =
\notag
\\
& 
\begin{bmatrix}
\IMEXBDFprevcoef_k \dx \Mhat & 0 & 0 & \dots & 0 \\
-\velangle_1 \dt \IMEXBDFexpcoef_k \Dminushat & \IMEXBDFprevcoef_k \vareps^2 \dx \Mhat - \vareps \dt \IMEXBDFexpcoef_k \widehat{UP}_1 & \vareps \dt \IMEXBDFexpcoef_k \Phat_2 & \dots & \vareps \dt \IMEXBDFexpcoef_k \Phat_{N_{\velangle}} \\
-\velangle_2 \dt \IMEXBDFexpcoef_k \Dminushat & \vareps \dt b_k \Phat_1 & \IMEXBDFprevcoef_k \vareps^2 \dx \Mhat - \vareps \dt \IMEXBDFexpcoef_k \widehat{UP}_2 & \dots & \vareps \dt \IMEXBDFexpcoef_k \Phat_{N_{\velangle}} \\
\vdots & \vdots & \vdots & \ddots & \vdots \\
-\velangle_{N_{\velangle}} \dt \IMEXBDFexpcoef_k \Dminushat & \vareps \dt \IMEXBDFexpcoef_k \Phat_1 & \vareps \dt \IMEXBDFexpcoef_k \Phat_2 & \dots & \IMEXBDFprevcoef_k \vareps^2 \dx \Mhat - \vareps \dt \IMEXBDFexpcoef_k \widehat{UP}_{N_{\velangle}}
\end{bmatrix}.
\label{eq:GRk}
\end{align}
\end{subequations}
Here $\Mhat, \Dminushat(\xi), \Dplushat(\xi), \Uhat_q, \Phat_q$, $q = 1, 2, \dots, N_{\velangle}$, are $\pd \times \pd$ matrices, with their $(i,j)$-th entry defined as 
%----------------------------------------------------------------------------------------------------------------
% Fourier analysis: GL, GR element matrices
%----------------------------------------------------------------------------------------------------------------
\begin{subequations}
\begin{align}
(\Mhat)_{ij} & = \frac{1}{2} \int_{-1}^{1} \Legendre_j(x) \Legendre_i(x) dx, \\
(\Dminushat(\xi))_{ij} & = -\int_{-1}^{1} \Legendre_j(x) \partial_x \Legendre_i(x) dx + \Legendre_j(1) \Legendre_i(1) - \exp(-\imag \xi) \Legendre_j(1) \Legendre_i(-1), \\
(\Dplushat(\xi))_{ij} & = -\int_{-1}^{1} \Legendre_j(x) \partial_x \Legendre_i(x) dx + \exp(\imag \xi) \Legendre_j(-1) \Legendre_i(1) - \Legendre_j(-1) \Legendre_i(-1), \\
(\Uhat_q(\xi))_{ij} 
& =\velangle_q \Big(\indicator_{\{\velangle_q \geq 0\}}(\velangle_q) (\Dminushat(\xi))_{ij} + \indicator_{\{\velangle_q < 0\}}(\velangle_q) (\Dplushat(\xi))_{ij}\Big),
\\
(\Phat_q(\xi))_{ij} & = \vweight_q (\Uhat_q(\xi))_{ij}, \quad (\widehat{UP}_q(\xi))_{ij} = (\Uhat_q(\xi))_{ij} -  (\Phat_q(\xi))_{ij},
\end{align}
\end{subequations}
where $\xi = \kappa \dx$ is the discrete wave number and $\indicator_{S}(y)$ is the indicator function of the set $S$.

We further introduce some block matrices, namely, 
$\Mhatblock  = \textrm{diag}(\Mhat,\dots,\Mhat) \in \mathbb{R}^{N_{rv} \times N_{rv}}$, 
$\Dplushatblock=[ 
\vweight_1 \velangle_1, \vweight_2 \velangle_2, \cdots, \vweight_{N_{\velangle}} \velangle_{N_{\velangle}}]
\otimes
\Dplushat
\in \mathbb{R}^{\pd \times N_{rv}} $,
$\Dminushatblock
 =
-
[
\velangle_1,  \velangle_2, \cdots,  \velangle_{N_{\velangle}}]^T\otimes\Dminushat\in \mathbb{R}^{N_{rv}\times \pd}
$, and 
$\UPhatblock
=
\textrm{diag}(\Uhat_1, \Uhat_2, \dots, \Uhat_{N_v})-
[1,1,\cdots,1]^T
\otimes
[\Phat_1, \Phat_2, \cdots, \Phat_{N_{\velangle}}]
\in \mathbb{R}^{N_{rv} \times N_{rv}}.
$
Here $\otimes$ is the Kronecker product and $N_{rv}=\pd N_{\velangle}$. With these, we have \eqref{eq:GL:GRk} more compactly as
%----------------------------------------------------------------------------------------------------------------
% Fourier analysis: GL, GR block matrix form
%----------------------------------------------------------------------------------------------------------------
\begin{equation}
\label{eq:GL:GRk:compact}
G_{L,\IMEXBDForder}
 =
\begin{bmatrix}
\dx \Mhat & \dt \IMEXBDFimpcoef_{\IMEXBDForder} \Dplushatblock \\
0 & \dx (\vareps^2 + \dt \IMEXBDFimpcoef_{\IMEXBDForder} \scatlower) \Mhatblock
\end{bmatrix},\quad
G_{R,n+k} 
 =
\begin{bmatrix}
\IMEXBDFprevcoef_k \dx \Mhat & 0 \\
\dt \IMEXBDFexpcoef_k \Dminushatblock & \IMEXBDFprevcoef_k \vareps^2 \dx \Mhatblock - \vareps \dt \IMEXBDFexpcoef_k \UPhatblock
\end{bmatrix}.
\end{equation}

\medskip
{\bf Step 2: Similarity transformation to $\mathbf{G}^{(\IMEXBDForder,\pd)}$.} We will show that $\mathbf{G}^{(\IMEXBDForder,\pd)}(\vareps,\scatlower,\dx,\dt;\xi)$ is similar to some  matrix $\mathbf{\hat{G}}^{(\IMEXBDForder,\pd)}(\frac{\vareps}{\scatlower \dx}, \frac{\dt}{\vareps \dx}; \xi)$, by following the analysis in \cite{peng2021asymptotic} (see the proof for Theorem 5.6). 
To this end, we introduce
\begin{equation}
\label{eq:fourier:4}
S_{\dx}
= 
\begin{bmatrix} 
\scatlower \dx \iden_{\pd} & 0 \\
0 & \iden_{N_{rv}}
\end{bmatrix},\;\;
T_{\dx}= 
\begin{bmatrix}
\frac{\scatlower}{\vareps} \iden_{\pd} & 0 \\
0 & \frac{1}{\vareps^2\dx} \iden_{N_{rv}} 
\end{bmatrix}
\in \mathbb{R}^{\big(\pd(N_{\velangle}+1)\big) \times \big(\pd(N_{\velangle}+1)\big)}.
\end{equation}
With a direct calculation, we have $\forall k = 0, 1, \dots, \IMEXBDForder-1$,
\begin{align*} 
&S_{\dx}^{-1} G_{L,\IMEXBDForder}^{-1} G_{R,n+k} S_{\dx}=S_{\dx}^{-1} G_{L,\IMEXBDForder}^{-1} T_{\dx}^{-1}S_{\dx} S_{\dx}^{-1}  T_{\dx} G_{R,n+k} S_{\dx}\\
=&S_{\dx}^{-1} \begin{bmatrix}
\frac{\scatlower \dx}{\vareps} \Mhat & \frac{\scatlower\dt}{\vareps } \IMEXBDFimpcoef_{\IMEXBDForder} \Dplushatblock \\
0 & (1+\frac{\scatlower\dt}{\vareps^2} \IMEXBDFimpcoef_{\IMEXBDForder}) \Mhatblock
\end{bmatrix}
^{-1}S_{\dx}  S_{\dx}^{-1}
\begin{bmatrix}
\frac{\scatlower \dx}{\vareps} \IMEXBDFprevcoef_k \Mhat & 0 \\
\frac{\dt}{\vareps^2 \dx} \IMEXBDFexpcoef_k \Dminushatblock & \IMEXBDFprevcoef_k \Mhatblock -\frac{\dt}{\vareps \dx} \IMEXBDFexpcoef_k \UPhatblock
\end{bmatrix} S_{\dx}\\
=&\begin{bmatrix}
\frac{\scatlower \dx}{\vareps} \Mhat & \frac{\dt}{\vareps \dx} \IMEXBDFimpcoef_{\IMEXBDForder} \Dplushatblock \\
0 & (1+(\frac{\scatlower \dx}{\vareps}) (\frac{\dt}{\vareps \dx}) \IMEXBDFimpcoef_{\IMEXBDForder}) \Mhatblock
\end{bmatrix}
^{-1}
\begin{bmatrix}
\frac{\scatlower \dx}{\vareps} \IMEXBDFprevcoef_k \Mhat & 0 \\
(\frac{\scatlower \dx}{\vareps}) (\frac{\dt}{\vareps \dx}) \IMEXBDFexpcoef_k \Dminushatblock & \IMEXBDFprevcoef_k \Mhatblock - \frac{\dt}{\vareps \dx} \IMEXBDFexpcoef_k \UPhatblock
\end{bmatrix},
\end{align*}
and the right hand side depends on $\vareps$, $\scatlower$, $\dx$, $\dt$ only through $\frac{\vareps}{\scatlower \dx}$ and $\frac{\dt}{\vareps \dx}$.
Now, we define $\mathbf{S}_{\dx} = \textrm{diag}(S_{\dx}, \dots, S_{\dx}) \in \mathbb{R}^{N_f\times N_f}$,
and come to
\begin{align}
\mathbf{S}_{\dx}^{-1} \mathbf{G}^{(\IMEXBDForder,\pd)} \mathbf{S}_{\dx}
& =
\begin{bmatrix}
S_{\dx}^{-1} G_{L,\IMEXBDForder}^{-1} G_{R,n+\IMEXBDForder-1} S_{\dx} & \dots &S_{\dx}^{-1} G_{L,\IMEXBDForder}^{-1} G_{R,n+1} S_{\dx} & S_{\dx}^{-1} G_{L,\IMEXBDForder}^{-1} G_{R,n} S_{\dx} \\
\iden_{\pd(N_{\velangle}+1)} & && \\
& \ddots& &\\
& & \iden_{\pd(N_{\velangle}+1)}&
\end{bmatrix}\notag
\\
&
 =
\mathbf{\hat{G}}^{(\IMEXBDForder,\pd)}\Big(\frac{\vareps}{\scatlower \dx}, \frac{\dt}{\vareps \dx}; \xi\Big).
\end{align}

This shows that  $\mathbf{G}^{(\IMEXBDForder,\pd)}$ is similar to $\mathbf{\hat{G}}^{(\IMEXBDForder,\pd)}(\frac{\vareps}{\scatlower \dx}, \frac{\dt}{\vareps \dx}; \xi)$, and they share the same eigenvalues. Therefore  the eigenvalues of $\mathbf{G}^{(\IMEXBDForder,\pd)}$ depend on $\vareps$, $\scatlower$, $\dx$, $\dt$ only through $\vareps/({\scatlower \dx})$ and $\dt/({\vareps \dx})$.

%%%%%%%%%%%%%%%%%%%%%%%%%%%%%%%%%%%%%%%%%%
\section{Proof of Theorem \ref{thm:condH}}
\label{sec:proof:CondH}
%%%%%%%%%%%%%%%%%%%%%%%%%%%%%%%%%%%%%%%%%%

Throughout the proof
we will use the fact that for a symmetric matrix $A$ of size $n\times n$, its maximum and minimum eigenvalues satisfy 
\begin{equation}\lambda_{\max}(A)=\max_{{\by}\ne 0}\frac{{\by}^TA{\by}}{{\by}^T{\by}}, \quad \lambda_{\min}(A)=\min_{{\by}\ne 0}\frac{{\by}^TA{\by}}{{\by}^T{\by}},
\label{eq:condH:p:-1}\end{equation}
hence 
\begin{equation} \lambda_{\min}(A)|{\by}|^2\leq {\by}^TA{\by}\leq \lambda_{\max}(A)|{\by}| ^2,\quad \forall {\by}\in\mathbb{R}^n.
\label{eq:condH:p:0}
\end{equation}
Here $|{\by}|$ stands for the 2-norm of the vector ${\by}$. 

With $\mH$  being SPD,  we have
\begin{equation}
    \textrm{cond}_2(\mH)=\frac{\lambda_{\max}(\mH)}{\lambda_{\min}(\mH)}.
    \label{eq:condH:p:01}
\end{equation} To establish the estimate in \eqref{eq:condH}, it boils down to estimating
\begin{align}
    {\bf z}^T\mH {\bf z}&=(1+\dt c_\IMEXBDForder\absorp){\bf z}^T\M {\bf z}
    +\lvg \velangle^2\rvg \dt^2c_\IMEXBDForder^2(\Dminus{\bz})^T(\Theta^{(3)})^{-1} (\Dminus{\bz})\notag\\
    &=(1+\dt c_\IMEXBDForder\absorp){\bf z}^T\M {\bf z}
    +\frac{\lvg \velangle^2\rvg \dt^2c_\IMEXBDForder^2}{\vareps^2+\dt c_\IMEXBDForder(\vareps^2\absorp+\scat)}(\Dminus{\bz})^T(\M)^{-1} (\Dminus{\bz})
    \label{eq:condH:p:1}
\end{align}
with respect to $|{\bz}|^2$ for any  ${\bf z}\in \mathbb{R}^{N_{rx}}$. We  have used the facts that   $\Dplus=-(\Dminus)^T$ in the case of periodic boundary conditions, and the scattering and absorption cross sections $\scat$ and $\absorp$ are constant.

\medskip
\noindent{\bf Step 1.} Consider any given ${\bf z}\in \mathbb{R}^{N_{rx}}$ and define  $\psi(x)={\bf z}^T\mathbf{\pdspacebasis}$, we have $\|\psi\|^2=(\psi, \psi)={\bf z}^T\M {\bf z}.$ It is more helpful to connect back to the way how 
the basis $\mathbf{\pdspacebasis}$ is defined, and rewrite $${\bf z} =(z_{1,0},z_{1,1},\dots,z_{1,\pd-1},z_{2,0},z_{2,1},\dots,z_{2,\pd-1},\dots,z_{N_x,0},z_{N_x,1},\dots,z_{N_x,\pd-1})^T$$ and 
$\psi(x)=\sum_{m=1}^{N_x}\sum_{l=0}^{\pd-1}z_{m,l} \Phi_l^m(x).$ 
We recall one property of the basis  $\{\Phi^m_l\}_{m,l}$,
\begin{equation}
    (\Phi_l^m, \Phi_{l'}^{m'})=\delta_{l l'}\delta_{mm'}\frac{\dx_m}{2l+1},
\end{equation}
inherited from  that of the Legendre polynomials. With this,
\begin{equation}{\bf z}^T\M {\bf z}=\|\psi\|^2= \sum_{m=1}^{N_x}\sum_{l=0}^{\pd-1} |z_{m,l}|^2  \frac{\dx_m}{2l+1},
\end{equation}
and, for a uniform spatial mesh (i.e. $\dx_m=\dx,\forall m$), we further obtain
\begin{equation}
    \frac{\dx}{2\pd-1}|{\bf z}|^2 \leq {\bf z}^T\M {\bf z}=\|\psi\|^2\leq \dx |{\bf z}|^2, \quad \forall {\bf z}\in \mathbb{R}^{N_{rx}}, \quad \psi(x)={\bf z}^T\mathbf{\pdspacebasis}.
    \label{eq:condH:p:2}
\end{equation}
This implies $\lambda_{\min}(\M)\geq h/(2r-1)$. Moreover, 
\begin{equation}
(\Dminus{\bz})^T(\M)^{-1} (\Dminus{\bz})
\leq \lambda_{\max}(M^{-1})|\Dminus{\bz}|^2=
\big(\lambda_{\min}(M)\big)^{-1}|\Dminus{\bz}|^2\leq \frac{2r-1}{h}|\Dminus{\bz}|^2.
\label{eq:condH:p:3}
\end{equation}
In general, one only has
\begin{equation}
(\Dminus{\bz})^T(\M)^{-1} (\Dminus{\bz})\geq 0,
\label{eq:condH:p:4}
\end{equation}
with the lower bound zero attainable. This is due to that constant functions are in the kernel of the advective operator $\Dminusop$, and this renders a nontrivial null space of  $\Dminus$.

%%%%%%%%%
\medskip
\noindent{\bf Step 2.} 
Next we turn to the estimate for the advection matrix $\Dminus$. Based on its  definition, we have 
$(\Dminus{\bz})_i=\sum_{j=1}^{N_{rx}}(\Dminusop e_j, e_i)z_j=(\Dminusop \psi, e_i)$, and 
\begin{align}
    |\Dminus{\bz}|^2&=\sum_{i=1}^{N_{rx}}(\Dminusop \psi, e_i)^2
    =\sum_{m=1}^{N_x}\sum_{l=0}^{\pd-1}\big(\Dminusop \psi, \Phi_l^m\big)^2\notag\\
    &=\sum_{m=1}^{N_x}\sum_{l=0}^{\pd-1}\Big(\int_{I_m}(\Dminusop \psi(x)) \Phi_l^m(x)dx\Big)^2\notag\\
    &\leq \sum_{m=1}^{N_x}\sum_{l=0}^{\pd-1} \|\Dminusop \psi\|_{L^2(I_m)}^2\|\Phi_l^m\|^2
    =\sum_{m=1}^{N_x}\sum_{l=0}^{\pd-1} \|\Dminusop \psi\|_{L^2(I_m)}^2\;\frac{h_m}{2l+1}\notag\\
     &=\left(\sum_{l=0}^{\pd-1} \frac{1}{2l+1}\right)\dx\|\Dminusop \psi\|^2, \quad (\textrm{using}\; \dx_m=\dx). 
\end{align}
By applying an estimate in Theorem 5 of \cite{reyna2015operator},
$${\|\Dminusop \psi\|\leq 2(\sqrt{3}+1)\frac{r^2}{h}\|\psi\|,\quad \psi\in \pdspace^{r-1},}$$
one can arrive at
\begin{equation}
    |\Dminus{\bz}|^2\leq C_r\frac{\|\psi\|^2}{h}\leq C_r|{\bz}|^2,\quad \forall  {\bz}\in \mathbb{R}^{N_{rx}}.\label{eq:condH:p:5}
\end{equation}
The upper bound in \eqref{eq:condH:p:2} is used, and $C_r=4(\sqrt{3}+1)^2\left(\sum_{l=0}^{\pd-1} \frac{1}{2l+1}\right)r^4$.

%%%%%%%%%
\noindent{\bf Step 3.} We are ready to assemble \eqref{eq:condH:p:1} and all the estimates established so far, particularly  \eqref{eq:condH:p:2}-\eqref{eq:condH:p:4}, \eqref{eq:condH:p:5}, and have
\begin{align*}
    \lambda_{\max}(\mH)&=\max_{{\bz}\ne 0} \frac{{\bf z}^T\mH {\bf z}}{{\bf z}^T{\bf z}}
   \leq(1+\dt c_\IMEXBDForder\absorp)h+\frac{(2r-1)C_r \lvg \velangle^2\rvg \dt^2c_\IMEXBDForder^2 }{h\big(\vareps^2+\dt c_\IMEXBDForder(\vareps^2\absorp+\scat)\big)},\\
   \lambda_{\min}(\mH)&=\min_{{\bz}\ne 0} \frac{{\bf z}^T\mH {\bf z}}{{\bf z}^T{\bf z}}\geq \min_{{\bz}\ne 0} \frac{{\bf z}^T\M{\bf z}}{{\bf z}^T{\bf z}} \geq \frac{\dx}{2\pd-1}.
\end{align*} 
This, together with \eqref{eq:condH:p:01}, leads to  the estimate in \eqref{eq:condH}.

%%%%%%%%%%%%%%%%%%%%%%%%%%%%%%%%%%%%%
\section{Proof of Lemma \ref{lem:AP}}
\label{sec:proof:AP}
%%%%%%%%%%%%%%%%%%%%%%%%%%%%%%%%%%%%%
%%%%%% St1 %%%
{\bf Strategy 1:} 
From the proposed methods with Strategy 1 in \eqref{eq:fully:st1}, we get
\begin{align}
\micro_{h,q}^{n+\IMEXBDForder} = &\frac{\vareps^2}{\vareps^2 + \dt \IMEXBDFimpcoef_{\IMEXBDForder} \scat } \sum_{k=0}^{\IMEXBDForder-1} (\IMEXBDFprevcoef_k-\dt\absorp\IMEXBDFexpcoef_k) \micro_{h,q}^{n+k} - \frac{\dt\velangle_q}{\vareps^2 + \dt \IMEXBDFimpcoef_{\IMEXBDForder} \scat}   \sum_{k=0}^{\IMEXBDForder-1} \IMEXBDFexpcoef_k\Dminusop \macro_h^{n+k}
 \notag \\
&- \frac{\dt\vareps}{\vareps^2 + \dt \IMEXBDFimpcoef_{\IMEXBDForder} \scat }\sum_{k=0}^{\IMEXBDForder-1} \IMEXBDFexpcoef_k \Big(\Dupwindop(\micro_{h}^{n+k};\velangle_q)-\lvg \Dupwindop(\micro_{h}^{n+k};\velangle)\rvg_h \Big)
, \;\;q = 1, \dots, N_{\velangle}. 
\label{eq:ap:1}
\end{align}
If 
$ \micro_{h,q}^{n+k}=\mO(1)$, $\forall k=0,\dots, \IMEXBDForder-1$, $\forall q=1,\dots, N_{\velangle}$, then  \eqref{eq:ap:1} together with \eqref{eq:fully:st1.a}  implies \eqref{eq:ap:01}.
On the other hand, if $ \micro_{h,q}^{n+k}=\mO(1/\vareps)$ for some integer $k\in [0, \IMEXBDForder-1]$ and some integer $q\in [1, N_\velangle]$, then
\begin{equation}
\micro_{h,q}^{n+\IMEXBDForder} =- \frac{\velangle_q}{ \IMEXBDFimpcoef_{\IMEXBDForder} \scat}   \sum_{k=0}^{\IMEXBDForder-1} \IMEXBDFexpcoef_k\Dminusop \macro_h^{n+k} 
- \frac{1}{ \IMEXBDFimpcoef_{\IMEXBDForder} \scat } \underbrace{\vareps \sum_{k=0}^{\IMEXBDForder-1} \IMEXBDFexpcoef_k \Big(\Dupwindop(\micro_{h}^{n+k};\velangle_q)-\lvg \Dupwindop(\micro_{h}^{n+k};\velangle)\rvg_h \Big)}_{\mO(1)} +\mO(\vareps),  
\label{eq:ap:2}
\end{equation}
for $q = 1, \dots, N_{\velangle}$, hence we have \eqref{eq:ap:03} and \eqref{eq:ap:05} regarding $\micro_{h,q}^{n+\IMEXBDForder}$. The boundedness of 
$\macro_{h}^{n+\IMEXBDForder}$ follows from \eqref{eq:fully:st1.a}.

%%%%%% St2 %%%
\medskip
\noindent{\bf Strategy 2:} With 
$ \micro_{h,q}^{n+k}=\mO(1)$, $\forall k=0,\dots, \IMEXBDForder-1$, $\forall q=1,\dots, N_{\velangle}$, we get $\macro_h^{n+\IMEXBDForder}=\mO(1)$ from the proposed methods with Strategy 2 in \eqref{eq:fully:st2} and
 \begin{align}
\micro_{h,q}^{n+\IMEXBDForder} =& \frac{\vareps^2}{\vareps^2 + \dt \IMEXBDFimpcoef_{\IMEXBDForder} \scat } \sum_{k=0}^{\IMEXBDForder-1} (\IMEXBDFprevcoef_k-\dt\absorp\IMEXBDFexpcoef_k) \micro_{h,q}^{n+k} - \frac{\dt\IMEXBDFimpcoef_{\IMEXBDForder}\velangle_q }{\vareps^2 + \dt \IMEXBDFimpcoef_{\IMEXBDForder} \scat}   \Dminusop \macro_h^{n+s}
 \notag \\
&- \frac{\dt\vareps}{\vareps^2 + \dt \IMEXBDFimpcoef_{\IMEXBDForder} \scat }\sum_{k=0}^{\IMEXBDForder-1} \IMEXBDFexpcoef_k \Big(\Dupwindop(\micro_{h}^{n+k};\velangle_q)-\lvg \Dupwindop(\micro_{h}^{n+k};\velangle)\rvg_h \Big)\notag\\
=&- \frac{\velangle_q}{\scat}   \Dminusop \macro_h^{n+s}+\mO(\vareps) 
, \quad q = 1, \dots, N_{\velangle}. 
\label{eq:ap:4}
\end{align}

%%%%%% St3 %%%
\noindent{\bf Strategy 3:} From the proposed methods with Strategy 3 in \eqref{eq:fully:st3}, we have $\forall q = 1, \dots, N_{\velangle}$,
 \begin{align}
\micro_{h,q}^{n+\IMEXBDForder} =& \frac{\vareps^2}{\vareps^2 + \dt \IMEXBDFimpcoef_{\IMEXBDForder} (\scat+\vareps^2\absorp)} \sum_{k=0}^{\IMEXBDForder-1} \IMEXBDFprevcoef_k \micro_{h,q}^{n+k} - \frac{\dt\IMEXBDFimpcoef_{\IMEXBDForder}\velangle_q }{\vareps^2 + \dt \IMEXBDFimpcoef_{\IMEXBDForder} (\scat+\vareps^2\absorp)}   \Dminusop \macro_h^{n+s}
 \notag \\
&- \frac{\dt\vareps}{\vareps^2 + \dt \IMEXBDFimpcoef_{\IMEXBDForder} (\scat+\vareps^2\absorp)}\sum_{k=0}^{\IMEXBDForder-1} \IMEXBDFexpcoef_k \Big(\Dupwindop(\micro_{h}^{n+k};\velangle_q)-\lvg \Dupwindop(\micro_{h}^{n+k};\velangle)\rvg_h \Big).
\label{eq:ap:5}
\end{align}
Under the assumption on $\micro_{h,q}^{n+k}$ for the relevant $k$ and $q$, we have
\begin{equation}
    \micro_{h,q}^{n+\IMEXBDForder}=- \frac{\dt\IMEXBDFimpcoef_{\IMEXBDForder}\velangle_q }{\vareps^2 + \dt \IMEXBDFimpcoef_{\IMEXBDForder} (\scat+\vareps^2\absorp)}   \Dminusop \macro_h^{n+s}+\mO(1).
\end{equation}
This, combined with \eqref{eq:fully:st3.a}, leads to 

\begin{equation}
(1+\dt \IMEXBDFimpcoef_{\IMEXBDForder}\absorp) \macro_h^{n+\IMEXBDForder} =   \frac{(\dt \IMEXBDFimpcoef_{\IMEXBDForder})^2}{\vareps^2 + \dt \IMEXBDFimpcoef_{\IMEXBDForder} (\scat+\vareps^2\absorp)} \lvg \velangle^2\rvg  \Dplusop \Dminusop \macro_h^{n+\IMEXBDForder}  +\mO(1).
\end{equation}
By multiplying this equation with $\macro_h^{n+\IMEXBDForder}$  and integrating in space, and using $\Dplusop=-(\Dminusop)^T$ in \eqref{eq:adjoint}, one gets
\begin{equation}
(1+\dt \IMEXBDFimpcoef_{\IMEXBDForder}\absorp) \|\macro_h^{n+\IMEXBDForder}\|^2 +   \frac{(\dt \IMEXBDFimpcoef_{\IMEXBDForder})^2}{\vareps^2 + \dt \IMEXBDFimpcoef_{\IMEXBDForder} (\scat+\vareps^2\absorp)} \lvg \velangle^2\rvg  \|\Dminusop \macro_h^{n+\IMEXBDForder}\|^2  =\mO(1)\|\macro_h^{n+\IMEXBDForder}\|,
\end{equation}
implying $\|\macro_h^{n+\IMEXBDForder}\|=\mO(1)$. Now if 
$ \micro_{h,q}^{n+k}=\mO(1)$, $\forall k=0,\dots, \IMEXBDForder-1$, $\forall q=1,\dots, N_{\velangle}$, then
\eqref{eq:ap:5} will give \eqref{eq:ap:02},
while if $ \micro_{h,q}^{n+k}=\mO(1/\vareps)$ for some integer $k\in [0, \IMEXBDForder-1]$ and some integer $q\in [1, N_\velangle]$, then
\eqref{eq:ap:5} will give the following instead
\begin{equation}
\micro_{h,q}^{n+\IMEXBDForder} =- \frac{\velangle_q}{\scat}   \Dminusop \macro_h^{n+s}-\frac{1}{ \IMEXBDFimpcoef_{\IMEXBDForder} \scat } \underbrace{\vareps \sum_{k=0}^{\IMEXBDForder-1} \IMEXBDFexpcoef_k \Big(\Dupwindop(\micro_{h}^{n+k};\velangle_q)-\lvg \Dupwindop(\micro_{h}^{n+k};\velangle)\rvg_h \Big)}_{\mO(1)} +\mO(\vareps) 
\end{equation}
for any $q = 1, \dots, N_{\velangle}$, and this will lead to \eqref{eq:ap:04} and \eqref{eq:ap:05} regarding $\micro_{h,q}^{n+\IMEXBDForder}$.

%%%%%%%%%%%%%%%%%%%%%%%%%%%%%%%%%%%%%%%%%%%%%%%%%%%%%%%%%%%}
\end{appendices}

\end{document}